\expandafter\ifx\csname amssym.def\endcsname\relax \else \fi
%
\expandafter\edef\csname amssym.def\endcsname{%
       \catcode`\noexpand\@=\the\catcode`\@\space}
\catcode`\@=11
%

\def\undefine#1{\let#1\undefined}
\def\newsymbol#1#2#3#4#5{\let\next@\relax
 \ifnum#2=\@ne\let\next@\msafam@\else
 \ifnum#2=\tw@\let\next@\msbfam@\fi\fi
 \mathchardef#1="#3\next@#4#5}
\def\mathhexbox@#1#2#3{\relax
 \ifmmode\mathpalette{}{\m@th\mathchar"#1#2#3}%
 \else\leavevmode\hbox{$\m@th\mathchar"#1#2#3$}\fi}
\def\hexnumber@#1{\ifcase#1 0\or 1\or 2\or 3\or 4\or 5\or 6\or 7\or 8\or
 9\or A\or B\or C\or D\or E\or F\fi}

\font\tenmsa=msam10
\font\sevenmsa=msam7
\font\fivemsa=msam5
\newfam\msafam
\textfont\msafam=\tenmsa
\scriptfont\msafam=\sevenmsa
\scriptscriptfont\msafam=\fivemsa
\edef\msafam@{\hexnumber@\msafam}
\mathchardef\dabar@"0\msafam@39
\def\dashrightarrow{\mathrel{\dabar@\dabar@\mathchar"0\msafam@4B}}
\def\dashleftarrow{\mathrel{\mathchar"0\msafam@4C\dabar@\dabar@}}

\def\ulcorner{\delimiter"4\msafam@70\msafam@70 }
\def\urcorner{\delimiter"5\msafam@71\msafam@71 }
\def\llcorner{\delimiter"4\msafam@78\msafam@78 }
\def\lrcorner{\delimiter"5\msafam@79\msafam@79 }
\def\yen{{\mathhexbox@\msafam@55 }}
\def\checkmark{{\mathhexbox@\msafam@58 }}
\def\circledR{{\mathhexbox@\msafam@72 }}
\def\maltese{{\mathhexbox@\msafam@7A }}

\font\tenmsb=msbm10
\font\sevenmsb=msbm7
\font\fivemsb=msbm5
\newfam\msbfam
\textfont\msbfam=\tenmsb
\scriptfont\msbfam=\sevenmsb
\scriptscriptfont\msbfam=\fivemsb
\edef\msbfam@{\hexnumber@\msbfam}
\def\Bbb#1{{\fam\msbfam\relax#1}}
\def\widehat#1{\setbox\z@\hbox{$\m@th#1$}%
 \ifdim\wd\z@>\tw@ em\mathaccent"0\msbfam@5B{#1}%
 \else\mathaccent"0362{#1}\fi}
\def\widetilde#1{\setbox\z@\hbox{$\m@th#1$}%
 \ifdim\wd\z@>\tw@ em\mathaccent"0\msbfam@5D{#1}%
 \else\mathaccent"0365{#1}\fi}
\font\teneufm=eufm10
\font\seveneufm=eufm7
\font\fiveeufm=eufm5
\newfam\eufmfam
\textfont\eufmfam=\teneufm
\scriptfont\eufmfam=\seveneufm
\scriptscriptfont\eufmfam=\fiveeufm

\csname amssym.def\endcsname

\expandafter\ifx\csname pre amssym.tex at\endcsname\relax \else  \fi
\expandafter\chardef\csname pre amssym.tex at\endcsname=\the\catcode`\@
\catcode`\@=11
\newsymbol\boxdot 1200
\newsymbol\boxplus 1201
\newsymbol\boxtimes 1202
\newsymbol\square 1003
\newsymbol\blacksquare 1004
\newsymbol\centerdot 1205
\newsymbol\lozenge 1006
\newsymbol\blacklozenge 1007
\newsymbol\circlearrowright 1308
\newsymbol\circlearrowleft 1309
\undefine\rightleftharpoons
\newsymbol\rightleftharpoons 130A
\newsymbol\leftrightharpoons 130B
\newsymbol\boxminus 120C
\newsymbol\Vdash 130D
\newsymbol\Vvdash 130E
\newsymbol\vDash 130F
\newsymbol\twoheadrightarrow 1310
\newsymbol\twoheadleftarrow 1311
\newsymbol\leftleftarrows 1312
\newsymbol\rightrightarrows 1313
\newsymbol\upuparrows 1314
\newsymbol\downdownarrows 1315
\newsymbol\upharpoonright 1316
 \let\restriction\upharpoonright
\newsymbol\downharpoonright 1317
\newsymbol\upharpoonleft 1318
\newsymbol\downharpoonleft 1319
\newsymbol\rightarrowtail 131A
\newsymbol\leftarrowtail 131B
\newsymbol\leftrightarrows 131C
\newsymbol\rightleftarrows 131D
\newsymbol\Lsh 131E
\newsymbol\Rsh 131F
\newsymbol\rightsquigarrow 1320
\newsymbol\leftrightsquigarrow 1321
\newsymbol\looparrowleft 1322
\newsymbol\looparrowright 1323
\newsymbol\circeq 1324
\newsymbol\succsim 1325
\newsymbol\gtrsim 1326
\newsymbol\gtrapprox 1327
\newsymbol\multimap 1328
\newsymbol\therefore 1329
\newsymbol\because 132A
\newsymbol\doteqdot 132B
 
\newsymbol\triangleq 132C
\newsymbol\precsim 132D
\newsymbol\lesssim 132E
\newsymbol\lessapprox 132F
\newsymbol\eqslantless 1330
\newsymbol\eqslantgtr 1331
\newsymbol\curlyeqprec 1332
\newsymbol\curlyeqsucc 1333
\newsymbol\preccurlyeq 1334
\newsymbol\leqq 1335
\newsymbol\leqslant 1336
\newsymbol\lessgtr 1337
\newsymbol\backprime 1038
\newsymbol\risingdotseq 133A
\newsymbol\fallingdotseq 133B
\newsymbol\succcurlyeq 133C
\newsymbol\geqq 133D
\newsymbol\geqslant 133E
\newsymbol\gtrless 133F
\newsymbol\sqsubset 1340
\newsymbol\sqsupset 1341
\newsymbol\vartriangleright 1342
\newsymbol\vartriangleleft 1343
\newsymbol\trianglerighteq 1344
\newsymbol\trianglelefteq 1345
\newsymbol\bigstar 1046
\newsymbol\between 1347
\newsymbol\blacktriangledown 1048
\newsymbol\blacktriangleright 1349
\newsymbol\blacktriangleleft 134A
\newsymbol\vartriangle 134D
\newsymbol\blacktriangle 104E
\newsymbol\triangledown 104F
\newsymbol\eqcirc 1350
\newsymbol\lesseqgtr 1351
\newsymbol\gtreqless 1352
\newsymbol\lesseqqgtr 1353
\newsymbol\gtreqqless 1354
\newsymbol\Rrightarrow 1356
\newsymbol\Lleftarrow 1357
\newsymbol\veebar 1259
\newsymbol\barwedge 125A
\newsymbol\doublebarwedge 125B
\undefine\angle
\newsymbol\angle 105C
\newsymbol\measuredangle 105D
\newsymbol\sphericalangle 105E
\newsymbol\varpropto 135F
\newsymbol\smallsmile 1360
\newsymbol\smallfrown 1361
\newsymbol\Subset 1362
\newsymbol\Supset 1363
\newsymbol\Cup 1264
 
\newsymbol\Cap 1265
 
\newsymbol\curlywedge 1266
\newsymbol\curlyvee 1267
\newsymbol\leftthreetimes 1268
\newsymbol\rightthreetimes 1269
\newsymbol\subseteqq 136A
\newsymbol\supseteqq 136B
\newsymbol\bumpeq 136C
\newsymbol\Bumpeq 136D
\newsymbol\lll 136E
 
\newsymbol\ggg 136F
 
\newsymbol\circledS 1073
\newsymbol\pitchfork 1374
\newsymbol\dotplus 1275
\newsymbol\backsim 1376
\newsymbol\backsimeq 1377
\newsymbol\complement 107B
\newsymbol\intercal 127C
\newsymbol\circledcirc 127D
\newsymbol\circledast 127E
\newsymbol\circleddash 127F
\newsymbol\lvertneqq 2300
\newsymbol\gvertneqq 2301
\newsymbol\nleq 2302
\newsymbol\ngeq 2303
\newsymbol\nless 2304
\newsymbol\ngtr 2305
\newsymbol\nprec 2306
\newsymbol\nsucc 2307
\newsymbol\lneqq 2308
\newsymbol\gneqq 2309
\newsymbol\nleqslant 230A
\newsymbol\ngeqslant 230B
\newsymbol\lneq 230C
\newsymbol\gneq 230D
\newsymbol\npreceq 230E
\newsymbol\nsucceq 230F
\newsymbol\precnsim 2310
\newsymbol\succnsim 2311
\newsymbol\lnsim 2312
\newsymbol\gnsim 2313
\newsymbol\nleqq 2314
\newsymbol\ngeqq 2315
\newsymbol\precneqq 2316
\newsymbol\succneqq 2317
\newsymbol\precnapprox 2318
\newsymbol\succnapprox 2319
\newsymbol\lnapprox 231A
\newsymbol\gnapprox 231B
\newsymbol\nsim 231C
\newsymbol\ncong 231D
\newsymbol\diagup 231E
\newsymbol\diagdown 231F
\newsymbol\varsubsetneq 2320
\newsymbol\varsupsetneq 2321
\newsymbol\nsubseteqq 2322
\newsymbol\nsupseteqq 2323
\newsymbol\subsetneqq 2324
\newsymbol\supsetneqq 2325
\newsymbol\varsubsetneqq 2326
\newsymbol\varsupsetneqq 2327
\newsymbol\subsetneq 2328
\newsymbol\supsetneq 2329
\newsymbol\nsubseteq 232A
\newsymbol\nsupseteq 232B
\newsymbol\nparallel 232C
\newsymbol\nmid 232D
\newsymbol\nshortmid 232E
\newsymbol\nshortparallel 232F
\newsymbol\nvdash 2330
\newsymbol\nVdash 2331
\newsymbol\nvDash 2332
\newsymbol\nVDash 2333
\newsymbol\ntrianglerighteq 2334
\newsymbol\ntrianglelefteq 2335
\newsymbol\ntriangleleft 2336
\newsymbol\ntriangleright 2337
\newsymbol\nleftarrow 2338
\newsymbol\nrightarrow 2339
\newsymbol\nLeftarrow 233A
\newsymbol\nRightarrow 233B
\newsymbol\nLeftrightarrow 233C
\newsymbol\nleftrightarrow 233D
\newsymbol\divideontimes 223E
\newsymbol\varnothing 203F
\newsymbol\nexists 2040
\newsymbol\Finv 2060
\newsymbol\Game 2061
\newsymbol\mho 2066
\newsymbol\eth 2067
\newsymbol\eqsim 2368
\newsymbol\beth 2069
\newsymbol\gimel 206A
\newsymbol\daleth 206B
\newsymbol\lessdot 236C
\newsymbol\gtrdot 236D
\newsymbol\ltimes 226E
\newsymbol\rtimes 226F
\newsymbol\shortmid 2370
\newsymbol\shortparallel 2371
\newsymbol\smallsetminus 2272
\newsymbol\thicksim 2373
\newsymbol\thickapprox 2374
\newsymbol\approxeq 2375
\newsymbol\succapprox 2376
\newsymbol\precapprox 2377
\newsymbol\curvearrowleft 2378
\newsymbol\curvearrowright 2379
\newsymbol\digamma 207A
\newsymbol\varkappa 207B
\newsymbol\Bbbk 207C
\newsymbol\hslash 207D
\undefine\hbar
\newsymbol\hbar 207E
\newsymbol\backepsilon 237F
\catcode`\@=\csname pre amssym.tex at\endcsname

%
%
\font\fivebi=cmmib5
\font\fivebsy=cmbsy5
\font\sixrm=cmr6
\font\sixi=cmmi6
\font\sixbf=cmbx6
\font\sixsy=cmsy6
\font\sixmsa=msam5 at 6pt
\font\sixmsb=msbm5 at 6pt
\font\sevenbi=cmmib7
\font\sevenbsy=cmbsy7
\font\eightrm=cmr8
\font\eightsl=cmsl8
\font\eightit=cmti8
\font\eighti=cmmi8
\font\eightbf=cmbx8
\font\eightsy=cmsy8
\font\eightmsa=msam7 at 8pt
\font\eightmsb=msbm7 at 8pt
\font\ninerm=cmr9
\font\ninesl=cmsl9
\font\nineit=cmti9
\font\ninei=cmmi9
\font\ninebf=cmbx9
\font\ninebi=cmmib10 scaled 900
\font\ninesy=cmsy9
\font\ninebsy=cmbsy10 scaled 900
\font\ninemsa=msam10 at 9pt
\font\ninemsb=msbm10 at 9pt
\font\tenbit=cmbxti10
\font\tenbsl=cmbxsl10
\font\tenbi=cmmib10
\font\tenbsy=cmbsy10
\font\twelvebf=cmbx12
\font\twelvebi=cmmib10 scaled 1200
\font\twelvebsy=cmbsy10 at 12pt

\let\sc=\sevenrm            
\def\eightpoint{%
     \def\rm{\fam0\eightrm}
     \textfont0=\eightrm \scriptfont0=\sixrm \scriptscriptfont0=\fiverm
     \textfont1=\eighti \scriptfont1=\sixi \scriptscriptfont1=\fivei
     \textfont2=\eightsy \scriptfont2=\sixsy \scriptscriptfont2=\fivesy
     \textfont3=\tenex \scriptfont3=\tenex \scriptscriptfont3=\tenex
     \textfont\itfam=\eightit \def\it{\fam\itfam\eightit}%
     \textfont\slfam=\eightsl \def\sl{\fam\slfam\eightsl}%
     \textfont\bffam=\eightbf \scriptfont\bffam=\sixbf
     \scriptscriptfont\bffam=\fivebf \def\bf{\fam\bffam\eightbf}%
     \textfont\msbfam=\eightmsb \textfont\msafam=\eightmsa
     \scriptfont\msafam=\sixmsa \scriptfont\msbfam=\sixmsb
     \scriptscriptfont\msafam=\fivemsa \scriptscriptfont\msbfam=\fivemsb
      \skewchar\eighti='177 \skewchar\sixi='177
      \skewchar\eightsy='60 \skewchar\sixsy='60
     \normalbaselineskip=10pt
     \setbox\strutbox=\hbox{\vrule height7pt depth3pt width0pt}
     \let\sc=\sixrm \let\big=\eightbig \normalbaselines\rm}
\def\ninepoint{%
     \def\rm{\fam0\ninerm}
     \textfont0=\ninerm \scriptfont0=\sixrm \scriptscriptfont0=\fiverm
     \textfont1=\ninei \scriptfont1=\sixi \scriptscriptfont1=\fivei
     \textfont2=\ninesy \scriptfont2=\sixsy \scriptscriptfont2=\fivesy
     \textfont3=\tenex \scriptfont3=\tenex \scriptscriptfont3=\tenex
     \textfont\itfam=\nineit \def\it{\fam\itfam\nineit}%
     \textfont\slfam=\ninesl \def\sl{\fam\slfam\ninesl}%
     \textfont\bffam=\ninebf \scriptfont\bffam=\sixbf
     \scriptscriptfont\bffam=\fivebf \def\bf{\fam\bffam\ninebf}%
     \textfont\msbfam=\ninemsb \textfont\msafam=\ninemsa
     \scriptfont\msafam=\sixmsa \scriptfont\msbfam=\sixmsb
     \scriptscriptfont\msafam=\fivemsa \scriptscriptfont\msbfam=\fivemsb
      \skewchar\ninei='177 \skewchar\sixi='177
      \skewchar\ninesy='60 \skewchar\sixsy='60
     \normalbaselineskip=11pt
     \setbox\strutbox=\hbox{\vrule height8pt depth3pt width0pt}%
     \let\sc=\sevenrm \let\big=\ninebig \normalbaselines\rm}
\def\tenpoint{%
     \def\rm{\fam0\tenrm}
     \textfont0=\tenrm \scriptfont0=\sevenrm \scriptscriptfont0=\fiverm
     \textfont1=\teni \scriptfont1=\seveni \scriptscriptfont1=\fivei
     \textfont2=\tensy \scriptfont2=\sevensy \scriptscriptfont2=\fivesy
     \textfont3=\tenex \scriptfont3=\tenex \scriptscriptfont3=\tenex
     \textfont\itfam=\tenit \def\it{\fam\itfam\tenit}%
     \textfont\slfam=\tensl \def\sl{\fam\slfam\tensl}%
     \textfont\bffam=\tenbf \scriptfont\bffam=\sevenbf
     \scriptscriptfont\bffam=\fivebf \def\bf{\fam\bffam\tenbf}%
     \textfont\msbfam=\tenmsb \textfont\msafam=\tenmsa
     \scriptfont\msafam=\sevenmsa \scriptfont\msbfam=\sevenmsb
     \scriptscriptfont\msafam=\fivemsa \scriptscriptfont\msbfam=\fivemsb
     \normalbaselineskip=12pt
     \let\sc=\sevenrm \let\big=\tenbig \normalbaselines\rm}
\catcode`@=11        
\def\tenbig#1{{\hbox{$\left#1\vbox to8.5pt{}\right.\n@space$}}}
\def\ninebig#1{{\hbox{$\textfont0=\tenrm\textfont2=\tensy
      \left#1\vbox to7.25pt{}\right.\n@space$}}}
\def\eightbig#1{{\hbox{$\textfont0=\ninerm\textfont2=\ninesy
      \left#1\vbox to6.5pt{}\right.\n@space$}}}
\catcode`@=12        
\def\bold{%
     \textfont0=\tenbf \scriptfont0=\sevenbf \scriptscriptfont0=\fivebf
     \textfont1=\tenbi \scriptfont1=\sevenbi \scriptscriptfont1=\fivebi
     \textfont2=\tenbsy \scriptfont2=\sevenbsy \scriptscriptfont2=\fivebsy
       \textfont\itfam=\tenbit \def\it{\fam\itfam\tenbit}%
       \textfont\slfam=\tenbsl \def\sl{\fam\slfam\tenbsl}%
       \textfont\bffam=\tenbf \scriptfont\bffam=\sevenbf
       \textfont\msbfam=\tenmsb \textfont\msafam=\tenmsa
   \fam0\tenbf}
\def\bigbold{%
     \textfont0=\twelvebf \scriptfont0=\ninebf \scriptscriptfont0=\sevenbf
     \textfont1=\twelvebi \scriptfont1=\ninebi \scriptscriptfont1=\sevenbi
     \textfont2=\twelvebsy \scriptfont2=\ninebsy \scriptscriptfont2=\sevenbsy
     \fam0\twelvebf}
%
%
%
%
\let\plainitem=\item
\let\plainitemitem=\itemitem
%

%
%
\def\openface{\Bbb}                
\def\N{{\openface N}}              
\def\Z{{\openface Z}}

\def\R{{\openface R}}
\def\C{{\openface C}}
%
%
%
\def\g{\hskip.17em\relax}               
\def\th{\thinspace}                     
\def\nl{\hfil\break}
\newskip\Bigskipamount
   \Bigskipamount=2\baselineskip plus.5\baselineskip minus.3\baselineskip
\def\Bigbreak{\removelastskip\vskip0pt plus .1\vsize\penalty-1000
              \vskip0pt plus-.1\vsize\vskip\Bigskipamount}
\def\Nobreak$$#1$${\postdisplaypenalty=10000$$#1$$\postdisplaypenalty=0}
%
%
%
\let\HHHHH=\H \def\Erdos{Erd\HHHHH os}    
\let\vvvvv=\v

\let\doublebar=\| 
\def\|{\!\!\restriction\!\!}
\def\bs{\backslash}

  \let\sub=\sube

\def\supe{\supseteq}

\def\sm{\smallsetminus}
\def\es{\emptyset}

\def\H{{\rm H}}                          
\def\wrt{with respect to}

\def\:{\colon}
\def\minor{\preccurlyeq} 
\def\Minor{\succcurlyeq}

\def\slt{\mathrel{\hbox{$\minor$\kern-.6em\lower.33ex\hbox{${}_s\;$}}}}
\def\sgt{\mathrel{\mathchoice                        
   {\hbox{$\Minor$\kern-.5em\lower.3ex\hbox{${}_s$}}}
   {\hbox{$\Minor$\kern-.5em\lower.3ex\hbox{${}_s$}}}
   {\hbox{$\scriptstyle\Minor\kern-.43em\lower.28ex\hbox{$\scriptstyle{}_s$}$}}
 {\hbox{$\scriptstyle\Minor\kern-.43em\lower.28ex\hbox{$\scriptstyle{}_s$}$}} }}    

\def\ucl(#1){\lfloor #1 \rfloor}
\def\dcl(#1){\lceil #1 \rceil}
\def\interior{\mathaccent"7017\relax}
%
%
\def\specrel#1#2{\mathrel{\mathop{\kern0pt #1}\limits_{#2}}}
\def\Specrel#1#2{\mathrel{\mathop{\kern0pt #1}\limits^{#2}}}
%
%
\def\alignspecrel#1#2{\mathrel{\mathop{\kern0pt #1}\limits_{\hbox
   to0pt{\hss$\scriptstyle#2$\hss}}}}
\def\alignSpecrel#1#2{\mathrel{\mathop{\kern0pt #1}\limits^{\hbox
   to0pt{\hss$\scriptstyle#2$\hss}}}}
\def\invlim{\specrel\lim{\raise 2pt\hbox{$\longleftarrow$}}}
\def\proof{\removelastskip\penalty55\medskip\noindent{\bf Proof. }}
\def\noproof{{\unskip\nobreak\hfill\penalty50\hskip2em\hbox{}\nobreak\hfill%
       $\square$\parfillskip=0pt\finalhyphendemerits=0\par}\goodbreak}
\def\endproof{\noproof\bigskip}
\newcount\refno
\def\ref#1#2\par{{\plainitem{[??]}#2\smallskip}}
\newtoks\thingtowrite 
\long\def\writerefnumber#1{%
    \thingtowrite={#1}%
    \immediate\write\refnumbersfile{\the\thingtowrite}%
    }
\newwrite\refnumbersfile
\def\makerefnumbers{\immediate\openout\refnumbersfile=RefNumbers%
  \refno=0 \writerefnumber{\refno=0}
  \def\ref##1##2\par{\global\advance\refno by 1
    \writerefnumber{\global\advance\refno by 1 \newcounter##1 ##1=\the\refno}%
    \plainitem{[\the\refno]}##2\smallskip
    }%
  }
\def\autorefnumbers{\refno=0
  \def\ref##1##2\par{\advance\refno by 1\plainitem{[\the\refno]}##2\smallskip}
  }
\def\userefnumbers{\refno=0
  \def\ref##1##2\par{\advance\refno by 1\plainitem{[\the##1]}##2\smallskip}
  }
%
%
\def\proclaimwithname #1. (#2) #3\par{{\bigbreak
  \clubpenalty=10000\noindent{\bf#1.\enspace}(#2)\nl
  {\sl #3}\par\bigbreak}}
\def\proposition (#1) #2\par{{\setbox0\hbox{(#1)\enspace}\bigbreak
   \sl\hangindent\the\wd0 \noindent\hskip\the\wd0
   \llap{\box0}\ignorespaces#2\par\bigbreak}}
\def\subsection #1\par{\vskip 3\medskipamount minus \smallskipamount\leftline{\bold #1}
        \penalty10000\smallskip\noindent}
      \def\section #1\par{\Bigbreak\centerline{\bf #1} 
              \penalty10000\bigskip\noindent}
%
\def\beginpsection #1\par{\Bigbreak\centerline{\bold #1}
        \penalty10000\bigskip}
\def\psubsection #1\par{\bigbreak\leftline{\bold #1}\penalty10000\bigskip}
%

%
%
\def\pitem#1{\smallskip\advance\parindent by 3mm
             \plainitem{\rm(#1)}\advance\parindent by-3mm}
\def\pitemitem#1{\smallskip\advance\parindent by 3mm
             \plainitemitem{\rm(#1)}\advance\parindent by-3mm}
\def\varitemitem#1{{\setbox0\hbox{\hskip\parindent#1\enskip}
           \smallbreak\hangindent\the\wd0 \noindent\hskip\the\wd0
           \llap{#1\enskip}\ignorespaces}}
%
%
\newdimen\newparindent
\def\iitem#1#2\par{\newparindent=\parindent \advance\newparindent by 3mm
           \smallbreak \hangindent\newparindent \noindent\hskip\newparindent
           \llap{{\rm #1}\enspace}\ignorespaces#2\par\smallbreak}
\def\iitemitem#1#2\par{\newparindent=\parindent \advance\newparindent by 3mm
           \smallbreak \hangindent2\newparindent \noindent\hskip2\newparindent
           \llap{{\rm #1}\enspace}\ignorespaces#2\par\smallbreak}
\def\varitem#1#2\par{{\setbox0\hbox{{\rm #1}\enspace}
           \smallbreak \hangindent\the\wd0 \noindent\hskip\the\wd0
           \llap{{\rm #1}\enspace}\ignorespaces#2\par\smallbreak}}
\def\enditem{\par}
%
\def\Textindent#1{\par \advance\parindent by 3mm
                  \textindent{{\rm #1}} \advance\parindent by -3mm}
\def\indentedline#1{\advance\hsize by -\parindent \line{#1}
                   \advance\hsize by \parindent}
\def\iindentedline#1{\advance\parindent by 3mm
                     \advance\hsize by -\parindent
                     \line{#1}
                     \advance\hsize by \parindent
                     \advance\parindent by -3mm}
\newdimen\margin   
\def\textdisplay#1&#2&#3$${\margin=\hsize
          \setbox1=\hbox{$\displaystyle#1\quad$}%
          \setbox2=\hbox{\quad#2\qquad$#3$}%
                     \advance\margin by-\wd1
                     \divide\margin by 2
   \ifdim\wd2 < \margin
      \box1\rlap{\quad#2}\eqno#3$$%
   \else
      \line{\qquad\hfil \box1\quad #2 \qquad $#3$}$$%
   \fi}
%
\def\ltextdisplay#1&#2&#3$${\margin=\hsize
           \setbox2=\hbox{$\displaystyle#2\quad$}
           \setbox3=\hbox{\quad#3\qquad}
                     \advance\margin by-\wd2
                     \divide\margin by 2
   \ifdim\wd3 < \margin
      \line{$#1$\hfil\box2\hbox to \margin{\box3\hfil}}$$%
   \else
      \line{$#1$\qquad\hfil\box2\quad #3\qquad} $$%
   \fi}
%
\def\textno#1&#2\par{%
   \margin=\hsize
   \advance\margin by -4\parindent
          \setbox1=\hbox{\sl#1}%
   \ifdim\wd1 < \margin
      $$\box1\eqno#2$$\endgraf%
   \else
      \bigbreak
      \line{\indent$\vcenter{\advance\hsize by -3\parindent
      \sl\noindent#1}\hfil#2$}%
      \bigbreak
   \fi}
%
\def\textlno#1&#2\par{%
   \margin=\hsize
   \advance\margin by -4\parindent
          \setbox1=\hbox{\sl#1}%
   \ifdim\wd1 < \margin
      $$\box1\leqno#2$$%
   \else
      \bigbreak
      \line{$#2\hfil\vcenter{\advance\hsize by -3\parindent
          \sl\noindent#1}\hskip\parindent$}%
      \bigbreak
   \fi}
%
%
%
\newcount\commentno
\def\COMMENT#1{$^{<\the\commentno>}$%
     \vadjust{\vbox to 0pt{\vss\vskip-8pt\rightline{%
     \rlap{\hbox{\hskip7mm \vbox{\pretolerance=-1
     \doublehyphendemerits=0 \finalhyphendemerits=0
     \hsize40mm\tolerance=10000\eightpoint
     \lineskip=0pt\lineskiplimit=0pt
     \rightskip=0pt plus16mm\baselineskip8pt\noindent
     \hskip0pt       
     {$\langle$\the\commentno. #1$\rangle$}\endgraf}}}}\vss}}%
     \global\advance\commentno by1}%
\def\writecommentsasfootnotes{%
 \def\COMMENT{\global\advance\commentno by1\footnote{$^{<\the\commentno>}$}}%
 }
\def\nocomments{\def\COMMENT##1{}}
%
%
\def\?#1{\vadjust{\vbox to 0pt{\vss\vskip-8pt\leftline{%
     \llap{\hbox{\vbox{\pretolerance=-1
     \doublehyphendemerits=0\finalhyphendemerits=0
     \hsize16truemm\tolerance=10000\eightpoint
     \lineskip=0pt\lineskiplimit=0pt
     \rightskip=0pt plus16truemm\baselineskip8pt\noindent
     \hskip0pt        
     #1\endgraf}\hskip7truemm}}}\vss}}}
\def\d{}
%
%
%
%
\def\ds#1{}
%
%
\long\def\indexwrite#1{%
    \thingtowrite={#1}%
    \immediate\write\index{\the\thingtowrite}%
    }
%
%
\newwrite\index
\def\makeindex{\immediate\openout\index=index%
   \immediate\write\index{\catcode`@=11}%
   \def\d##1 {\ifmmode
     \write\index{$##1$, }%
     \write\index{\the\count0}\write\index{}
   \else
     \write\index{{##1}, }%
     \write\index{\the\count0}\write\index{}
   \fi {##1} }
      \def\ds##1{\ifmmode
     \write\index{$##1$, }%
     \write\index{\the\count0}\write\index{}
   \else
     \write\index{##1, }%
     \write\index{\the\count0}\write\index{}
   \fi}}
%
\newdimen\gap
\gap=3truemm
\newdimen\hackwidth
\hackwidth=15truemm
\def\disablems{\def\mos##1{\strut}}
\def\m#1{\ds{#1}\mo{#1}}
\def\mo#1{\ifmmode {#1}\else {\it#1}\fi\mos{#1}}
\def\mos#1{\ifmmode
     \strut\vadjust{\vbox to 0pt{\vss\kern-11pt\leftline{%
     \llap{\hbox{\vbox{\pretolerance=-1
     \doublehyphendemerits=0\finalhyphendemerits=0
     \hsize\hackwidth\tolerance=10000\eightpoint
     \lineskip=0pt\lineskiplimit=0pt
     \rightskip=0pt plus\hsize\baselineskip8pt\noindent
     $#1$\strut\endgraf}\hskip\gap }}}\vss}}%
   \else
     \strut\vadjust{\vbox to 0pt{\vss\kern-11pt\leftline{%
     \llap{\hbox{\vbox{\pretolerance=-1
     \doublehyphendemerits=0\finalhyphendemerits=0
     \hsize\hackwidth\tolerance=10000\eightpoint
     \lineskip=0pt\lineskiplimit=0pt
     \rightskip=0pt plus\hsize\baselineskip8pt\noindent
     \hskip0pt    
     {\sl#1}\strut\endgraf}\hskip\gap }}}\vss}}%
   \fi}%
\newcount\remarkno
\def\REMARK#1{{\footnote{${}^{\the\remarkno}$}{{#1}}%
   \global\advance\remarkno by1}}
\def\noremarks{\def\REMARK##1{}}
%
%
%
%
\def\picture #1 by #2 (#3){
  \vbox to #2{
          \vfill
          \special{picture #3}
          \hrule width #1 height 0pt depth 0pt
           }}
\newdimen\topfiguremargin
   \topfiguremargin=0pt                                  
\newdimen\bottomfiguremargin
   \bottomfiguremargin=\medskipamount                    
\newdimen\normalpictureheight
\normalpictureheight=40mm
\def\Fig.#1 (#2by#3; heightfactor:#4; caption:#5) {{%
   \dimen2=\normalpictureheight
   \dimen0=#2                          
      \divide\dimen2 by 1000
      \multiply\dimen2 by#4              
   \count2=\dimen2                  
      \dimen1=#3                             
   \count1=\dimen1
   \divide\count1 by 1000
   \divide\count2 by \count1          
   \divide\dimen0 by 1000
   \multiply\dimen0 by \count2      
         \dimen1=\hsize
         \advance\dimen1 by -\dimen0
         \divide\dimen1 by 2               
   \midinsert
   \vbox to \topfiguremargin{\vfil}
   \noindent\hskip\dimen1
   \picture\dimen0 by \dimen2  (Fig.#1 scaled \the\count2)%
   \vskip\bottomfiguremargin                     
      \ninepoint
      \parindent=.1\hsize\narrower\narrower
      \setbox0\hbox{#5}
      \ifdim\wd0 < .6\hsize
           \centerline{F{\sc IGURE} #1.\hskip1em#5}
       \else
           \plainitem{F{\sc IGURE} #1. }#5\par
       \fi
   \vskip0pt\endinsert}}
%
\def\textpicture #1(#2by#3; #4width#5lower#6){{%
      \dimen0=#5\count2=\dimen0                    
      \dimen0=#2\count1=\dimen0                    
   \divide\count1 by 1000
   \divide\count2 by \count1                 
   \hbox{\vrule #4width0pt\vbox to 0pt{\vss\vskip#6%
      \special{picture #1 scaled \the\count2}\hrule width#5 height0pt\vss}}}}
%
%
%
%
%
\def\figure #1. #2 (#3; #4) {{%
   \def\bigskip{\par\ifdim\lastskip<\bigskipamount\removelastskip
      \vskip\bigskipamount\fi}
   \midinsert\vskip\topfiguremargin
   \dimen0=\normalpictureheight
      \divide\dimen0 by 1000
      \multiply\dimen0 by#4        
   \centerline{\epsfbox{#3.eps}}
   \vskip\bottomfiguremargin                     
      \ninepoint
      \parindent=.1\hsize\narrower\narrower
      \setbox0\hbox{#2}
      \ifdim\wd0 < .6\hsize
           \centerline{F{\sc IGURE} #1.\hskip1em#2}
       \else
           \plainitem{F{\sc IGURE} #1. }#2\par
       \fi
  \endinsert}}
%
%

%
%
\def\Abh#1 {{\sl Abh.\g Math.\g Sem.\g Univ.\g Hamburg\penalty100\ \bf#1\ }}
\def\AMASH#1 {{\sl Acta Math.\g Acad.\g Sci.\g Hung.\penalty100\ \bf#1\ }}
\def\Advances#1 {{\sl Adv.\g Math.\penalty100\ \bf#1\ }}
\def\Annals#1 {{\sl Ann.\g Math.\penalty100\ \bf#1\ }}
\def\AnnComb#1 {{\sl Ann.\g Comb.\penalty100\ \bf#1\ }}
\def\AMM#1 {{\sl Amer.\g Math.\g Monthly\penalty100\ \bf#1\ }}
\def\Archiv#1 {{\sl Arch.}\g {\sl Math.\penalty100\ \bf#1\ }}
\def\ArsComb#1 {{\sl Ars Comb.\penalty100\ \bf#1\ }}
\def\CJM#1 {{\sl Can.\g J.\th Math.\penalty100\ \bf#1\ }}
\def\Comb#1 {{\sl Com\-bi\-na\-to\-ri\-ca\penalty100\ \bf#1\ }}
\def\CPC#1 {{\sl Comb.\g Probab.\g Comput.\penalty100\ \bf#1\ }}
\def\Crelle#1 {{\sl J.}\th {\sl Reine Angew.}\g
    {\sl Math.\penalty100\ \bf#1\ }}
\def\DM#1 {{\sl Discrete Math.\penalty100\ \bf#1\ }}
\def\DAM#1 {{\sl Discrete Appl.\g Math.\penalty100\ \bf#1\ }}
\def\EJC#1 {{\sl Eur.}\g{\sl J.}\g{\sl Comb.\penalty100\ \bf#1\ }}
\def\EJ#1 {{\sl Electronic.}\g{\sl J.}\g{\sl Comb.\penalty100\ \bf#1\ }}
\def\GC#1 {{\sl Graphs Comb.\penalty100\ \bf#1\ }}
\def\IJ#1 {{\sl Isr.\g J.\th Math.\penalty100\ \bf#1\ }}
\def\Inv#1 {{\sl In\-vent.\g math.\penalty100\ \bf#1\ }}
\def\JAlg#1 {{\sl J.}\th {\sl Algorithms\penalty100\ \bf#1\ }}
\def\JCTA#1 {{\sl J.}\th {\sl Comb.}\g {\sl Theory~A\penalty100\ \bf#1\ }}
\def\JCTB#1 {{\sl J.}\th {\sl Comb.}\g {\sl Theory~B\penalty100\ \bf#1\ }}
\def\JGT#1 {{\sl J.}\th {\sl Graph Theory\penalty100\ \bf#1\ }}
\def\BLMS#1 {{\sl Bull.\g Lond.\g Math.\g Soc.\penalty100\ \bf#1\ }}
\def\JLMS#1 {{\sl J.\g Lond.\g Math.\g Soc.\penalty100\ \bf#1\ }}
\def\PLMS#1 {{\sl Proc.\g Lond.\g Math.\g Soc.\penalty100\ \bf#1\ }}
\def\Order#1 {{\sl Order\ \bf#1\ }}
\def\Random#1 {{\sl Random Struct.\g Alg.\penalty100\ \bf#1\ }}
\def\MA#1 {{\sl Math.}\g {\sl Ann.\penalty100\ \bf#1\ }}
\def\MN#1 {{\sl Math.}\g {\sl Nachr.\penalty100\ \bf#1\ }}
\def\MPCPS#1 {{\sl Math.\g Proc.\g Camb.\g Phil.\g Soc.\penalty100\ \bf#1\ }}
\def\MS#1 {{\sl Math.}\g {\sl Scand.\penalty100\ \bf#1\ }}
\def\MZ#1 {{\sl Math.}\g {\sl Zeit.\penalty100\ \bf#1\ }}
\def\BAMS#1 {{\sl Bull.\th Amer.\g Math.\g Soc.\penalty100\ \bf#1\ }}
\def\JAMS#1 {{\sl J.\th Amer.\g Math.\g Soc.\penalty100\ \bf#1\ }}
\def\MAMS#1 {{\sl Mem.\g Amer.\g Math.\g Soc.\penalty100\ \bf#1\ }}
\def\PAMS#1 {{\sl Proc.\g Amer.\g Math.\g Soc.\penalty100\ \bf#1\ }}
\def\SIAM#1 {{\sl SIAM J.\g Discrete Math.\penalty100\ \bf#1\ }}
\def\SLNM#1 {{\sl Springer Lecture Notes in Mathematics\penalty100\ \bf#1\ }}
\def\TAMS#1 {{\sl Trans.\g Amer.\g Math.\g Soc.\penalty100\ \bf#1\ }}
\def\TCSA#1 {{\sl Theor.\g Comput.\g Sci.~A\penalty100\ \bf#1\ }}
%
%
%
%
%
%
%
%
%
%
%
%
%
%
%
\bigskipamount=1\baselineskip plus.3\baselineskip minus.3\baselineskip
\medskipamount=\bigskipamount\divide\medskipamount by 2
\smallskipamount=\medskipamount\divide\smallskipamount by 2 
\medmuskip = 3mu plus 2mu minus 1mu
\thickmuskip = 6mu plus 4mu minus 2mu 
\def\smallbreak{\par \ifdim\lastskip<\smallskipamount
   \removelastskip \penalty-100 \smallskip \fi}
\def\medbreak{\par \ifdim\lastskip<\medskipamount
   \removelastskip \penalty-250 \medskip \fi}
\def\bigbreak{\par \ifdim\lastskip<\bigskipamount
   \removelastskip \penalty-500 \bigskip \fi}
\catcode`@=11        
  \def\raggedbottom{\topskip10pt plus20pt \r@ggedbottomtrue} 
\catcode`@=12        
\def\ge{\geqslant}
\def\le{\leqslant}
\let\elt=\in
\def\in{\mathrel{\mathchoice
   {\raise .7pt \hbox{$\scriptstyle\elt$}}
   {\raise .7pt \hbox{$\scriptstyle\elt$}}
   {\raise .5pt \hbox{$\hskip .5pt\scriptscriptstyle\elt\hskip .5pt$}}
   {\raise.35pt \hbox{$\scriptscriptstyle\elt$}} }}
\let\hasaselt=\owns
\def\owns{\mathrel{\mathchoice
   {\raise .7pt \hbox{$\scriptstyle\hasaselt$}}
   {\raise .7pt \hbox{$\scriptstyle\hasaselt$}}
   {\raise .5pt \hbox{$\hskip .5pt\scriptscriptstyle\hasaselt\hskip .5pt$}}
   {\raise.35pt \hbox{$\scriptscriptstyle\hasaselt$}} }}
\let\exis=\exists
   \def\exists{\exis\>}
\let\nexis=\nexists
   \def\nexists{\nexis\>}
\let\foral=\forall
   \def\forall{\foral\>}
\let\Rightarro=\Rightarrow
   \def\Rightarrow{\>\Rightarro\>}
\let\mi=\min
   \def\min{\mi\>}
\let\ma=\max
   \def\max{\ma\>}
\let\su=\sup
   \def\sup{\su\>}
\let\inff=\inf
   \def\inf{\inff\>}
\mathchardef\to="2221   
\def\proclaim #1.#2 #3\par{\bigbreak
   \noindent{\bf#1.}#2\enspace{\sl#3}\par\bigbreak}
%
\newskip\sectionheadlineskipamount
\sectionheadlineskipamount=8pt plus 2pt minus 1pt
\def\beginsection #1\par{\Bigbreak\centerline{\bold #1}
        \penalty10000\vskip\sectionheadlineskipamount\noindent}
\let\ffootnote=\footnote
\def\footnote#1#2{\ffootnote{#1}{\eightpoint#2\vskip-12pt}}
%
%
\def\item#1#2\par{\parindent=10mm\smallbreak\hang\indent
                  \llap{{\rm #1}\enspace}\ignorespaces#2\par\smallbreak
                  \parindent=7mm}
\def\itemitem#1#2\par{\parindent=10mm\smallbreak
                  \indent\hangindent2\parindent\indent
                  \llap{{\rm #1}\enspace}\ignorespaces#2\par\smallbreak
                  \parindent=7mm}
%
%
%
\pretolerance=0 
\tolerance=2000
\baselineskip=13pt                 
\vsize=200mm                   
\hsize=120mm                   
\hoffset=9mm                   
\parindent=7mm
\relpenalty=2000 \binoppenalty=5000  
\hyphenpenalty=100
\abovedisplayskip=12pt plus3pt minus4pt
\belowdisplayskip=12pt plus3pt minus4pt    
%
%
\belowdisplayshortskip=9pt plus3pt minus3pt
%
%
%
%
%
 \hyphenation{Baum-ord-nung Baum-ord-nun-gen End-ecke End-ecken kur-zen
Kur-zen Graphen-ei-gen-schaft Graphen-ei-gen-schaften he-raus he-raus-ar-bei-ten
he-raus-zu-ar-bei-ten Schnitt-raum}%
 \hyphenation{ac-cess-ible ana-log-ous ana-log-ous-ly ana-lyze ana-lyse
ana-ly-sis answer answers aver-age bundle bundles Buch-ge-sell-schaft col-our
col-ours col-oured col-our-ing col-our-ings con-struct-ible con-struct-ive
con-struct-ive-ly co-rol-lary Co-rol-lary des-cend des-cend-ing Deut-sche
end-li-cher de-fi-ni-tion de-fi-ni-tions De-fi-ni-tion equi-val-ent
equi-val-ence Euler-ian exist-ence every Gra-phen Hamil-ton-ian homeo-mor-phic
homeo-mor-phism homeo-mor-phisms hy-po-thesis hy-po-theses in-ac-cess-ible
ir-regu-lar ir-regu-lar-ity method methods modi-fi-ca-tion mono-chro-matic par-ticu-lar
pro-po-si-tion pro-po-si-tions Pro-po-si-tion regu-lar regu-lar-ity regu-lar-ly
sig-ni-fi-cant sig-ni-fi-cant-ly sig-ni-fi-cance to-po-lo-gical to-po-lo-gical-ly
un-at-tached un-end-li-cher using Using Wis-sen-schaft-li-che}
%

\userefnumbers\let\newcounter=\newcount\refno =0
\global \advance \refno by 1 \newcounter \refAharoniCtbleEM  \refAharoniCtbleEM =\the \refno 
\global \advance \refno by 1 \newcounter \refAharoniBergerEM  \refAharoniBergerEM =\the \refno 
\global \advance \refno by 1 \newcounter \refMFMC  \refMFMC =\the \refno 
\global \advance \refno by 1 \newcounter \refAharoniCT  \refAharoniCT =\the \refno 
\global \advance \refno by 1 \newcounter \refABLLeftRightTours  \refABLLeftRightTours =\the \refno 
\global \advance \refno by 1 \newcounter \refAsratianKhachatrianLocalization  \refAsratianKhachatrianLocalization =\the \refno 
\global \advance \refno by 1 \newcounter \refBeanFinitary  \refBeanFinitary =\the \refno 
\global \advance \refno by 1 \newcounter \refBenjaminiSchramm  \refBenjaminiSchramm =\the \refno 
\global \advance \refno by 1 \newcounter \refBergerBruhnEndDegrees  \refBergerBruhnEndDegrees =\the \refno 
\global \advance \refno by 1 \newcounter \refBiggsPotential  \refBiggsPotential =\the \refno 
\global \advance \refno by 1 \newcounter \refEGT  \refEGT =\the \refno 
\global \advance \refno by 1 \newcounter \refBruhnPeripheral  \refBruhnPeripheral =\the \refno 
\global \advance \refno by 1 \newcounter \refBruhnPersonal  \refBruhnPersonal =\the \refno 
\global \advance \refno by 1 \newcounter \refDuality  \refDuality =\the \refno 
\global \advance \refno by 1 \newcounter \refMacLaneArbitrarySurfaces  \refMacLaneArbitrarySurfaces =\the \refno 
\global \advance \refno by 1 \newcounter \refInfiniteMatroidAxioms  \refInfiniteMatroidAxioms =\the \refno 
\global \advance \refno by 1 \newcounter \refTreeEnds  \refTreeEnds =\the \refno 
\global \advance \refno by 1 \newcounter \refCyCoCy  \refCyCoCy =\the \refno 
\global \advance \refno by 1 \newcounter \refBruhnSteinDiestelEM  \refBruhnSteinDiestelEM =\the \refno 
\global \advance \refno by 1 \newcounter \refAgelosHenningLA  \refAgelosHenningLA =\the \refno 
\global \advance \refno by 1 \newcounter \refBruhnBicycles  \refBruhnBicycles =\the \refno 
\global \advance \refno by 1 \newcounter \refBruhnSteinMacLane  \refBruhnSteinMacLane =\the \refno 
\global \advance \refno by 1 \newcounter \refBruhnSteinEndDeg  \refBruhnSteinEndDeg =\the \refno 
\global \advance \refno by 1 \newcounter \refBruhnSteinEndDuality  \refBruhnSteinEndDuality =\the \refno 
\global \advance \refno by 1 \newcounter \refBruhnYuHamilton  \refBruhnYuHamilton =\the \refno 
\global \advance \refno by 1 \newcounter \refCasteelsRichterBicycles  \refCasteelsRichterBicycles =\the \refno 
\global \advance \refno by 1 \newcounter \refCatlinLineGraphs  \refCatlinLineGraphs =\the \refno 
\global \advance \refno by 1 \newcounter \refCoornaertDelzantPapadopoulos  \refCoornaertDelzantPapadopoulos =\the \refno 
\global \advance \refno by 1 \newcounter \refYuHamilton  \refYuHamilton =\the \refno 
\global \advance \refno by 1 \newcounter \refCtbleEM  \refCtbleEM =\the \refno 
\global \advance \refno by 1 \newcounter \refCyclesExpository  \refCyclesExpository =\the \refno 
\global \advance \refno by 1 \newcounter \refSpanningTrees  \refSpanningTrees =\the \refno 
\global \advance \refno by 1 \newcounter \refTopSurveyI  \refTopSurveyI =\the \refno 
\global \advance \refno by 1 \newcounter \refTopSurveyII  \refTopSurveyII =\the \refno 
\global \advance \refno by 1 \newcounter \refBook  \refBook =\the \refno 
\global \advance \refno by 1 \newcounter \refMinorUniversal  \refMinorUniversal =\the \refno 
\global \advance \refno by 1 \newcounter \refCyclesOne  \refCyclesOne =\the \refno 
\global \advance \refno by 1 \newcounter \refCyclesTwo  \refCyclesTwo =\the \refno 
\global \advance \refno by 1 \newcounter \refTST  \refTST =\the \refno 
\global \advance \refno by 1 \newcounter \refNST  \refNST =\the \refno 
\global \advance \refno by 1 \newcounter \refHomologyGraphs  \refHomologyGraphs =\the \refno 
\global \advance \refno by 1 \newcounter \refHomotopy  \refHomotopy =\the \refno 
\global \advance \refno by 1 \newcounter \refHomologySpaces  \refHomologySpaces =\the \refno 
\global \advance \refno by 1 \newcounter \refHomSurvey  \refHomSurvey =\the \refno 
\global \advance \refno by 1 \newcounter \refEilenbergSteenrod  \refEilenbergSteenrod =\the \refno 
\global \advance \refno by 1 \newcounter \refAgelosPathConnected  \refAgelosPathConnected =\the \refno 
\global \advance \refno by 1 \newcounter \refAgelosOWReportFleischner  \refAgelosOWReportFleischner =\the \refno 
\global \advance \refno by 1 \newcounter \refAgelosHotchpotch  \refAgelosHotchpotch =\the \refno 
\global \advance \refno by 1 \newcounter \refAgelosInfiniteFleischner  \refAgelosInfiniteFleischner =\the \refno 
\global \advance \refno by 1 \newcounter \refAgelosLTop  \refAgelosLTop =\the \refno 
\global \advance \refno by 1 \newcounter \refAgelosUniqueFlows  \refAgelosUniqueFlows =\the \refno 
\global \advance \refno by 1 \newcounter \refAgelosPersonal  \refAgelosPersonal =\the \refno 
\global \advance \refno by 1 \newcounter \refAgelosPhilippGeodetic  \refAgelosPhilippGeodetic =\the \refno 
\global \advance \refno by 1 \newcounter \refGromovHyperbolicGroups  \refGromovHyperbolicGroups =\the \refno 
\global \advance \refno by 1 \newcounter \refHahnEdgeEnds  \refHahnEdgeEnds =\the \refno 
\global \advance \refno by 1 \newcounter \refHalinInfGrid  \refHalinInfGrid =\the \refno 
\global \advance \refno by 1 \newcounter \refHalinMinVertex  \refHalinMinVertex =\the \refno 
\global \advance \refno by 1 \newcounter \refHalinMinimization  \refHalinMinimization =\the \refno 
\global \advance \refno by 1 \newcounter \refHalinInfMinimization  \refHalinInfMinimization =\the \refno 
\global \advance \refno by 1 \newcounter \refHalinProblems  \refHalinProblems =\the \refno 
\global \advance \refno by 1 \newcounter \refHatcher  \refHatcher =\the \refno 
\global \advance \refno by 1 \newcounter \refHiggsMatroidsDuality  \refHiggsMatroidsDuality =\the \refno 
\global \advance \refno by 1 \newcounter \refHiggsBMatroids  \refHiggsBMatroids =\the \refno 
\global \advance \refno by 1 \newcounter \refHiggsMatroids  \refHiggsMatroids =\the \refno 
\global \advance \refno by 1 \newcounter \refLick  \refLick =\the \refno 
\global \advance \refno by 1 \newcounter \refMaderHomEigenschaften  \refMaderHomEigenschaften =\the \refno 
\global \advance \refno by 1 \newcounter \refMaderMinVertex  \refMaderMinVertex =\the \refno 
\global \advance \refno by 1 \newcounter \refMaderMinVertices  \refMaderMinVertices =\the \refno 
\global \advance \refno by 1 \newcounter \refMaderMinVerticesInf  \refMaderMinVerticesInf =\the \refno 
\global \advance \refno by 1 \newcounter \refMaderReduktion  \refMaderReduktion =\the \refno 
\global \advance \refno by 1 \newcounter \refMaderEdgeConPres  \refMaderEdgeConPres =\the \refno 
\global \advance \refno by 1 \newcounter \refTheoFlows  \refTheoFlows =\the \refno 
\global \advance \refno by 1 \newcounter \refNadlerContinuumTheory  \refNadlerContinuumTheory =\the \refno 
\global \advance \refno by 1 \newcounter \refNWTreePacking  \refNWTreePacking =\the \refno 
\global \advance \refno by 1 \newcounter \refNWArboricity  \refNWArboricity =\the \refno 
\global \advance \refno by 1 \newcounter \refOberlySumner  \refOberlySumner =\the \refno 
\global \advance \refno by 1 \newcounter \refOxleyInfiniteMatroidsPaper  \refOxleyInfiniteMatroidsPaper =\the \refno 
\global \advance \refno by 1 \newcounter \refOxleyInfiniteMatroidsSurvey  \refOxleyInfiniteMatroidsSurvey =\the \refno 
\global \advance \refno by 1 \newcounter \refRichterGraphLikeSpacesSurvey  \refRichterGraphLikeSpacesSurvey =\the \refno 
\global \advance \refno by 1 \newcounter \refRichterVella  \refRichterVella =\the \refno 
\global \advance \refno by 1 \newcounter \refSchulzEdgeEnds  \refSchulzEdgeEnds =\the \refno 
\global \advance \refno by 1 \newcounter \refPhilippNormal  \refPhilippNormal =\the \refno 
\global \advance \refno by 1 \newcounter \refMayaTreePacking  \refMayaTreePacking =\the \refno 
\global \advance \refno by 1 \newcounter \refMayaEndDeg  \refMayaEndDeg =\the \refno 
\global \advance \refno by 1 \newcounter \refMayaBanff  \refMayaBanff =\the \refno 
\global \advance \refno by 1 \newcounter \refCTFleischner  \refCTFleischner =\the \refno 
\global \advance \refno by 1 \newcounter \refCTPlanDualInf  \refCTPlanDualInf =\the \refno 
\global \advance \refno by 1 \newcounter \refCTConPres  \refCTConPres =\the \refno 
\global \advance \refno by 1 \newcounter \refCTLineGraphs  \refCTLineGraphs =\the \refno 
\global \advance \refno by 1 \newcounter \refCTVella  \refCTVella =\the \refno 
\global \advance \refno by 1 \newcounter \refTutteTreePacking  \refTutteTreePacking =\the \refno 
\global \advance \refno by 1 \newcounter \refVellaThesis  \refVellaThesis =\the \refno 
\global \advance \refno by 1 \newcounter \refLasVergnasMatroids  \refLasVergnasMatroids =\the \refno 
\global \advance \refno by 1 \newcounter \refWagnerBook  \refWagnerBook =\the \refno 
\global \advance \refno by 1 \newcounter \refWhyburnMenger  \refWhyburnMenger =\the \refno 
\global \advance \refno by 1 \newcounter \refWoessDirichlet  \refWoessDirichlet =\the \refno 
\global \advance \refno by 1 \newcounter \refWoessBook  \refWoessBook =\the \refno 
\global \advance \refno by 1 \newcounter \refZhanLineGraphs  \refZhanLineGraphs =\the \refno 

\hbox{}\vskip2pt\centerline{\bigbold Locally finite graphs with ends: a topological approach$^*$}
\vskip 4mm
\centerline{Reinhard Diestel}
\def\problem#1{\?{PROBLEM #1}}
   \def\problem#1{}\def\?#1{}
\disablems
\noremarks
\nocomments
\input epsf.tex 
\def\secConcepts{1}%
   \def\xxxCompactMetricG{\secConcepts.1}
   \def\xxxArcs{\secConcepts.2}
   \def\xxxArcComponents{\secConcepts.3}
   \def\xxxJumpingArc{\secConcepts.4}
   \def\xxxConnected{\secConcepts.5}
   \def\xxxTSTexistence{\secConcepts.6}
   \def\xxxTSTeq{\secConcepts.7}
   \def\xxxFundamentalBasic{\secConcepts.8}
   \def\xxxFundamentalDuality{\secConcepts.9}
\def\secTheory{2}%
   \def\xxxFundamental{\secTheory.1}
   \def\xxxNSTgenerators{\secTheory.2}
   \def\xxxGenerating{\secTheory.3}
   \def\xxxGeodetic{\secTheory.4}
   \def\xxxCBasics{\secTheory.5}
   \def\xxxOrthogonality{\secTheory.6}
   \def\xxxIsos{\secTheory.7}
   \def\xxxExtensions{\secTheory.8}
   \def\xxxOrthDec{\secTheory.9}
   \def\xxxInfOrth{\secTheory.10}
   \def\xxxBicycles{\secTheory.11}
   \def\xxxDec{\secTheory.12}
   \def\xxxBoundaryOperators{\secTheory.13}
   \def\xxxDualTrees{\secTheory.14}
   \def\xxxDualTopTrees{\secTheory.15}
   \def\xxxGeomDuals{\secTheory.16}
\def\secTechniques{3}%
   \def\xxxkCon{\secTechniques.1}
   \def\xxxTST{\secTechniques.2}
   \def\xxxDirectArc{\secTechniques.3}
\def\secApplications{4}%
   \def\xxxEndConnStein{\secApplications.1}
   \def\xxxEndDegH{\secApplications.2}
   \def\xxxEndDegSubspaces{\secApplications.3}
   \def\xxxEndDegInf{\secApplications.4}
   \def\xxxEdgeConnPresCircles{\secApplications.5}
   \def\xxxVxConnPresCircles{\secApplications.6}
   \def\xxxMinkCon{\secApplications.7}
   \def\xxxDegk{\secApplications.8}
   \def\xxxCircleDegk{\secApplications.9}
   \def\xxxVxMinimal{\secApplications.10}
   \def\xxxTreePacking{\secApplications.11}
   \def\xxxMayaArboricity{\secApplications.12}
   \def\xxxLocalHamilton{\secApplications.13}
   \def\xxxLocalHamiltonOS{\secApplications.14}
   \def\xxxInfiniteFleischner{\secApplications.15}
   \def\xxxLineGraphs{\secApplications.16}
   \def\xxxInfiniteTutte{\secApplications.17}
   \def\xxxToughness{\secApplications.18}
   \def\xxxMacLane{\secApplications.19}
   \def\xxxKelmans{\secApplications.20}
   \def\xxxWhitney{\secApplications.21}
   \def\xxxFCT{\secApplications.22}
   \def\xxxFlows{\secApplications.23}
   \def\xxxBoundarySourceSink{\secApplications.24}
   \def\xxxDirichlet{\secApplications.25}
   \def\xxxNflows{\secApplications.26}
   \def\xxxHflows{\secApplications.27}
\def\secOutlook{5}
   \def\xxxVTop{\secOutlook.1}
   \def\xxxRichterVellaPeano{\secOutlook.2}
   \def\xxxMetricNST{\secOutlook.3}
   \def\xxxLtopModG{\secOutlook.4}
   \def\xxxFiniteSums{\secOutlook.5}
   \def\xxxMetricHomology{\secOutlook.6}
   \def\xxxCasInverseLimit{\secOutlook.7}
   \def\xxxSingular{\secOutlook.8}
   \def\xxxSingularAlt{\secOutlook.9}
   \def\xxxMatroidWhitney{\secOutlook.10}
   \def\xxxMatroidExamples{\secOutlook.11}
\def\refBookSection{\refBook, Ch.~8.5}
\def\figThreeEnds{0}
\def\figWildCircle{\secConcepts.1}
\def\figTSTs{\secConcepts.2}
\def\figPedestrian{\secTheory.1}
\def\figIToPDuals{\secTheory.2}
\def\figSerpents{\secTheory.3}
\def\figMinK{\secTechniques.1}
\def\figBadLimit{\secTechniques.2}
\def\figTKfour{\secApplications.1}
\def\figMinKalt{\secApplications.2}
\def\figBadEdge{\secApplications.3}
\def\figElusive{\secApplications.4}
\def\figFlow{\secApplications.5}
\def\figFans{\secOutlook.1}
\def\figHyperbolic{\secOutlook.2}
\def\figKringelsimplex{\secOutlook.3}
\def\em{\sl}
\def\NST{normal spanning tree}
\def\TST{topological spanning tree}
\def\ucl(#1){\lfloor #1 \rfloor}
\def\tie{\th }
\def\vxspace{{\cal V}}
\def\edgespace{{\cal E}}
\def\cyclespace{{\cal C}}
\def\cutspace{{{\cal C}^*}{}}
\def\fin{_{\rm fin}}
\def\zh{zu\-sam\-men\-h\"an\-gend}
%
\def\lowfwd #1#2#3{{\setbox0\hbox{$#1$}\setbox1\hbox{$E'\!$}
            \mathchoice
            {{\mathop{\kern0pt #1}\limits^{\kern#2pt\raise.#3ex
     \vbox to 0pt{\hbox{$\scriptscriptstyle\rightarrow$}\vss}}}}%
            {{\mathop{\kern0pt #1}\limits^{\kern#2pt\raise.#3ex
     \vbox to 0pt{\hbox{$\scriptscriptstyle\rightarrow$}\vss}}}}%
            {\ifdim\wd0<\wd1{\,\vec{#1}\,}\else
     {\mathop{\kern0pt #1}\limits^{\kern#2pt\raise.0ex
     \vbox to 0pt{\hbox{$\scriptscriptstyle\rightarrow$}\vss}}}\fi}%
            {{\vec{#1}}}%
            }}
\def\nlowfwd #1#2#3{{\setbox0\hbox{$#1$}\setbox1\hbox{$E'\!$}
            \mathchoice
            {{\mathop{\kern0pt #1}\limits^{\kern#2pt\raise.#3ex
     \vbox to 0pt{\hbox{$\scriptscriptstyle\Rightarrow$}\vss}}}}%
            {{\mathop{\kern0pt #1}\limits^{\kern#2pt\raise.#3ex
     \vbox to 0pt{\hbox{$\scriptscriptstyle\Rightarrow$}\vss}}}}%
            {\ifdim\wd0<\wd1{\,\vec{#1}\,}\else
     {\mathop{\kern0pt #1}\limits^{\kern#2pt\raise.0ex
     \vbox to 0pt{\hbox{$\scriptscriptstyle\Rightarrow$}\vss}}}\fi}%
            {{\vec{#1}}}%
            }}
\def\fwd #1#2{{\lowfwd{#1}{#2}{15}}}
\def\lowbkwd #1#2#3{{\mathop{\kern0pt #1}\limits^{\kern#2pt\raise.#3ex
     \vbox to 0pt{\hbox{$\scriptscriptstyle\leftarrow$}\vss}}}}

\def\vC{\kern-1pt\fwd C3\kern-.5pt}

\def\vd{\kern-1pt\lowfwd d2{10}\kern-1pt}
\def\vD{\kern-.7pt\fwd D3\kern-.5pt}
\def\ve{\kern-1pt\lowfwd e{1.5}1\kern-1pt}
\def\nve{\kern-1pt\nlowfwd e{1.5}2\kern-1pt}
\def\vf{\kern-1pt\lowfwd f{1.5}1\kern-1pt}
\def\ev{\kern-1pt\lowbkwd e{1.5}1\kern-1pt}
\def\veStar{{\mathop{\kern0pt e\lower1.5pt\hbox{${}^*$}}\limits^{\kern0pt
   \raise.02ex\vbox to 0pt{\hbox{$\scriptscriptstyle\rightarrow$}\vss}}}}
\def\eStarv{{\mathop{\kern0pt e\lower1.5pt\hbox{${}^*$}}\limits^{\kern0pt
   \raise.02ex\vbox to 0pt{\hbox{$\scriptscriptstyle\leftarrow$}\vss}}}}
\def\vedash{{\mathop{\kern0pt e\lower.5pt\hbox{${}
     \scriptstyle'$}}\limits^{\kern0pt\raise.02ex
     \vbox to 0pt{\hbox{$\scriptscriptstyle\rightarrow$}\vss}}}}

\def\dm{\mathop{\kern 0pt -}\limits^{\textstyle\vbox to 0pt{\vskip3pt\hbox{.}\vss}}}

\def\vE{\kern-.7pt\fwd E3\kern-.7pt}

\def\vF{\kern-.7pt\fwd F3\kern-.7pt}

\def\vG{\kern-.7pt\fwd G3\kern-.7pt}
\def\vH{\kern-.5pt\fwd H3\kern-.5pt}
\def\vP{\kern-.7pt\fwd P3\kern-.6pt}


\def\vQ{\kern-.7pt\fwd Q3\kern-.6pt}

\def\B{{\cal B}}
\def\H{{\cal H}}

\def\C{{\cal C}}
\def\D{{\cal D}}
\def\FF{{\cal F}}
\def\F{{\openface F}}
\lineskiplimit=-6pt
\newcount\footnoteno \footnoteno=1
\def\Footnote#1{{\footnote{${}^{\the\footnoteno}$}{#1}%
   \global\advance\footnoteno by 1}}
\def\Top{{\tenrm T\eightrm OP}}
\def\MTop{{\tenrm MT\eightrm OP}}
\def\VTop{{\tenrm VT\eightrm OP}}
\def\ITop{{\tenrm IT\eightrm OP}}
\def\ETop{{\tenrm ET\eightrm OP}}
\def\fnTop{{\eightrm T\sevenrm OP}}
\def\fnMTop{{\eightrm MT\sevenrm OP}}
\def\fnVTop{{\eightrm VT\sevenrm OP}}
\def\fwd #1#2{{\lowfwd{#1}{#2}{15}}}
\def\vCC{\kern-.7pt\fwd{\C}3\kern-.7pt}
\def\vBB{\kern-.7pt\fwd{\B}3\kern-.7pt}
\def\vEE{\kern-.7pt\fwd{\edgespace}3\kern-.7pt}
\def\cutspace{\B}
\bigskip\medskip
{\narrower\narrower\ninepoint\noindent
   This paper is intended as an introductory survey of a newly emerging field: a topological approach to the study of locally finite graphs that crucially incorporates their ends. Topological arcs and circles, which may pass through ends, assume the role played in finite graphs by paths and cycles.\nl\indent
   This approach has made it possible to extend to locally finite graphs many classical theorems of finite graph theory that do not extend verbatim. The shift of paradigm it proposes is thus as much an answer to old questions as a source of new ones; many concrete problems of both types are suggested in the paper. \nl\indent
   This paper attempts to provide an entry point to this field for readers that have not followed the literature that has emerged in the last 10 years or so. It takes them on a quick route through what appear to be the most important lasting results, introduces them to key proof techniques, identifies the most promising open problems, and offers pointers to the literature for more detail.\par}

\vskip-\medskipamount\vskip0pt

  \vfootnote*{\eightpoint This paper has appeared in two parts~[\the\refTopSurveyI,\th \the\refTopSurveyII]. It is complemented by~[\the\refHomSurvey], the third part, which studies the algebraic-topological aspects of the theory.\vskip-12pt}

\begingroup\ninepoint
\beginsection Contents

\vskip-\smallskipamount

\def\n#1\par {\noindent #1\smallskip}
\def\nn#1\par {#1\par}

\n Introduction

\n 1.\ Concepts and Basic Theory

\n 2.\ The Topological Cycle Space

\nn 2.1 Generating sets

\nn 2.2 Characterizations of algebraic cycles

\nn 2.3 Cycle-cut orthogonality

\nn 2.4 Orthogonal decomposition

\nn 2.5 Duality\medbreak

\n 3.\ Proof Techniques

\nn 3.1 The direct approach

\nn 3.2 The use of compactness

\nn 3.3 Constructing arcs directly\medbreak

\n 4.\ Applications

\nn 4.1 Extremal infinite graph theory

\nn 4.2 Cycle space applications

\nn 4.3 Flows in infinite graphs and networks\medbreak

\n 5.\ Outlook

\nn 5.1 Graphs with infinite degrees

\nn 5.2 The identification topology

\nn 5.3 Compactification versus metric completion

\nn 5.4 Homology of locally compact spaces with ends

\nn 5.5 Infinite matroids

\endgroup

\beginsection Introduction

This paper describes a topological framework in which many well-known theorems about finite graphs that appear to fail for infinite graphs have a natural infinite analogue. It has been realised in recent years that many such theorems, especially about paths and cycles, do work in a slightly richer setting: not in the (locally finite) graph itself, but in its compactification obtained by adding its {\it ends\/}. For example, the plane graph $G$ in Figure~\figThreeEnds\ has three ends. When we add these, we obtain a compact space~$|G|$ in which the fat edges form a \hbox{\it circle\/}---a subspace homeomorphic to the standard topological circle~$S^1$. Allowing such circles as `infinite cycles', and allowing topological arcs through ends as `infinite paths', we can restore the truth of many well-known theorems about finite graphs whose infinite analogues would fail if we allowed only the usual finite paths and cycles familiar from finite graphs.

 \figure \figThreeEnds.
 A circle through three ends
 (ThreeCycle; 1000)

The aim of this paper is to provide a reasonably complete but readable introduction to this new approach, offering a fast track to its current state of the art. It describes all the fundamental concepts, all the main results, and many open problems. Given this aim, there will not be the space for many proofs. However, there is a complementary source~[\the\refBookSection] that offers detailed proofs of all the most basic facts (which are here included in the narrative but not all proved), and which the interested reader may wish to consult early for some more `feel' for the subject. For readers already familiar with those basic techniques, the current paper also describes some more advanced but fundamental proof techniques that cannot be found explicitly in~[\the\refBook].

We begin in Section~\secConcepts\ with the definition of the space $|G|$ and an overview of its {\it basic properties\/}. These properties mostly concern the topological analogues of familiar finite concepts involving paths and cycles, such as spanning trees, degrees (of ends), connectivity and so on. It is these analogues that we shall need in the place of their finite counterparts when we wish to extend theorems from finite graph theory to infinite graphs whose naive extension fails.

We continue in Section~\secTheory\ with some {\it basic theory\/}, a body of results that are not yet applications of the new topological concepts but relate them to each other, much in the way their finite counterparts are related. Most results in this section concern the homology of a graph, i.e., its {\it cycle\/} and {\it cut space\/} and the way they interact. These homology aspects have, so far, been the prime field of application for those topological notions involving ends. But there are other applications too, and reading Section~\secTheory\ will not be a prerequisite for reading the rest of the paper.

Section~\secTechniques\ explains some {\it proof techniques\/}, e.g.\ for the construction of topological arcs and circles as limits of finite paths and cycles, that have evolved over the past 10 years. These years have seen some considerable simplification of the techniques required to deal with the problems one usually encounters in this area. The aim of Section~\secTechniques\ is to describe the state of the art here, so as to equip those new to the field quickly with the main techniques now available.

Section~\secApplications\ is devoted to {\it applications\/}: theorems that involve our new topological concepts but answer questions that could be asked without them. Often, these questions simply ask for an extension to infinite graphs of well known finite theorems whose naive infinite extensions are either false or trivially true. Examples include all the standard planarity criteria in terms of the cycle space, some classical theorems about Euler tours and Hamilton cycles, connectivity results such as the tree-packing theorem of Nash-Williams and Tutte, and infinite electrical networks. The consideration of arcs and circles instead of paths and cycles has made it possible to ask extremal-type questions about the interaction of graph invariants that were so far meaningful only for finite graphs. Some of these are indicated in Section~\secApplications; more are given in Stein~[\the\refMayaBanff].

Finally, there is an {\it outlook\/} in Section~\secOutlook\ to new horizons: extensions to graphs that are not locally finite, and implications of our findings in related fields, such as geometric group theory and infinite matroids. Homological aspects from an algebraic point of view are also indicated; these will be explored more fully in~[\the\refHomSurvey].

{\it Open problems\/} are not collected at the end but interspersed within the text. There are plenty of these, as well as the overall quest to push the general approach further: to identify more theorems about paths and cycles in finite graphs that do not extend naively, and to find the correct topological analogue that does extend.

\section\secConcepts. Concepts and Basic Facts

Throughout this section, let $G$ be a fixed infinite, locally finite, connected graph. This section serves to introduce the concepts on which our topological approach to the study of such graphs is based: the space $|G|$ formed by $G$ and its ends; topological paths, circles and spanning trees in this space; notions of connectivity in~$|G|$. The style will be descriptive and informal, aiming for overall readability; should any technical points remain unclear, the reader is referred to~[\the\refBookSection] for more formal definitions of the concepts introduced here, and to [\the\refBook] in general for graph-theoretic terms and notation.

Terms such as `path' or `connected', which formally have different meanings in topology and in graph theory, will be used according to context: in the graph-theoretical sense for graphs, and in the usual topological sense for topological spaces. If the context is ambiguous, the two meanings will probably coincide, making a formal distinction unnecessary.

We call 1-way infinite paths {\it rays\/}, and 2-way infinite paths {\it double rays\/}. An {\it end\/} of $G$ is an equivalence class of rays in~$G$, where two rays are considered equivalent if no finite set of vertices separates them in~$G$. The graph shown in Figure~\figThreeEnds\ has three ends; the $\Z\times\Z$ grid has only one, the infinite binary tree has continuum many. We write $\Omega(G)$ for the set of ends of~$G$.

Topologically, we view $G$ as a cell complex with the usual topology. Adding its ends compactifies it, with the topology generated by the open sets of $G$ and the following additional basic open sets. For every finite set $S$ of vertices and every end~$\omega$, there is a unique component $C$ of $G-S$ in which every ray of $\omega$ has a tail. We say that $\omega$ {\it lives in\/}~$C$ and write $C =: C(S,\omega)$. Now for every such~$S$ and every component $C$ of~$G-S$, we declare as open the union $\hat C$ of $C$ with the set of ends living in $C$ and with all the `open' $S$--$C$ edges of~$G$ (i.e., without their endpoints in~$S$). We denote the space just obtained by~$|G|$.%
   \Footnote{The open neighbourhoods of ends are defined slightly more generally in~[\the\refBook], but for locally finite graphs the two definitions are equivalent. Topologies for graphs with infinite degrees are discussed in Section~\secOutlook.}

\proclaim Theorem \xxxCompactMetricG.~[\the\refSpanningTrees,\th\the\refBook]
The space $|G|$ is compact and metrizable.%
   \COMMENT{}

The space $|G|$ is known as the {\it Freudenthal compactification\/} of~$G$. The main feature of its topology is that rays in $G$ converge as they should: to the end of which they are an element. One can show~[\the\refCyclesOne,~Prop.\th 4.5] that every Hausdorff topology on $G\cup\Omega(G)$ with (essentially)%
   \COMMENT{}
   this feature (and which induces the 1-complex topology on~$G$) refines the topology of~$|G|$. This identifies $|G|$ as the unique most powerful Hausdorff topology on $G\cup\Omega(G)$, in a sense that can be made quite precise~[\the\refCyclesOne].%
   \Footnote{The identity to~$|G|$ from $G\cup\Omega(G)$ with any finer topology is continuous, so in $|G|$ there will be at least as many arcs and circles (and possibly more). Arcs and circles in $|G|$ will be our main tools.}

Of the many natural aspects of this topology let us mention just two more, which relate it to better-known objects. Consider the binary tree~$T_2$, and think of its rays from the root as 0--1 sequences. The resulting bijection between the ends of $T_2$ and these sequences is a homeomorphism between $\Omega(T_2)$, as a subspace of~$|T_2|$, and $\{0,1\}^\N$ with the product topology. Identifying pairs of ends whose sequences specify the same rational (one sequence ending on zeros, the other on~1s) turns this bijection into a homeomorphism from the resulting identification space of $\Omega(T_2)$ to~$[0,1]$. Without such identification, on the other hand, $\Omega(G)$~is always a subset of a Cantor set.

\medbreak

Instead of paths and cycles in $G$ we can now consider {\it arcs\/} and {\it circles\/} in~$|G|$: homeomorphic images of the real interval $[0,1]$ and of the complex unit circle~$S^1$. While paths and cycles are examples of arcs and circles, Figure~\figThreeEnds\ shows a circle that is not a cycle. Arcs and circles that are not paths or cycles must contain ends. An arc containing uncountably many ends always induces the ordering of the rationals on a subset of its vertices~[\the\refTreeEnds]. Such arcs, and circles containing them, are called {\it wild\/} (Figure~\figWildCircle), but they are quite common.

 \figure \figWildCircle. 
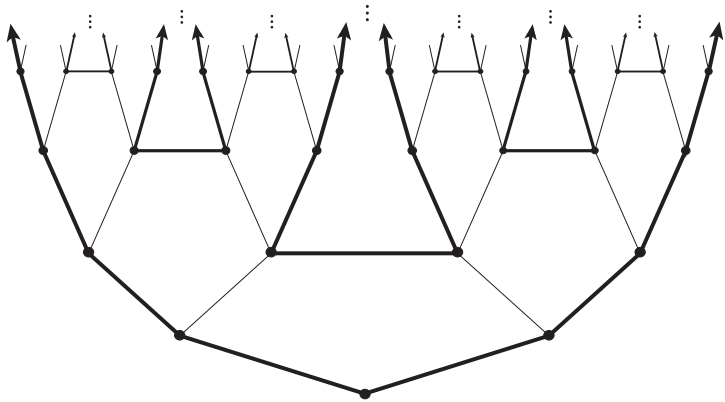
 The heavy edges form a {\it wild circle\/}.
 (WildCircle; 1000)

Arcs and circles are examples of a natural type of subspace of~$|G|$: subspaces that are the closure in $|G|$ of some subgraph of~$G$.%
  \Footnote{It takes a line of proof that arcs and circles do indeed have this property, i.e., that the union of their edges is dense in them. This is because any arc between distinct ends must contain an edge incident with finite set $S$ that separates them; this will follow from Lemma~\xxxJumpingArc\ below.}
  We call such a subspace $X$ of $|G|$ a {\it standard subspace\/}, and write $V(X)$ and $E(X)$ for the set of vertices or edges it contains. Note that the ends in $X$ are ends of~$G$, not of the subgraph that gave rise to~$X$; in particular, ends in $X$ need not have a ray in that subgraph. Which ends of $G$ are in $X$ is determined just by $V(X)$: they are precisely the ends that are limits of vertices in~$X$.

Given a standard subspace~$X$ and an end $\omega\in X$, the maximum number of arcs in $X$ that end in $\omega$ but are otherwise disjoint is the {\it (vertex-) degree\/} of $\omega$ in~$X$; the maximum number of edge-disjoint arcs in $X$ ending in~$\omega$ is its {\it edge-degree\/} in~$X$. Both maxima are indeed attained, but it is non-trivial to prove this~[\the\refBruhnSteinEndDeg].%
   \COMMENT{}
   End degrees behave largely as expected; for example, the connected standard subspaces in which every vertex and every end has (vertex-) degree~2 are precisely the circles. (Use Lemma~\xxxArcs\ to prove this.) In Section~\secApplications.1 we shall define a third type of end degrees, their {\it relative degree\/}, which is useful for the application of end degrees to extremal-type problems about infinite graphs.

Standard subspaces have the important property that connectedness and arc-connectedness are equivalent for them. This will often be convenient: while connectedness is much easier to prove (see Lemma~\xxxConnected), it is usually arc-connectedness that we need.

\proclaim Lemma \xxxArcs.~[\the\refBook,\th\the\refTST,\th\the\refCTVella] Connected standard subspaces of $|G|$ are locally connected and arc-connected.

\noindent
   The proof that a connected standard subspace $X$ is locally connected is not hard: an open neighbourhood $\hat C\cap X$ of an end $\omega$ will be connected if we choose the set $S$ in its definition so as to minimize the number of $C$--$S$ edges in~$X$.
 But local connectedness is not a property we shall often use directly. Its role here is that it offers a convenient stepping stone towards the proof of arc-connectedness.%
   \Footnote{By general `continuum theory'~[\the\refNadlerContinuumTheory], compact, connected and locally connected metric spaces are arc-connected.}
   Direct proofs that $X$ is arc-connected can be found in~[\the\refTST,\th\the\refAgelosPathConnected], and we shall indicate one in Section~\secTechniques. Connected subspaces of $|G|$ that are neither open%
   \COMMENT{}
   nor closed need not be arc-connected~[\the\refAgelosPathConnected].

\proclaim Corollary \xxxArcComponents. The arc-components of a standard subspace are closed.%
   \COMMENT{}

\proof The closure in $X$ of an arc-component of a standard subspace $X$ is connected and itself standard, and hence arc-connected by Lemma~\xxxArcs.
\endproof

The edge set $E(C)$ of any circle $C$ will be called a {\it circuit\/}. Given any set $F$ of edges in~$G$, we write $\overline F$ for the closure of $\bigcup F$ in~$|G|$, and call $\overline F$ the standard subspace {\it spanned by~}$F$. (This is with slight abuse of our usual notation, in which we write $\overline X$ for the closure in $|G|$ of a subset $X\sub |G|$.) Similarly, we write $\interior F$ for the set of all inner points of edges in~$F$ (while usually we write $\interior X$ for the interior of a subset $X\sub |G|$).

The set of edges of $G$ across a partition $\{V_1,V_2\}$ of $V(G)$ is a {\it cut\/} of~$G$; the sets $V_1,V_2$ are the {\it sides\/} of this cut. A minimal non-empty cut is a {\it bond\/}.

The following lemma is one of our basic tools for handling arcs. It says that an arc cannot `jump across' a {\em finite\/} cut without containing an edge from it:\looseness=-1

\proclaimwithname Lemma \xxxJumpingArc. (Jumping Arc Lemma {\rm [\the\refBook]})
   Let $F$ be a cut of $G$ with sides $V_1,V_2$. Let $X$ be a standard subspace of~$|G|$, and put $X_i := X\cap V_i$ $(i=1,2)$.
   \pitem{i} If $F$ is finite, then $\overline{V_1}\cap\overline{V_2}=\es$, and there is no arc in $|G|\sm\interior F$ with one endpoint in ${V_1}$ and the other in~${V_2}$.%
   \COMMENT{}
   \pitem{ii} If $F\cap E(X)$ is infinite, then $\overline{X_1}\cap\overline{X_2}\ne\es$, and there may be such an arc in~$X$.
   \enditem

\noindent
   The proof of Lemma \xxxJumpingArc\th (i) is straightforward from the definition of the topology of~$|G|$: deleting the edges of a finite cut splits $|G|$ into two disjoint open sets. When $F\cap E(X)$ is infinite, an intersection point $\omega\in \overline{X_1}\cap\overline{X_2}$ can be obtained as the limit of two vertex sequences, one in $X_1$ and the other in~$X_2$, that are joined by infinitely many cut edges of~$X$.%
   \COMMENT{}

\medbreak

Although the `jumping arc' is a nice way to memorize Lemma~\xxxJumpingArc, its main assertion is not about arcs but about connectedness. It implies that connectedness for standard subspaces can be characterized in graph-theoretical terms alone, without any explicit mention of ends or the topology of~$|G|$:

\proclaim Lemma \xxxConnected. {\rm [\the\refBook]}
A standard subspace of $|G|$ is connected if and only if it contains an edge from every finite cut of $G$ of which it meets both sides.

We shall say that a standard subspace $X$ of $|G|$ is {\it $k$-edge-connected\/} if the deletion of fewer than $k$ edges will not make it disconnected. Similarly, $X$~is {\it $k$-vertex-connected\/} if $V(X) > k$ and the deletion of fewer than $k$ vertices and their incident edges does not leave a disconnected space. Note that for $k=1$ both notions coincide with ordinary topological connectedness, and that for $X=|G|$ the space $X$ is $k$-edge-connected or $k$-vertex-connected if and only if the graph $G$ is $k$-edge-connected or $k$-connected (by Lemma~\xxxConnected).

How about deleting ends as well as vertices and/or edges? For $X=|G|$, this will never help to disconnect~$X$: if deleting a finite set $U$ of vertices and any set of ends disconnects~$|G|$, then so does the deletion of $U$ alone, and similarly for edges.%
   \COMMENT{}
   For arbitrary standard subspaces~$X$, however, deleting ends can make sense.%
   \COMMENT{}
   It will normally result in a subspace that is no longer standard, but the main reason for primarily considering standard subspaces, that connectedness in them is equivalent to arc-connectedness (Lemma~\xxxArcs), is preserved.\looseness=-1
   \Footnote{Georgakopoulos's characterization of the subspaces that are connected but not arc-connected~[\the\refAgelosPathConnected] implies that any such space has uncountably many arc-components consisting of one end only. One clearly cannot obtain such a space by deleting finitely many ends from a connected standard subspace.}%
   \COMMENT{}

So let us call a subspace $X\sub |G|$ {\it substandard\/} if its closure in $|G|$ is standard (so that $X$ contains no partial edges), and {\it $k$-connected\/} if after the deletion of fewer than $k$ vertices, edges or ends it will still be arc-connected. This makes sense in the context of Menger's theorem, which Thomassen and Vella~[\the\refCTVella] proved for topological spaces $X$ that%
   \vadjust{\penalty-200}
    include all subspaces of~$|G|$: given any two points $a,b\in X$ and $k\in\N$, if for every set $S\sub X$ of fewer than $k$ points there is an $a$--$b$ arc in $X\sm S$, then $X$ contains $k$ arcs from $a$ to~$b$ that pairwise meet only in $a$ and~$b$.%
   \Footnote{There is also version of Menger's theorem for disjoint arcs between sets~$A,B\sub X$ in~[\the\refCTVella], which is easier to prove.}%
   \COMMENT{}
   Hence in a $k$-connected standard or substandard subspace any two vertices or ends can be linked by $k$ independent arcs:%
   \COMMENT{}
   a useful property that can fail in standard subspaces that are merely $k$-vertex-connected.%
   \COMMENT{}

\medbreak

To work with arcs and circles in a way that resembles finite graph theory, we need one more addition to our topological toolkit: the notion of a \TST. A~{\it topological spanning tree\/} of a connected standard subspace $X$ of $|G|$ is an (arc-) connected standard subspace $T\sub X$ of $|G|$ that contains every vertex (and hence every end) of $X$ but contains no circle. A~\TST\ of $|G|$ will also be called a {\it \TST\ of~}$G$.

The closure $\overline T$ of an ordinary spanning tree $T$ of $G$ is not normally a \TST: as soon as $T$ contains disjoint rays from the same end, $\overline T$~will contain a circle. Conversely, the subgraph of $G$ underlying a \TST\ need not be a graph-theoretical tree: it will be acyclic, of course, but it need not be connected. Figure~\figTSTs\ shows examples of both these phenomena in the double ladder.

 \figure \figTSTs.
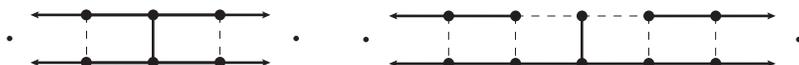
 An ordinary spanning tree of $G$ {(left)}, and\nl a \TST\ of $G$ {(right)}
 (TSTs; 1000)
 
Ordinary spanning trees whose closures are \TST s do always exist, however: all {\em normal\/} spanning trees have this property, and all countable connected graphs have \NST s. (A~spanning tree $T$ of $G$ is {\it normal\/} if, for a suitable choice of a root, the endvertices of every edge of $G$ are comparable in the tree-order of~$T$. See~[\the\refBook].) Often, therefore, \NST s are the best choice of a spanning tree for our purposes.

More generally, we have the following existence lemma for connected standard subspaces $X$ of~$|G|$:

\proclaim Lemma \xxxTSTexistence.
Every standard subspace $Z\sub X$ not containing a circle extends to a \TST\ of~$X$.

\proof
We begin by enumerating the edges in $E(X)\sm E(Z)$. We then go through these edges one by one, considering each for deletion from~$X$. We delete an edge if this does not disconnect the space $X\sm\interior F$, where $F$ is the set of edges already deleted. Having considered every edge in $E(X)\sm E(Z)$, we are left with a standard subspace $T$ that contains~$V(X)$ but contains no circle: this would have an edge in $E(X)\sm E(Z)$, which was considered for deletion and should have been deleted.%
   \COMMENT{}
   The space $T$~is connected (cf.\ Lemma~\xxxConnected),%
   \COMMENT{}
   and hence arc-connected (Lemma~\xxxArcs). By construction, $Z\sub T\sub X$ as desired.\looseness=-1
   \endproof

\noindent
   Unlike in finite graphs, it is considerably harder to construct a \TST\ `from below' (maintaining acyclicity) than, as we did just now, `from above' (maintaining connectedness, and using the non-trivial Lemma~\xxxArcs). A~proof `from below' will be indicated in the proof of Lemma~\xxxTST.

\medbreak

The properties of \TST s resemble those of spanning trees of finite graphs. For example:

\proclaim Lemma \xxxTSTeq. {\rm [\the\refBook]}\COMMENT{}
The following assertions are equivalent for a standard subspace $T$ of~$|G|$ contained in a standard subspace $X$:
\smallskip
   \pitem{i} $T$ is a \TST\ of~$X$.
   \pitem{ii} $T$ is maximally acirclic, that is, it contains no circle but adding any edge of $E(X)\sm E(T)$ creates one.%
   \COMMENT{}
   \pitem{iii} $T$ is minimally connected, that is, it is connected but deleting any  edge of $T$ disconnects it.
\enditem

\noindent
The proof of Lemma~\xxxTSTeq\ needs Lemma~\xxxArcs: we use connectedness in the form of arc-connectedness, but prove only ordinary topological connectedness.

\medbreak

It is not hard to show that the arcs which a \TST\ $T$ of $X$ contains between any two of its points are unique.%
   \COMMENT{}
   Hence every {\it chord\/} $e\in E(X)\sm E(T)$ creates a well-defined {\it fundamental circuit\/} $C_e$ in $T\cup e$, while every edge $f\in E(T)$ lies in a well-defined {\it fundamental cut\/} $D_f$ of~$X$, the set of edges in $X$ between the two arc-components of $T\sm\interior f$.%
   \Footnote{For $X=|G|$ this is a cut of~$G$ in the usual sense.}
   Note that such fundamental cuts are finite by Lemma~\xxxJumpingArc\th(ii), since the two arc-components of $T\sm\interior f$ together contain all the vertices of $X$ and are closed (Corollary~\xxxArcComponents) but disjoint.%
   \COMMENT{}%
   \looseness=-1

For $X=|G|$, \TST s compare with ordinary spanning trees as follows:

\proclaim Lemma \xxxFundamentalBasic. Let $G$ be a locally finite connected graph.
   \pitem{i} The fundamental circuits of any ordinary spanning tree of $G$ are finite, but its fundamental cuts may be infinite.
   \pitem{ii} The fundamental circuits of any \TST\ of $G$ may be infinite, but its fundamental cuts are finite.
   \pitem{iii} The fundamental circuits and cuts of normal spanning trees of $G$ are finite.
   \noproof\enditem


The fact that the fundamental cuts of a \TST\ of $G$%
   \COMMENT{}
   are finite implies by Lemma~\xxxJumpingArc\ that they are in fact bonds.%
   \COMMENT{}

\penalty-500\medbreak

Fundamental circuits and cuts of \TST s are subject to the same duality as for ordinary spanning trees:

\proclaim Lemma \xxxFundamentalDuality.
Let $T$ be a \TST\ of a standard subspace $X$ of~$|G|$, and let $f\in E(T)$ and $e\in E(X)\sm E(T)$. Then $e\in D_f\Leftrightarrow f\in C_e$.\noproof

\section\secTheory. The Topological Cycle Space

The way in which cycles and cuts interact in a graph can be described algebraically: in terms of its `cycle space', its `cut space', and the duality between them.
In this section we show how the cycle space theory of finite graphs extends to locally finite graphs in a way that encompasses infinite circuits. The fact that this can be done, that our topological circuits, cuts and spanning trees interact in the same way as ordinary cycles, cuts and spanning trees do in a finite graph, is by no means clear but rather surprising. For example, there is nothing visibly topological about a finite cut in an infinite graph,%
   \Footnote{With hindsight, of course, Lemma~\xxxConnected\ shows that this impression is wrong.}
   so the fact that the edge sets orthogonal to its finite cuts are precisely its topological circuits and their sums (Theorem~\xxxOrthogonality) comes as a pleasant surprise: it provides a natural answer to a natural question, but not by design---it was not `built into' the definition of a circle.

As it turns out, extending finite cycle space theory in this way is not only possible but also necessary: it is the `topological cycle space' of a locally finite graph, not its usual finitary cycle space, that interacts with its other structural features, such as planarity, in the way we know it from finite graphs. We shall discuss this in some depth in Section~\secApplications. In Section~\secOutlook\ we shall see how our theory integrates with the broader topological context of (singular) homology in more general spaces.

As before, let $G$ be a fixed infinite, connected, locally finite graph. We start by defining the `topological cycle space' $\C(G)$ of~$G$ in analogy to the \hbox{mod-2} (or `unoriented') cycle space of a finite graph:%
   \Footnote{One can also define `oriented' versions $\vCC(G)$ of $\C(G)$, with integer or real coefficients~[\the\refHomologyGraphs].%
   \COMMENT{}
   Some of the theorems listed below have obvious oriented analogues, with very similar proofs. But the differences have not yet been investigated systematically and may well be worth further study; see Sections~\secTheory.3--4 and~\secApplications.3 for some good problems.}
   its elements will be sets of edges (that is to say, maps $E(G)\to\F_2$, or formal sums of edges with coefficients in~$\F_2$) generated from circuits by taking symmetric differences of edge sets. These edge sets, the circuits, and the sums may be infinite.

Let us make this more precise. Let the {\it edge space\/} $\edgespace(G)$ of $G$ be the $\F_2$-vector space of all maps $E(G)\to\F_2$, which we think of as subsets of $E(G)$ with symmetric difference as addition. The {\it vertex space\/}~$\vxspace(G)$ is defined likewise. Call a family $(D_i)_{i\in I}$ of elements of $\edgespace(G)$ \m{thin\/} if no edge lies
in $D_i$ for infinitely many~$i$. Let the \m{(thin) sum\/} $\sum_{i\in I} D_i$ of
this family be the set of all edges that lie in $D_i$ for an odd number of
indices~$i$. Given any subset $\D\sub\edgespace(G)$, the edge sets that are sums of sets in%
   \COMMENT{}
   $\D$ form a subspace of~$\edgespace(G)$. The {\it (topological) cycle space\/} $\cyclespace(G)$ of $G$ is the subspace of $\edgespace(G)$ consisting of the sums of circuits. The {\it cut space\/} $\cutspace(G)$ of $G$ is the subspace of $\edgespace(G)$ consisting of all the cuts in~$G$ and the empty set. (Unlike the circuits, these already form a subspace.) We sometimes call the elements of $\cyclespace(G)$ {\it algebraic cycles\/} in~$G$.

\subsection \secTheory.1 Generating sets

The sums of elements of $\D\sub\edgespace(G)$, and the subspace of $\edgespace(G)$ consisting of all those sums, are said to be {\it generated\/} by~$\D$.%
   \COMMENT{}
   For example, the cycle space of the graph in Figure~\figThreeEnds\ is generated by its central hexagon and its squares, or by the infinite circuit consisting of the fat edges and
all the squares. Bruhn and Georgakopoulos~[\the\refAgelosHenningLA] proved that if $\D$ is thin,%
   \COMMENT{}
   the space it generates is closed under thin sums. As we shall see, this applies to both $\C(G)$ and $\cutspace(G)$.

As in Lemma \xxxFundamentalBasic, the duality between $\C(G)$ and~$\cutspace(G)$---which we \hbox{address} more thoroughly later---is reflected by a switch between topological and ordinary spanning trees:

\vbox{%
\proclaim Theorem \xxxFundamental.~[\the\refBook]
   \pitem {i} Given an ordinary spanning tree of~$G$, its fundamental cuts generate~$\cutspace(G)$, but its fundamental circuits need not generate~$\C(G)$.
\pitem {ii} Given a \TST\ of~$G$, its fundamental circuits generate~$\cyclespace(G)$,%
   \COMMENT{}
   but its fundamental cuts need not generate~$\cutspace(G)$.
   \pitem {iii} Given a \NST\ of~$G$, its fundamental circuits generate~$\cyclespace(G)$, and its fundamental cuts generate~$\cutspace(G)$.
   \enditem

}\noindent
   To prove the first assertion in Theorem \xxxFundamental\th(ii), one shows that a given set $D\in\cyclespace(G)$ equals the sum $\sum C_e$ taken over all chords $e$ of the \TST\ such that $e\in D$. One has to show that these $C_e$ form a thin family (use Lemmas \xxxFundamentalDuality\ and~\xxxFundamentalBasic\th(ii)), but also that $D = \sum C_e$. While for finite $G$ one just notes that $D + \sum C_e$ consists of tree edges and deduces that $D + \sum C_e = \es$ (yielding $D = \sum C_e$), this last implication is now non-trivial: it is not clear that the tree, which by assumption contains no circuit, cannot contain a sum of circuits. The proof in~[\the\refBook] deduces this, i.e.\ that $D$ and $\sum_e C_e$ also agree on tree edges, from Theorem \xxxCBasics\ (i)$\Rightarrow$(iv) below, which in turn is an easy consequence of the jumping arc lemma.%
   \COMMENT{}

For the second assertion in~(ii), recall that the edge set of a \TST\ can miss an infinite cut (as in Figure~\figTSTs). Such a cut will not be a sum of fundamental cuts~$D_f$, because any such sum contains all those tree-edges~$f$.

The proof of assertion~(i) is analogous to that of~(ii), though in the first statement one also has to show that the fundamental cuts do not generate {\em more\/} than~$\cutspace(G)$: that a thin sum of cuts is again a cut. (See the discussion after Theorem~\xxxNSTgenerators.) For the second statement of~(i), recall that the edges of an ordinary spanning tree may contain a circuit, which will not be a sum of fundamental circuits $C_e$ because any such sum contains all those chords~$e$.

Assertion (iii) follows from (i) and~(ii).

\medbreak

   By Lemmas \xxxFundamentalDuality\ and~\xxxFundamentalBasic, the fundamental circuits and cuts of normal spanning trees form thin families.%
   \COMMENT{}
  Hence by the result of~[\the\refAgelosHenningLA] cited earlier, the subspaces they generate in~$\edgespace(G)$ are closed under thin sums. By Theorem~\xxxFundamental\th(iii) these are the spaces $\cyclespace(G)$ and~$\cutspace(G)$:%
   \COMMENT{}

\proclaim Theorem \xxxNSTgenerators.
   \pitem{i} $\cyclespace(G)$ is generated by finite circuits and is closed under infinite thin sums.
   \pitem{ii} $\cutspace(G)$ is generated by finite bonds and is closed under infinite thin sums.

Since we already used the closure statement of~(ii) for cuts in our proof of Theorem~\xxxFundamental, our proof of this part of Theorem~\xxxNSTgenerators\th(ii) has been circular. But one does not need Theorem~\xxxFundamental, or indeed~[\the\refAgelosHenningLA], to prove that a sum of cuts is again a cut, a set of edges across a vertex partition.%
   \COMMENT{}
One way to do this directly is to construct that partition explicitly.%
   \COMMENT{}
   The simplest way, however, is indirect: to show first that the cuts are precisely the edge sets that meet every finite circuit in an even number of edges,%
   \Footnote{This is Theorem~\xxxOrthogonality\th(ii); its proof will be independent of what we are proving here.}
   and then to observe that the set of those edge sets is closed under taking thin sums.%
   \COMMENT{}

But also the closure statement of~(i), for~$\C(G)$, is non-trivial. Indeed, if $\sum D_i$ is a thin sum of sets $D_i = \sum_j C_i^j\in\C(G)$, where the $C_i^j$ are circuits, the combined sum $\sum_{i,j} C_i^j$ need not be thin. By Theorem~\xxxGenerating\ below one can assume for each~$i$ that its own $C_i^j$ are disjoint, in which case the combined family $(C_i^j)$ will be thin (because the family $(D_i)$ is thin by assumption). But Theorem~\xxxGenerating\ is hard, and too big a tool for this proof. As for cuts, the quickest way to prove that $\C(G)$ is closed under thin sums is to use Theorem~\xxxOrthogonality\th(i): that the elements of $\C(G)$ are precisely those sets of edges that meet every finite cut evenly, a class of edge sets that is obviously closed under sums. See also the discussion after Theorem~\xxxOrthogonality.

\medbreak

A bond is {\it atomic\/} if all its edges are incident with some fixed vertex.%
   \COMMENT{}
     Following Tutte, let us call a circuit in $|G|$ {\it peripheral\/} if it contains every edge between its incident vertices (i.e., has no chord in~$G$) and the set of these vertices does not separate~$G$.

\proclaim Theorem \xxxGenerating.
   \pitem{i} Every element of $\cyclespace(G)$ is a disjoint union of circuits.
   \pitem{ii} If $G$ is 3-connected, its peripheral circuits generate~$\cyclespace(G)$.
   \pitem{iii} Every element of $\cutspace(G)$ is a disjoint union of bonds.
   \pitem{iv} The atomic bonds of $G$ generate~$\cutspace(G)$.

\noindent
   Statement (i) was first established in~[\the\refCyclesTwo], with a long and involved proof. As techniques developed, simpler and more elegant proofs were found [\the\refAgelosHotchpotch,\th\the\refCTVella,\th\the\refVellaThesis]. We shall discuss these techniques, and the different proofs they lead to, in Section~\secTechniques. Statement~(ii), which extends a theorem of Tutte for finite graphs (see~[\the\refBook]), is due to Bruhn~[\the\refBruhnPeripheral]. It fails if we allow only finite peripheral circuits or only finite sums. Statements~(iii) and (iv) are easy, with proofs as for finite graphs.%
   \COMMENT{}

Theorem \xxxGenerating\ implies that, as in finite graphs, the circuits in $|G|$ are the minimal non-empty elements of~$\cyclespace(G)$, while the bonds in $G$ are the minimal non-empty elements of~$\cutspace(G)$.

Finally, there is a generating theorem for the cycle space of a finite graph whose extension to $|G|$ and~$\C(G)$ requires an interesting twist. When $G$ is finite, its cycle space is generated by the edge sets of its {\it geodetic\/} cycles: those that contain a shortest path in~$G$ between any two of their vertices~[\the\refBook]. This is still true for locally finite graphs (with `arc' instead of `path'), but only if we measure the length of an arc in the right way: rather than by counting its edges, we have to assign lengths to the edges of~$G$, and then measure the length of an arc by adding up the lengths of its edges. By giving edges shorter lengths if they lie `far out near the ends', it is possible to do this in such a way that the resulting metric on $|G|$ induces its original topology.%
   \Footnote{This will be discussed in detail in Section~\secOutlook.2.}
   And then, no matter how exactly we choose the edge lengths, Georgakopoulos and Spr\"ussel proved that the geodetic circles do what they should:

\proclaim Theorem \xxxGeodetic.~[\the\refAgelosPhilippGeodetic]
Given any edge-length function $\ell\: E(G)\to [0,\infty)$ whose resulting metric on $|G|$ induces the original topology of~$|G|$, the circuits of the geodetic circles in $|G|$ generate~$\C(G)$.

\subsection \secTheory.2 Characterizations of algebraic cycles

There are various equivalent ways to describe the elements of~$\cyclespace(G)$ and of~$\cutspace(G)$, each extending a similar statement about finite graphs. Let us list these now, beginning with~$\cyclespace(G)$.

A~closed topological path in a standard subspace $X$ of~$|G|$, based at a vertex, is a {\it topological Euler tour\/} of $X$ if it traverses every edge in $X$ exactly once. One can show that if $|G|$ admits a topological Euler tour it also has one that traverses every end at most once~[\the\refAgelosInfiniteFleischner]. For arbitrary standard subspaces this is false: consider the closure of two disjoint double rays in the $\Z\times\Z$ grid.%
   \COMMENT{}

Recall that, given a set $D$ of edges, $\overline D$~denotes the closure of the union of all the edges in~$D$, the standard subspace of $|G|$ {\it spanned by~$D$.}

\penalty-1000

\proclaim Theorem \xxxCBasics. 
The following statements are equivalent for sets $D\sub E(G)$:
   \smallskip
   \pitem{i} $D\in\cyclespace(G)$, that is to say, $D$ is a sum of circuits in~$|G|$.
   \pitem{ii} Every component of $\overline D$ admits a topological Euler tour.
   \pitem{iii} Every vertex and every end has even (edge-) degree in~$\overline D$.
   \pitem{iv} $D$ meets every finite cut in an even number of edges.
   \enditem

\noindent
   The equivalence of (i) and~(ii) was proved in~[\the\refCyclesOne] for $D=E(G)$ and extended to arbitrary $D$ by Georgakopoulos~[\the\refAgelosHotchpotch]; we shall meet the techniques needed for the proof in Section~\secTechniques. The equivalence with~(iii) is a deep theorem, due to Bruhn and Stein~[\the\refBruhnSteinEndDeg] for $D=E(G)$ and to Berger and Bruhn~[\the\refBergerBruhnEndDegrees] for arbitrary~$D$. Note that (iii) assumes that end degrees have a parity even when they are infinite. Finding the right way to divide ends of infinite edge-degree into `odd' and `even' was one of the major difficulties to overcome for this characterization of~$\cyclespace(G)$.

The equivalence of (i) with~(iv), again from~[\the\refCyclesOne], is one of the cornerstones of topological cycle space theory: its power lies in the fact that the finitary statement in~(iv) is directly compatible with compactness proofs (see Section~\secTechniques). Its implication (i)$\to$(iv) follows from the jumping arc lemma, applied to the circles whose circuits sum to the given set $D\in\C(G)$. For the converse implication one compares a given set $D$ as in~(iv) with the sum $\sum C_e$ of fundamental circuits of a \TST\ taken over all chords $e\in D$. It is clear that $D$ agrees with this sum on chords. Using Lemmas \xxxFundamentalDuality\ and~\xxxFundamentalBasic\ one proves that they also agree on tree edges.

\subsection \secTheory.3 Cycle-cut orthogonality

Given a set $\FF\sub\edgespace(G)$, let us write $\FF\fin$ for the set of finite elements of~$\FF$. Call two sets $D,F\sub E(G)$ {\it orthogonal\/} if $|D\cap F|$ is finite and even, and put
 $$\FF^\bot := \{\, D\sub E(G) \mid \hbox{$D$ is orthogonal to every $F\in\FF$}\,\}\,.$$%
   \COMMENT{}
 Let us abbreviate $\vxspace:= \vxspace(G)$, \ $\edgespace:= \edgespace(G)$, \ $\C := \C(G)$ and $\cutspace := \cutspace(G)$.

The following theorem describes the duality between the topological \hbox{cycle} and cut space of $G$ in terms of orthogonal sets. Its four assertions show some interesting symmetries: they can be summarized neatly as `all equations containing each of the symbols $=$, $\C$, $\cutspace$, ${}\fin$ and~${}^\bot$ exactly once (and no others)'.%
   \COMMENT{}

\proclaim Theorem~\xxxOrthogonality. 
   \COMMENT{}
   \smallskip
   \pitem{i} $\C = \cutspace\fin^\bot$. 
   \pitem{ii} $\cutspace = \C\fin^\bot$. 
   \pitem{iii} If $G$ is 2-edge-connected, then $\C^\bot = \cutspace\fin$.%
      \COMMENT{}
   \pitem{iv} $\cutspace^\bot = \C\fin$.
\enditem

\noindent
Theorem~\xxxOrthogonality\th(i) is just a reformulation of Theorem~\xxxCBasics\ (i)$\Leftrightarrow$(iv).

Assertion (ii) is the dual (see Section~\secOutlook.4) finitary characterization of the cut space. An interesting feature of its seemingly trivial%
   \COMMENT{}
   inclusion $\cutspace\sub \C\fin^\bot$ is that it requires an application of Theorem \xxxGenerating\th(i). Indeed, any cut $F$ is clearly orthogonal to every finite circuit. But an arbitrary finite element $D$ of~$\cyclespace$ might come as an infinite sum of circuits (finite or infinite), and even if $F$ is orthogonal to every term in the sum it need not be orthogonal to~$D$.%
   \COMMENT{}
   However, the circuits into which $D$ decomposes {\em disjointly\/} by Theorem \xxxGenerating\th(i) will be finite%
       \COMMENT{}
and yield $D$ in a finite sum, which preserves orthogonality.
    The proof of $\cutspace\supe \C\fin^\bot$ is the same as for finite graphs~[\the\refBook]: contracting the edges not in a given set $D\in \C\fin^\bot$ leaves a graph with edge set~$D$, whose finite circuits have to be even by definition of~$D$, because they expand to finite circuits of~$G$.%
   \COMMENT{}
   This contracted bipartite graph defines a partition of~$V(G)$ showing that $D$ is a cut.%
   \COMMENT{}

The inclusion $\C^\bot\supe \cutspace\fin$ in (iii) is equivalent to $\C\sub \cutspace\fin^\bot$ of~(i). The converse inclusion is due to Richter and Vella~[\the\refRichterVella]. Its interesting part is that an edge set orthogonal to every element of~$\C$ (not only to every circuit)%
   \COMMENT{}
   cannot be infinite.%
   \COMMENT{}
   (It will be a cut by~(ii), since in particular it is orthogonal to every finite circuit.)

In~(iv), the weaker inclusion of $\cutspace^\bot\sub \C$ follows from~(i). To see that an edge set orthogonal to all cuts cannot be infinite, assume it is, pick infinitely many independent edges, and obtain a contradiction by extending these to a cut.%
   \COMMENT{}
   The inclusion $\cutspace^\bot\supe \C\fin$ is equivalent to the inclusion $\cutspace\sub \C\fin^\bot$ of~(ii).

\medbreak

Theorem~\xxxOrthogonality\th(ii) has a nice corollary (which we already used in our proofs of Theorems \xxxFundamental\ and~\xxxNSTgenerators): that $\C(G)$ and $\cutspace(G)$ are closed under infinite sums. While this is not immediate from the definition of either $\C(G)$ or~$\cutspace(G)$, the proof that $\cutspace\fin^\bot$ ($=\C(G)$) and $\cyclespace\fin^\bot$ ($=\cutspace(G)$) are closed under thin sums is immediate from the definition of these sets.%
   \COMMENT{}

We remark further that, unlike in finite graphs, there can be edge sets that are orthogonal to all circuits but not to all elements of~$\C$, and edge sets that are orthogonal to all bonds but not to all cuts~[\the\refBook] (even when $G$ is 2-connected).

\medbreak

For finite~$G$ it is well known that $\dim(\edgespace) = \dim(\cyclespace) + \dim(\cutspace)$ (see e.g.~[\the\refBook]), so $\edgespace/\cutspace\simeq\cyclespace$ and $\edgespace/\cyclespace\simeq \cutspace$.%
   \COMMENT{}
   All this is still true for infinite~$G$ and our topological cycle space~$\C$. But it is more instructive to find canonical%
   \COMMENT{}
   such isomorphisms: to define an epimorphism $\sigma\:\edgespace\to\cyclespace$ with kernel~$\cutspace$, and an epimorphism $\tau\:\edgespace\to\cutspace$ with kernel~$\cyclespace$. This can be done with Theorem~\xxxOrthogonality. Indeed, let $T$ be a \NST%
   \Footnote{For each of $\sigma$ and $\tau$ defined below we shall need that $T$ is both an ordinary and a \TST.}%
   \COMMENT{}
   of~$G$, and write $C_e$ and $D_f$ for its fundamental circuits and cuts. By Lemma~\xxxFundamentalBasic, all these are finite. Given $E\sub E(G)$, define
 $$\eqalignno{\sigma(E) &:= \sum\ \{\,C_e : e\in E(G)\sm E(T)\hbox{ and $|E\cap C_e|$ is odd}\,\}\cr
\noalign{\noindent and}
 \tau(E) &:= \sum\ \{\,D_f : f\in E(T)\hbox{ and $|E\cap D_f|$ is odd}\,\}\,.\cr}$$
\bigskip

\proclaim Corollary \xxxIsos.
\pitem{i} The map $\sigma$ is an epimorphism $\edgespace\to\cyclespace$ with kernel~$\cutspace$.
\pitem{ii} The map $\tau$ is an epimorphism $\edgespace\to\cutspace$ with kernel~$\cyclespace$.%
   \COMMENT{}

\smallskip\proof
   (i) To see that $\sigma$ is onto, let $D\in\C$ be given and let $E$ consist of those chords $e$ whose fundamental circuits $C_e$ sum to~$D$.%
   \COMMENT{}
   The kernel of $\sigma$ consists of those sets $E$ that meet every $C_e$ evenly.%
   \COMMENT{}
   These are precisely the sets $E$ that meet every finite $D\in\cyclespace$ evenly,%
   \COMMENT{}
   since these are generated by finitely many~$C_e$.%
   \COMMENT{}
   By Theorem~\xxxOrthogonality\th(ii), these sets $E$ are precisely those in~$\cutspace$.

(ii) is analogous, using Theorem~\xxxOrthogonality\th(i).
   \endproof%
   \COMMENT{}

By Theorem~\xxxOrthogonality, subsets of sets in $\C$ cannot contain an odd cut, and subsets of sets in $\cutspace$ cannot contain an odd circuit. These properties are the only obstructions to being a subset of a set in $\C$ or~$\cutspace$:

\proclaim Theorem \xxxExtensions.~[\the\refBook,~\the\refAgelosHotchpotch]
Let $E\sub E(G)$ be any set of edges.\smallskip
 \pitem{i} $E$ extends to an element of $\cutspace$ if and only if $E$ contains no odd circuit.
 \pitem{ii} $E$ extends to an element of $\C$ if and only if $E$ contains no odd bond.
 \enditem

\smallskip

\proof
We prove the backward implications.%
   \COMMENT{}

(i) Pick a spanning tree in each component of the graph $(V(G),E)$. Adding edges of~$E(G)\sm E$, extend the union of these trees to a spanning tree $T$ of~$G$. Consider the sum
 $$F:= \!\!\sum_{f\in E\cap E(T)} D_f\ \in\>\cutspace$$
of fundamental cuts of~$T$. (The sum is thin by Lemmas \xxxFundamentalDuality\ and~\xxxFundamentalBasic.) We claim that $E\sub F$. Let $e\in E$ be given. If $e\in E(T)$ then $e\in F$, because $D_e$ is the only fundamental cut of $T$ containing~$e$. If $e$ is a chord of~$T$, then its fundamental circuit $C_e$ lies in~$E$, by construction of~$T$. As $C_e$ is finite and hence even (by assumption), Lemma~\xxxFundamentalDuality\ implies that $e\in D_f$ for an odd number of edges $f\in E\cap E(T)$, giving $e\in F$ as desired.

(ii) In every (arc-) component of the standard subspace $\overline{E(G)\sm E}$%
  \COMMENT{}
   pick a \TST\ (Lemma~\xxxTSTexistence). Use Lemma~\xxxTSTexistence\ again to extend their union to a \TST\ $T$ of~$G$; note that the additional edges will be in~$E$, since adding them did not create a circle.%
   \COMMENT{}
   Consider the sum
 $$D:= \!\!\sum_{e\in E\sm E(T)} C_e\ \in\>\C$$
of fundamental circuits of~$T$. (The sum is thin by Lemmas \xxxFundamentalDuality\ and~\xxxFundamentalBasic.) We claim that $E\sub D$. Let $d\in E$ be given. If $d\notin E(T)$ then $d\in D$, because $C_d$ is the only fundamental circuit of $T$ containing~$d$. If $d\in E(T)$, then its fundamental cut $D_d$ lies in~$E$, by construction of~$T$.%
   \COMMENT{}
   As $D_d$ is finite (Lemma~\xxxFundamentalBasic) and hence even (by assumption; recall that fundamental cuts are in fact bonds), Lemma~\xxxFundamentalDuality\ implies that $d\in C_e$ for an odd number of edges $e\in E\sm E(T)$, giving $d\in D$ as desired.
   \endproof

\subsection \secTheory.4 Orthogonal decomposition

If we can write a set $E\sub E(G)$ as a sum $E=D+F$ with $D\in\C$ and $F\in\cutspace$, we call this a {\it decomposition\/} of~$E$.%
   \COMMENT{}
   If the sets $D$ and~$F$ are orthogonal (which for infinite sets is not automatic, their intersection could be infinite or even odd;%
   \COMMENT{}
   cf.\ Lemma~\xxxJumpingArc), we call this decomposition of~$E$ {\it orthogonal\/}.

\proclaim Problem \xxxOrthDec. 
\pitem{i} For which $G$ does every set $E\in\edgespace(G)$ have a decomposition?
\pitem{ii} For which $G$ does every set $E\in\edgespace(G)$ have an orthogonal decomposition?
\pitem{iii} For which $G$ does every set $E\in\edgespace(G)$ have a decomposition $E=D+F$ with either $D\in\C(G)$ and $F\in\cutspace\fin (G)$ or $D\in\C\fin(G)$ and $F\in\cutspace (G)$?
\pitem{iv} When are such decompositions unique?\problem{}

\noindent
   Note that (iii)$\Rightarrow$(ii)$\Rightarrow$(i) for $G$ fixed, by Theorem~\xxxOrthogonality.

\medbreak

If $G$ is finite, a look at the dimensions of the subspaces involved shows that Problem~\xxxOrthDec\ has a clear solution: the answer is `yes' on all counts if $\C(G)\cap\cutspace(G) = \{\es\}$, while otherwise it is `no' on all counts.%
   \COMMENT{}

For $G$ finite, non-empty elements of $\C(G)\cap\cutspace(G)$---edge sets that are both (algebraic) cycles and `cocycles' (cuts)---are called {\it bicycles\/}. If $G$ has no bicycle then, naturally, it is called {\it pedestrian\/}. We shall adopt these definitions verbatim for infinite graphs---Figure~\figPedestrian\ shows a pedestrian graph%
   \Footnote{To see it's pedestrian, consider the leftmost edge of a hypothetical bicycle, and check on which sides of the cut that it defines the vertices near it come to lie.}%
   \COMMENT{}
   from~[\the\refCyCoCy]---but remark that only the finite bicycles are orthogonal to themselves.%
   \COMMENT{}
   In view of Theorem~\xxxOrthogonality\ we would, ideally, hope for a positive answer to Problem~\xxxOrthDec\th(ii) for all $G$ that have no finite bicycles. As we shall see, this is too much to ask---but what exactly is possible is an open problem.

 \figure \figPedestrian.
 An infinite pedestrian graph
 (pedestrian; 1000)

So what can we say when $G$ is infinite? It turns out that neither (ii) nor~(iii) of Problem~\xxxOrthDec\ holds even for pedestrian graphs---except for those that we might call {\em essentially finite\/}: graphs whose spanning trees have only finitely many chords.%
   \Footnote{It is easy to see that this property does not depend on the spanning tree chosen.\COMMENT{}}

\proclaim Theorem \xxxInfOrth.~[\the\refBruhnPersonal,\th\the\refAgelosPersonal]
$G$ satisfies the assertion of {\rm Problem~\xxxOrthDec\th(ii)} if and only if $G$ is pedestrian and essentially finite.

\proof%
   \COMMENT{}
    If $G$ is essentially finite, the union $D$ of all its circuits is the edge set of a finite subgraph~$H$. Its algebraic cycles coincide with those of~$H$, and $\cutspace(G)$ consists of those edge sets whose intersection with $D$ lies in~$\cutspace(H)$. (To see this, consider a spanning tree of $G$ that extends spanning trees of the components of~$H$.)%
   \COMMENT{}
   The assertion now follows easily from the corresponding assertion for~$H$.%
   \COMMENT{}

Before we consider the case that $G$ is not essentially finite, let us outline some general conditions which an infinite set $E\sub E(G)$ must satisfy in order to have an orthogonal decomposition. We shall then show that if $G$ is not essentially finite it will always contain an infinite set of edges violating those conditions.%
   \COMMENT{}
   So let $E$ have an orthogonal decomposition, $E = D+F$ say. Then $D\cap F$ is finite (and even), and $E$ is the disjoint union of a set $E_\C := D\sm (D\cap F)$ and a set $E_\cutspace := F\sm (D\cap F)$ that can be turned into elements of $\C(G)$ or~$\cutspace(G)$, respectively, by adding a finite set (the set $D\cap F$). Moreover, at least one of $E_\C$ and $E_\cutspace$ must be infinite. Hence, for $G$ to be a counterexample to Problem~\xxxOrthDec\th (ii), all it takes is an infinite set $E\sub E(G)$ that has no infinite subset $E_\C$ that can be turned into an algebraic cycle by adding finitely many edges, and no infinite subset $E_\cutspace$ that can be turned into a cut by adding finitely many edges.

Now assume that $G$ is not essentially finite. Let $T$ be a \NST. Then $T$ has an infinite independent set $E$ of chords. Since algebraic cycles induce even degrees, we cannot add finitely many edges to an infinite subset $E'$ of $E$ to obtain an algebraic cycle. But neither can we add finitely many edges to $E'$ to obtain a cut~$F$: since $E'$ is infinite it would take infinitely many fundamental cuts $D_f$ to generate~$F$ (Lemma \xxxFundamentalBasic), and $F$ would differ from $E'$ by at least those (infinitely many) tree edges~$f$.
\endproof  

Two alternative approaches to Problem~\xxxOrthDec\ remain: to settle for~(i), or to prove the assertions from (ii) and~(iii) for specific edge sets $E\sub E(G)$.\problem{in Comment}%
   \COMMENT{}
   The entire set $E=E(G)$, for example, always has an orthogonal decomposition, even into disjoint sets: this is a theorem of Gallai for finite graphs, whose extension to locally finite graphs~[\the\refCyCoCy] has an easy compactness proof of the kind we shall discuss in Section~\secTechniques.%
   \problem{in Comment (Easy BSc)}%
   \COMMENT{}

For singleton sets $\{e\}$, the following is known:

\proclaim Theorem \xxxBicycles.~[\the\refBruhnBicycles,\th\the\refCasteelsRichterBicycles]%
   \COMMENT{}
   Let $e$ be any edge of~$G$.
\pitem{i} Either $e\in E$ for some $E\in\C\cap\cutspace$, or there are $D\in\C\fin$ and $F\in\cutspace\fin$ such that $\{e\} = D+F$, but not both.
\pitem{ii} Either $e\in E$ for some $E\in\C\fin\cap\cutspace\fin$, or there are $D\in\C$ and $F\in\cutspace$ such that $\{e\} = D+F$, but not both.
\enditem

Theorem~\xxxBicycles\th(i) implies a positive answer to Problem~\xxxOrthDec\th(i) for pedestrian graphs~$G$:%
   \COMMENT{}
   the sets from $\C\fin$ and $\cutspace\fin$ needed to generate the singleton sets $\{e\}$ form a thin family~[\the\refBruhnBicycles,~Prop.\th 13], and their sum over all $e\in E$ is exactly~$E$.

\proclaim Corollary~\xxxDec.
If $\C(G)\cap\cutspace(G) = \es$, then every set $E\sub E(G)$ has a decomposition $E=D+F$ with $D\in\C(G)$ and $F\in\cutspace(G)$.
   \COMMENT{}

What can we say if all we assume is that $G$ has no {\em finite\/} bicycle? Then every singleton set $\{e\}$ has a decomposition $\{e\} = D + F$ as in Theorem~\xxxBicycles\th(ii). But unlike for pedestrian graphs,%
   \COMMENT{}
   this cannot always be chosen orthogonal: there can be edges $e$ all whose decompositions $\{e\} = D+F$ are such that both $D$ and~$F$ (and hence also $D\cap F$) are infinite~[\the\refBruhnBicycles]. Moreover, the family of {\em all\/} the $D$ and~$F$ that can be used in such singleton decompositions will not be thin, as soon as $G$ has an infinite bicycle~[\the\refBruhnPersonal].%
   \COMMENT{}
  However we might still hope that for every concrete $E\sub E(G)$ we can choose these decompositions, one for every $e\in E$, so that the sets $D$ and $F$ used form a thin family.\problem{} This would yield the result, stronger than Corollary~\xxxDec, that every edge set in a graph without a finite bicycle has a decomposition.

Finally, uniqueness. For pedestrian graphs, any decomposition $E=D+F$ of an edge set $E$ will be unique: if $E=D'+F'$ was another, then $D-D' = F'-F$ would be a bicycle. For those pedestrian graphs that have orthogonal decompositions, i.e., for the essentially finite ones, we thus have linear%
   \COMMENT{}
   projections $\edgespace\to\C$ and $\edgespace\to\cutspace$ mapping $E$ to $D$ and to~$F$, respectively. There are explicit descriptions of these projections in terms of the number of spanning trees of the graph involved~[\the\refBiggsPotential] (when the graph is finite)---not very tangible but perhaps unavoidably so.

\medbreak

Finally, homology theory appears to suggest the following related question:

\proclaim Problem \xxxBoundaryOperators.
For which $G$ is there a natural `boundary operator' $\partial\:\edgespace(G)\to\vxspace(G)$ with kernel~$\C(G)$, and a natural `coboundary operator' $\delta\:\vxspace(G)\to\edgespace(G)$ with image~$\cutspace\fin(G)$?%
   \COMMENT{}

I suspect that this simplicial attempt to marry our topological cycle space to the usual homological setup is in fact inadequate, and that the way to do this lies not in changing the boundary operators but in restricting the chains. A~singular approach that does capture the topological cycle space as its first homology group, just as in finite graphs, and where some infinite chains are allowed but not every subset (or sum) of edges is a legal 1-chain, has been constructed in~[\the\refHomologySpaces]. It allows for some infinite chains, but not all; see Section~\secOutlook\ for more.

It may also be instructive to study how the answer to the questions in this section, in particular to Problem~\xxxOrthDec,\problem{} changes if we replace our coefficient ring $\F_2$ with $\Z$ or~$\R$.%
   \COMMENT{}
    There will be no bicycles then,%
   \COMMENT{} 
   which simplifies things. But other problems occur. With integer coefficients, for example, 1-dimensional simplicial \hbox{(co-)} chains even of a finite graph need not decompose into an algebraic cycle and a cut (both oriented):%
   \COMMENT{}
   as one can readily check, the triangle with one oriented edge mapped to~1 and the others to~0 is a counterexample.%
   \COMMENT{}

For infinite graphs and real coefficients,%
   \COMMENT{}
   let $\vEE(G)$ consist of only those functions $\psi\: \vE(G)\to \R$ for which \big($\psi(\ve) = \psi(\ev)$ for all~$\ve\in\vE(G)$ and\big)
 $$\sum_{e\in E(G)}\psi^2(e) < \infty\,,$$
 where $\psi^2(e) := \psi(\ve)^2\ (= \psi(\ev)^2)$. For such $\psi,\psi'\in\vEE(G)$ we may define%
   \COMMENT{}
 $$\langle\psi,\psi'\rangle := \sum_{e\in E(G)} \psi(\nve)\psi'(\nve)
       \quad\big(<\infty\big)\,,$$
 where $\nve$ is the {\it natural orientation\/} of~$e$, one of its two orientations that was picked and fixed once and for all when $G$ was first defined. Now $\psi$ and~$\psi'$ can be called {\it orthogonal\/} if $\langle\psi,\psi'\rangle = 0$, and we can study Problem~\xxxOrthDec.\problem{}

\subsection \secTheory.5 Duality

For finite multigraphs $G$ and~$G^*$ one calls $G^*$ a {\it dual\/} of~$G$ if there is a bijection $e\mapsto e^*$ between their edge sets $E := E(G)$ and~$E^* := E(G^*)$ such that every set $F\sub E$ satisfies
 \textno $F\sub\C(G)$ if and only if $F^*\in\cutspace(G^*)$, & (\dag)

\noindent
   where $F^* := \{\,e^*\in E^*\mid e\in F\,\}$. If $G$ is 3-connected, there is at most one such~$G^*$ (up to isomorphism), which is again 3-connected. In this case we think of $G^*$ as being {\em defined\/} by~$(\dag)$, given~$G$. By the finite version of Theorem~\xxxOrthogonality,%
   \COMMENT{}
   condition $(\dag)$~implies the corresponding condition with $\C$ and $\cutspace$ swapped. Hence if $G^*$ is a dual of $G$ then $G$ is a dual of~$G^*$, and for 3-connected $G$ we have $G^{**} = G$.

Can such duals exist also for locally finite infinite graphs? Geometric duals of plane graphs suggest not: the geometric dual of the $\N\times\N$ grid,%
   \Footnote{\dots naively defined just for this example; see Problem~\xxxGeomDuals\ below}
   for example, has a vertex of infinite degree. And indeed, Thomassen~[\the\refCTPlanDualInf] proved%
   \COMMENT{}
   that if $G$ is a locally finite 3-connected graph that (even only) has a locally finite `finitary dual'---a graph $G^*$ for which only the finite sets of edges in $G$ and $G^*$ are required to satisfy~$(\dag)$---then every edge of $G$ must lie in exactly two finite peripheral circuits. In particular, the $\N\times\N$ grid has no such dual.

Hence in order to be able to define duals for all locally finite graphs that should have a dual---that is, for the locally finite planar graphs---we shall have to allow as duals at least some graphs with infinite degrees. And if we want duals to have duals, we have to define duals also for such graphs with infinite degrees. As we shall see in Section~\secOutlook, allowing arbitrary countable graphs with infinite degrees leads to considerable difficulties. However, Thomassen~[\the\refCTPlanDualInf] proved that any graph with (even only) a finitary dual must be {\it finitely separable\/}: every two vertices must be separated by some finite set of edges.

Thus, no matter how we define duals: as soon as their finite edge sets are required to satisfy~$(\dag)$---which is certainly the minimum we shall have to ask---these graphs will have to be finitely separable. On the other hand, the class of finitely separable graphs is not closed under taking finitary duals, and finitary duals are not unique~[\the\refDuality,\th\the\refCTPlanDualInf].%
   \COMMENT{}
   Extending duality to infinite graphs thus seems to be fraught with problems.

In our topological setup, however, the problem has a most elegant solution. Let us say that vertex $v$ {\it dominates\/} a ray~$R$ in $G$ if $G$ contains an infinite $v$--$R$ fan (Fig.~\figFans). It is easy to see that if $v$ dominates one ray of an end it dominates all its rays; we then also say that its {\it dominates\/} that end~$\omega$. Assuming that $G$ is finitely separable, we can then define a space $\tilde G$ similar to~$|G|$, but with the difference that new points at infinity are added only for the undominated ends, while the dominated ends are made to converge to the vertex dominating them.%
   \Footnote{Since $G$ is finitely separable, this vertex is unique. We shall define $\tilde G$ more formally in Section~\secOutlook.}
   Figure~\figIToPDuals\ shows a graph $G$ whose two ends are both dominated by the vertex~$v$. Thus, $\tilde G$~has no added ends, but the double rays `bend round' so that both their tails converge to~$v$.

 \figure \figIToPDuals.
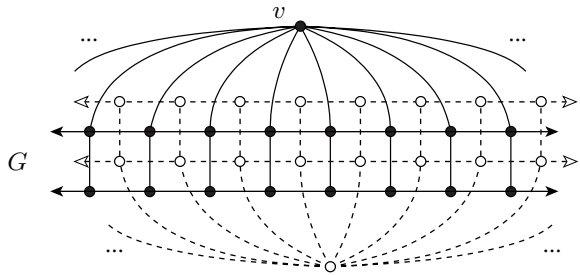
 A self-dual graph with infinite circuits but no added ends
 (selfdual; 1000)

In this space~$\tilde G$, duality works as if by magic. We define $\C(G)$ as the set of circuits (the edge sets of circles) in~$\tilde G$, define $\cutspace(G^*)$ as the set of cuts of~$G^*$ (as before), and call $G^*$ a {\it dual\/} of~$G$ if there is a bijection $e\mapsto e^*$ between their edge sets satisfying $(\dag)$ for all sets $F\sub E(G)$, finite or infinite. If $G$ is 3-connected (and finitely separable), then any dual $G^*$ of $G$ is unique and also finitely separable, and we have $G^{**} = G$ witnessed by the map $e^*\mapsto e=: e^{**}$~[\the\refDuality]. 

\goodbreak

For the purpose only of the next theorem, let us call a (graph-theoretical) spanning tree of a finitely separable graph $G$ {\it acirclic\/} if its closure in $\tilde G$ contains no circle (and similarly for~$G^*$). Note that unless $G$ has no dominated ends, \NST s of $G$ will {\em not\/} be acirclic: they contain a ray from every end, and a ray from a dominated ends that starts at its dominating vertex forms a circle in~$\tilde G$. However, one can show that every finitely separable connected graph has an acirclic spanning tree~[\the\refDuality].%
   \COMMENT{}

\proclaim Theorem~\xxxDualTrees.~[\the\refDuality] 
Let $G=(V,E)$ and $G^*=(V^*,E^*)$ be connected finitely separable multigraphs, and let $^*\colon E\to E^*$
be a bijection. Then the following two assertions are equivalent:
\pitem{i} $G$ and $G^*$ are duals of each other, and this is witnessed by the
 map~$^*$ and its inverse.
\pitem{ii} Given a set $F\subseteq E$, the graph $(V,F)$ is an acirclic
spanning tree of~$G$ if and only if $(V^*,E^* \sm F^*)$ is an
acirclic spanning tree\ of~$G^*$.%
   \COMMENT{}
   \enditem

\noindent
   (For~(ii), define $F^* := \{\,e^*\in E^*\mid e\in F\,\}$ as before.) 

\medbreak

Can we generalize Theorem~\xxxDualTrees\ to \TST s that are not necessarily graph-theoretical trees? A~positive answer to this question would be in line with the philosophy of our theory, but the answer is not obvious. Once we use ends as `connectors' (as we do in a \TST), we should treat them on a par with edges: we should have a bijection between the ends of $\tilde G$ and those of~$\tilde{\hbox{${G^*}$}}$, and use each end from this common set of ends in precisely one of the two trees. But no matter how one defines the `ends of~$\tilde G$', whether one takes just the undominated ends or the undominated ends plus the dominating vertices (i.e., the set of all limit points of rays), there need not be a bijection between the ends of $\tilde G$ and those of~$\tilde{\hbox{${G^*}$}}$: the $3\times\Z$ grid, for example, is a 3-connected graph with two ends, whose dual%
   \COMMENT{}
   has no undominated end and only one dominating vertex.

\penalty-5000

Surprisingly, things work better if we, initially, continue to consider the ends of $G$ and $G^*$ in the original spaces $|G|$ and~$|G^*|$. Bruhn and Stein~[\the\refBruhnSteinEndDuality] have shown that, if $G$ is 2-connected, there is a homeomorphism between the subspaces $\Omega$ of~$|G|$ and $\Omega^*$ of~$|G^*|$ that is {\it compatible\/} with the duality map $e\mapsto e^*$ from $E$ to~$E^*$ in the sense that, for 
every set $F\subseteq E$, an end $\omega\in\Omega$ is an 
accumulation point of $\bigcup F$ in~$|G|$ if and only if 
$\omega^*$ is an accumulation point of $\bigcup F^*$ in~$|G^*|$.%
   \COMMENT{}
   The question, then, could be roughly as follows: is it true that no matter how we form a spanning tree in $\tilde G$ from some of the edges and ends of~$G$ (taken in~$|G|$), the other edges and ends (taken in~$|G^*|$) will form a spanning tree in~$\tilde{\hbox{${G^*}$}}$? And that, conversely, if this is the case for two graphs $G$ and $G^*$ then these form a dual pair? 

Figure~\figSerpents\ shows a dual pair of graphs%
   \COMMENT{}
   illustrating this. The grey tree in $G^*$ needs the `shared' end in order to be connected, and adding the end to it does not create a circle. The dual black tree, however, is already connected (without the end), and adding the end there would create a circle. Hence there is a unique way of assigning the end to one of the two trees, so that both are spanning and acirclic in their respective graphs.

 \figure \figSerpents.
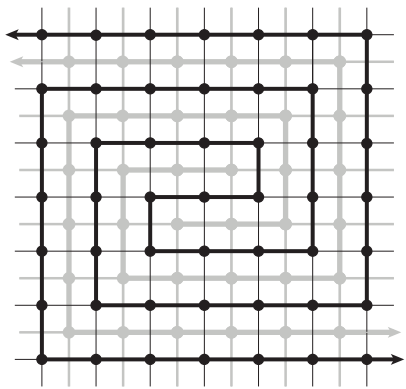
 Dual `spanning trees': the grey tree contai\rlap{ns} the shared end, the black tree does not%
   \COMMENT{}
 (serpents; 1000)

Let us make this precise in greater generality. Given a set $F\sub E$ of edges and a set $\Psi\sub\Omega$ of ends, let us say that the subspace $\bigcup F\cup\Psi$ of $|G|$ is {\it formed by\/} $F$ and~$\Psi$. Given a subspace $X$ of~$|G|$, let $\tilde X$ denote the quotient space%
   \COMMENT{}
   of $X$ obtained by identifying every vertex in $X$ with all the ends in $X$ that it dominates (in~$G$). (For $X = G$ this yields precisely the definition of $\tilde G$ given informally earlier.) For example, if $G$ is the fan shown on the left in Figure~\figFans, we might take as $X$ the space formed by the set $F$ of fat edges and the unique end of~$G$. Then $\tilde X$ is a circle. However, if we take as $X$ the space formed only by the set $F$ of fat edges, without the end, then $\tilde X$ is homeomorphic to the half-open interval~$[0,1)$. Call $\tilde X$ a {\it spanning tree\/} of~$\tilde G$ if it contains all the vertices of~$\tilde G$, is arc-connected, and contains no circle.%
   \COMMENT{}

\goodbreak

\proclaim Theorem~\xxxDualTopTrees.~[\the\refTreeEnds]
Let $G=(V,E,\Omega)$ and $G^*=(V^*,E^*,\Omega^*)$ be 2-connected finitely separable multigraphs, and let $^*\colon E\to E^*$ and  $^*\colon \Omega\to \Omega^*$
be compatible\?{Check whether (ii) implies compatibility, in which case we can move this assumption to~(i).} bijections. Then the following two assertions are equivalent:
   \pitem{i} $G$ and $G^*$ are duals of each other, and this is witnessed by the maps $E\to E^*$ and $\Omega\to\Omega^*$ and their inverses.
   \pitem{ii} Whenever the subspace $X$ of $|G|$ formed by subsets $F\sub E$ and $\Psi\sub\Omega$ is such that $\tilde X$ is a spanning tree of~$\tilde G$, the sets $E^*\sm F^*$ and $\Omega^*\sm\Psi^*$ form a subspace $Y$ of $|G^*|$ such that $\tilde Y$ is a spanning tree of~$\tilde{\hbox{${G^*}$}}$.
   \enditem

\noindent
   Here, as before, $F^* := \{\,e^*\in E^*\mid e\in F\,\}$ and ${\Omega^* := \{\,\omega^*\in \Omega^*\mid \omega\in \Omega\,\}}$. The obvious implications for infinite matroids which Theorem~\xxxOrthogonality, Theorem~\xxxDualTrees\ and~\xxxDualTopTrees\ suggest will be a topic in Section~\secOutlook.

\medbreak

For finite graphs, the theory of algebraic duality outlined above corresponds to a theory of {\it geometric duality\/}: two connected graphs are duals of each other if and only if they can be drawn in the plane so that precisely the corresponding edges cross and vertices correspond bijectively to faces of the dual.%
   \Footnote{In particular, graphs with a dual are planar. This is Whitney's theorem, which extends to locally finite finitely separable graphs if and only if duality is defined as we did: with $\C$ as the topological cycle space of $G$ based on~$\tilde G$. See Section~\secApplications.}

For our finitely separable graphs, this is not so obvious: how, for instance, should we draw the two dual graphs of Figure~\figIToPDuals\ so that their rays converge as they do in~$\tilde G$ and~$\tilde{\hbox{${G^*}$}}$? In particular, where in such a drawing can we put the two vertices of infinite degree, the limit points of those rays? For now, setting up geometric duality in a satisfactory way remains an unsolved problem:

\proclaim Problem \xxxGeomDuals. Define geometric duals for finitely separable graphs compatibly with their algebraic duals.%
   \COMMENT{}

\noindent
   Perhaps the solution to Problem~\xxxGeomDuals\ lies in an approach similar to that for Theorem~\xxxDualTopTrees: to apply duality `before identification' rather than after, that is to say, to look for geometrically dual embeddings of $|G|$ and~$|G^*|$ rather than of $\tilde G$ and~$\tilde{G^*}$.%
   \COMMENT{}
   The drawback of such an approach is that $|G|$ and~$|G^*|$ need not be compact as soon as we refine their topologies enough to make their embeddings in the plane continuous.%
   \Footnote{Since $|G|$ fails to be Hausdorff as soon as $G$ has a dominated end, we need to refine our current topology, known as \fnVTop, to a topology called~\fnMTop; see Section~\secOutlook.1--2. Under \fnMTop, $|G|$~is compact only if $G$ is locally finite (in which case \fnMTop\ agrees with~\fnVTop).}%
   \COMMENT{}
   This, then, raises issues such as how exactly to define a face.

\section\secTechniques. Proof Techniques

Once more, let $G$ be an infinite, locally finite, connected graph.%
   \COMMENT{}
   The aim of this section is to bring the reader new to the field up to date with its collective memory of techniques. Although each of the deeper proofs of existing theorems naturally has its own difficulties and ways of overcoming them, there is by now a body of basic approaches that can be described---as well as some common pitfalls the novice might like to hear about before plummeting their depths themselves.

We shall concentrate on proofs that establish the existence of infinite substructures of~$|G|$ such as arcs, circles, or \TST s. Such a proof usually has two parts: the construction of this structure, usually by a limit process, and the proof that the object thus constructed does what it is intended to do.%
   \Footnote{For example, we might be constructing a {\it Hamilton circle\/}, a circle in $|G|$ passing through all the vertices. The construction part might define a map $\sigma\: [0,1]\to |G|$ that is, by definition, injective at vertices and inner points of edges (except that $\sigma(0) = \sigma(1)$), while we might still have to show that $\sigma$ is continuous and injective at ends.}
   We shall only discuss the construction part. If the construction is done well, it can happen that the proof of correctness becomes easy. This does not always happen, and if it does, it may simply mean that the problem is intricate enough that `blind' constructions will fail: that any difficulties that may arise in the proof have to be anticipated by the construction. Still, the basic techniques that are common to many proofs tend to lie at the construction level, while difficulties arising at proof level tend to be individual to the problem. Our focus on the construction part of proofs reflects this phenomenon; it is not a deliberate restriction.

\subsection \secTechniques.1 The direct approach

The direct, or naive, approach to the construction of standard subspaces with certain desired properties, such as an arc, a circle or a \TST, is to obtain them as the union of nested finite subgraphs defined explicitly, or as the intersection of a nested sequence of explicitly defined superspaces. For example, we might try to construct a spanning tree as a union of finite trees $T_0\sub T_1\sub\dots$ eventually covering all the vertices of~$G$. The union $T$ of these~$T_n$ will be a tree, but its closure in $|G|$ may well contain a circle: unlike finite cycles, this will arise only at the limit step, and it may not be clear how to choose the $T_n$ so that this does not happen.

If we want this naive approach to work, we have to express the statement to be proved in a finitary way---for example, as statements about finite sets of edges. The characterization of topological connected\-ness in Lemma~\xxxConnected\ is a case in point; the characterization of cycle space elements in Theorem~\xxxOrthogonality\th(i) is another. Unlike the circle in our example, finitary properties do not appear or disappear at limits, which makes it easy to prove such statements by induction or an application of Zorn's lemma.

For instance, we have already seen the construction of a \TST\ `from above', in the proof of Lemma \xxxTSTexistence. In that proof we used the fact that connectedness of a standard subspace is a finitary property (Lemma~\xxxConnected), and that \TST s are edge-minimal connected standard subspaces of $|G|$ containing all its vertices. The proof of the following lemma, which we shall need in Section~\secApplications, is very similar:

\proclaim Lemma~\xxxkCon. Let $X$ be a standard subspace of~$G$.
\pitem{i} If $X$ is $k$-edge-connected, we can delete edges from $X$ to obtain a standard subspace $Y$ that is edge-minimally $k$-edge-connected.
\pitem{ii} If $X$ is $k$-vertex-connected, we can delete edges from $X$ to obtain a standard subspace $Y$ that is edge-minimally $k$-vertex-connected.
\enditem

\proof (i) To obtain $Y$ from~$X$, we go through the edges of $X$ in turn, deleting the edge under consideration if and only if this does not spoil the $k$-edge-connectedness of the current space. To show that $Y$ is still $k$-edge-connected (in which case it will clearly be minimally so), we consider a hypothetical set $F$ of fewer than $k$ edges whose removal disconnects~$Y$. By Lemma~\xxxConnected, the space $Y\sm\interior F$ meets both sides of some finite cut of $G$ in which it has no edge. But $X\sm\interior F$~has an edge in this cut, because it is connected. All its edges in the cut, however, were deleted in the construction of~$Y$, and the edge deleted last should not have been deleted.

(ii) The proof of the vertex-connectivity case is analogous: we again delete edges one by one, but instead of an edge set~$F$ and the space~$Y\sm\interior F$ we consider a separating set $U$ of fewer than $k$ vertices and the standard space obtained from $Y$ by deleting $U$ and its incident edges.%
   \COMMENT{}
\endproof

The analogue of Lemma~\xxxkCon\ for $k$-connected subspaces---those that cannot be made disconnected by the deletion of fewer than~$k$ vertices, edges or ends---cannot be proved in the same way. For arbitrary substandard subspaces it is false (consider the double ladder, without the ends, for $k=2$), but for standard subspaces it is an open problem.%
   \Footnote{We shall discuss this problem in Section~\secApplications.1.}
   For example, if we delete all the broken edges in the graph of Figure~\figMinK\ in $\omega$~steps, the space will lose its 3-connectedness (since deleting the two ends will then disconnect it), but does so only at the limit step.

 \figure \figMinK.
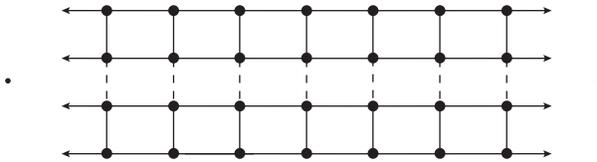
 The subspace obtained by deleting all the broken edges is no longer 3-connected.
 (MinK; 1000)

\medbreak

As we shall see in Section~\secApplications, many applications of the object $|G|$ have the following form: there is an existence statement about finite graphs that fails for infinite graphs---for example, a theorem asserting the existence of a Hamilton cycle---but a topological analogue works in~$|G|$. In our example, it might be that the same conditions (such as being planar and 4-connected) that force a Hamilton cycle in a finite graph force a Hamilton circle in~$|G|$. Rather than proving the infinite theorem from first principles, we might try to use the finite result for the proof of its infinite analogue. We would then want to express $|G|$ as a limit of finite graphs $G_0, G_1,\dots$ that satisfy the assumptions (such as planarity and 4-connectedness) if $G$ does, use the finite theorem to find the desired structures in the~$G_n$, and then take a limit of these to obtain the analogous structure in~$|G|$.

Experience has shown that defining these $G_n$ as induced subgraphs of~$G$, e.g.\ those on its first $n$ vertices, is {\em not\/} often a good approach: a badly chosen $G_n\sub G$ is too oblivious of how it lies inside~$G$. (There might, for example, exist $G_n$-paths%
   \Footnote{See~[\the\refBook] for the notion of an $H$-path in $G$ for a subgraph $H\sub G$.}
   in $G$ between vertices of $G_n$ that are easily separated in~$G_n$.) However, there is a standard approach using finite `nested' minors rather than subgraphs as~$G_n$, which we define now.

Let $v_0,v_1,\dots$ be an enumeration of the vertices of~$G$. For each $n\in\N$ let $S_n := \{v_0,\dots,v_n\}$, and write $G_n$ for the minor of $G$ obtained by contracting each component of $G-S_n$ to a vertex, deleting any loops that arise in the contraction but keeping multiple edges. The vertices of $G_n$ contracted from components of $G-S_n$ will be called the {\it dummy vertices\/} of~$G_n$.

Note that every cut of $G_n$ is also a cut of~$G$. Conversely, a cut of $G$ is a cut of every $G_n$ that contains all its edges. This is an important feature, which distinguishes these minors from a similar sequence of finite subgraphs exhausting~$G$. A~topologically connected standard subspace of~$|G|$, for example, will induce a connected subgraph in each~$G_n$ (cf.\ Lemma~\xxxConnected), an element of the cycle space of~$G$ will induce a cycle space element of $G_n$ (cf.\ Theorem \xxxCBasics.), and so on. We shall appeal to these properties again in our discussion of homology in Section~\secOutlook.3.

As a first application of these $G_n$ let us show that they allow us to construct a \TST\ directly `from below', using the direct approach.%
   \Footnote{The proof of Lemma~\xxxTST\ can be generalized to a proof `from below' of the full statement of Lemma~\xxxTSTexistence\ for the existence of \TST s in standard subspaces~[\the\refTreeEnds].}

\proclaim Lemma \xxxTST. $G$ has a \TST.

\removelastskip\medskip\noindent{\bf Proof \it(by explicit construction from below).}\nl
Pick spanning trees $T_n$ of $G_n$ recursively so that $E(T_{n+1})\cap E(G_n) = E(T_n)$.%
   \COMMENT{}%
   \Footnote{There is a canonical%
   \COMMENT{}
   way of doing this. Unless $v_{n+1}$ forms a component of $G-S_n$ by itself (in which case $T_{n+1} = T_n$), the graph $G_{n+1}$~is obtained from $G_n$ by expanding the dummy vertex of $G_n$ contracted from the component containing~$v_{n+1}$ to a star with centre $v_{n+1}$ and some new dummy vertices of~$G_{n+1}$ as leaves. Note that while the vertices of this star are unique, its edges are not, since $G_{n+1}$ has multiple edges.%
   \COMMENT{}
   We  obtain $T_{n+1}$ from $T_n$ by adding the edges of this star. In other words, we extend $T_n$ to $T_{n+1}$ in any way that does not add an edge at a vertex in~$S_n$.}
   We claim that $T:=\overline{\bigcup_n E(T_n)}$ is a \TST\ of~$G$.%
   \COMMENT{}

To show that $T$ is connected (and hence arc-connected, by Lemma~\xxxArcs), we have to check that every finite cut of $G$ contains an edge from~$T$. It does, because it is also a cut of some~$G_n$,%
   \COMMENT{}
   and $T_n$ has an edge in this cut.

Now suppose that $T$ contains a circle, $C$~say. Pick vertices $u,v\in C$, and choose $n$ large enough that $u,v\in S_n$. Let $F$ be a fundamental cut of $T_n$ separating $u$ from~$v$. Since this is also a cut of~$G$ and $C$ meets both its sides, $C$~has an edge in~$F$ (Lemma~\xxxConnected). Deleting this edge from $C$ leaves an arc, which still has an edge in~$F$ (by the same argument). Hence $C\sub T$ has two edges in~$F$, contradicting the fact that $E(T)\cap E(G_n) = E(T_n)$ and $T_n$ has only one edge in~$F$.%
   \COMMENT{}
   \endproof

Finally, there is a type of construction, very common in finite graph theory, which is rather naive in the context of~$|G|$ but should not go unmentioned: the construction of paths or arcs by `greedily moving along'. As an example, consider the proof of Theorem~\xxxGenerating\th(i), that every edge set $D\in\cyclespace(G)$ is a disjoint union of circuits. Let us just try to find the first circuit. If $G$ is finite, we start from any edge in~$D$ and `greedily move along' further edges of~$D$: since $D$ has an even number of edges at every vertex, we cannot get stuck until we reach our starting vertex again. We stop, however, when we hit the first previously visited vertex, having found our first circuit. We then delete this and repeat, eventually decomposing $D$ into circuits.

When $G$ and $D$ are infinite, this approach can fail quite spectacularly. For example, $D$~might be the wild circuit of Figure~\figWildCircle. By `moving along' its edges we would run into an end, from which there would be no escape, not even by an inverse ray.
The linear order structure of the double rays on the wild circle indicates that the only way to obtain it as a limit is to `approach all its points at once'---which is exactly what we shall do using our sequence of~$G_n$.

\subsection \secTechniques.2 The use of compactness

A typical compactness proof, such as of the \Erdos\th-\th De\th Bruijn theorem on colouring, is a way of making consistent choices for all the finite substructures of an infinite structure, and then deducing from the fact that all the finite substructures have a certain property that so does the entire structure.

In our context, only the simplest results have compactness proofs quite like this. But we often use compactness for the construction part of a proof, including the more difficult ones: to define a limit object (e.g., a subspace~$X$) in such a way that all of certain finite substructures (e.g., the finite cuts of~$G$) relate to the limit object in a certain desired way (contain an edge of~$X$). This may be less than we actually want to prove (e.g.\ that $X$ is arc-connected)---but it may be a good start.

The best examples for this are infinite analogues of theorems asserting for a finite graph the existence of a path or cycle that can be found by a process of `focussing in' on some part of a graph: in such cases, the same process may converge to a limit object in~$|G|$ with similar properties. Bruhn's proof~[\the\refBruhnPeripheral] of Theorem~\xxxGenerating\th(ii), for example, works in principle like Tutte's own proof of his finite theorem. Or think of the `focussing in' on a component in Thomassen's proof that every finite 3-connected graph has an edge whose contraction preserves 3-connectedness. (See~[\the\refBook] for both.)\problem{find and prove an infinite analogue!}

To illustrate how such a compactness proof typically works, or might fail, let us take another look at the proof of Theorem~\xxxGenerating\th(i). We were trying to find a single circuit $C$ in a given set $D\in\cyclespace(G)$ (which might be the wild circle, so that we cannot find $C$ by `moving along'). Our approach now is to construct $C$ as a limit of circuits $C_n$ in the~$G_n$. Since $D$ meets every cut of~$G_n$ evenly (because it is also a finite cut of~$G$), we have $D_n := D\cap E(G_n)\in\cyclespace(G_n)$. By the trivial finite version of Theorem~\xxxGenerating\th(i), we find a circuit $C_n\sub D_n$. Will these~$C_n$, or at least a subsequence, tend to a limit $C\sub D$ that is a circle in~$|G|$?

Not necessarily. For a start, the limit might exist%
   \Footnote{We have not made this precise yet.}
   but be empty. Indeed, depending on the structure of~$D$,%
   \COMMENT{}
   we might have chosen the $C_n$ so badly that no edge lies eventually in~$C_n$.%
   \Footnote{For example, let $G$ be the plane $\Z\times\Z$ grid, and choose as $C_n$ the boundary circuit of the outer face of~$G_n$.}%
   \COMMENT{}
   But there is a simple way to prevent this: by insisting that every $C_n$ contain some fixed edge $e\in D$.

Assuming this, let us use the infinity lemma to construct a limit of these~$C_n$. To apply the lemma we need to define `predecessors', which is easy: let the $n$th finite set in the setup of the lemma contain not only $C_n$ but also all edge sets of the form $C_n^m := C_m\cap E(G_n)$ with $m > n$. Since the cuts of $G_n$ are also cuts of~$G_m$, these sets will meet every such cut evenly and hence lie in~$\C(G_n)$. We may now define the predecessor relation for the infinity lemma as inclusion (of edge sets), and obtain the desired limit $C$ as a union of certain nested sets $C^m_n\in \cyclespace(G_n)$ found by the lemma, one for each~$n$.%
   \COMMENT{}
   Then~$C$, too, meets every finite cut of $G$ evenly, because that cut is also a cut of some~$G_n$ and the selected $C^m_n$ meets this cut evenly. Hence $C\in\cyclespace(G)$, by Theorem~\xxxOrthogonality\th(i).

But is $C$ a circuit? Unfortunately, it need not be: we cannot in general prevent the $C^m_n$ from traversing dummy vertices more than once,%
   \COMMENT{}
   in which case the closure $\overline C$ of their union can be a union of circles that are vertex-disjoint but meet at ends. Figure~\figBadLimit, where $D = E(G)$,%
   \COMMENT{}
   indicates an example where this happens: the sequence of $C_n$ there is formed by alternately extending them in two different places, so that both extensions tend to the unique end, yielding a limit in which that end has degree~4.%
   \COMMENT{}

 \figure \figBadLimit.
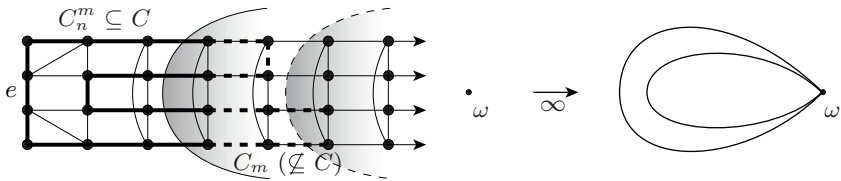
 Finite cycles tending to a double circle. \rlap{The}\nl shaded regions contract to dummy verti\rlap{ces.}
 (BadLimitNEWLettering; 1000)

\goodbreak
 
So what did we achieve? As the limit of our circuits $C_n$ we found just another algebraic cycle in~$\cyclespace(G)$, but not a circuit. For our present problem this seems like next to nothing: after all, we started with an algebraic cycle $D\in\cyclespace(G)$, looking for a circuit $C\sub D$. But in general, the method does have its uses: it shows that from algebraic cycles $C_n\in\cyclespace(G_n)$---we never used that they are actually circuits---we can always obtain an algebraic cycle $C\sub\bigcup_n C_n$ by compactness, even if the $C_n$ were not induced by any previously given $D\in\cyclespace(G)$. And we can say a little more about $C$ too: there is a subsequence of the~$C_n$ (consisting of those $C_m$ that gave rise to one of the $C^m_n$ in our nested sequence picked by the infinity lemma) such that every edge of $C$ lies eventually on every $C_n$ in this subsequence.

Moreover, the set $C\sub\cyclespace(G)$ we found (now using that the $C_n$ are circuits) does have one advantage over~$D$: the standard subspace $\overline C\sm\interior e$ is connected. Indeed, every $C_m\sm\{e\}$ spans a connected subgraph in~$G_m$, and hence so does $C^m_n\sm\{e\}$ in~$G_n$. The connectedness of $\overline C\sm\interior e$ now follows from $C^m_n\sub C$ and Lemma~\xxxConnected.%
   \COMMENT{}
   We can therefore use Lemma~\xxxArcs\ to find an arc in $\overline C\sm\interior e$ between the endvertices of~$e$, which together with $e$ will form the desired single circuit in~$D$ containing~$e$---solving our problem.%
   \Footnote{The proof of Theorem~\xxxGenerating\th(i) follows easily now: we just collect circuits $C\sub D$ inductively and delete them, to eventually decompose~$D$, starting each new circuit with the next uncovered edge of $D$ in some fixed enumeration.}

\subsection \secTechniques.3 Constructing arcs directly

In our compactness proof discussed above, we obtained as a limit of finite cycles $C_n\sub G_n$, all containing some fixed edge~$e$, an edge set $C\in\cyclespace(G)$ with $\overline C\sm\interior e$ connected, and we had to invoke Lemma~\xxxArcs\ in order to find inside $\overline C\sm\interior e$ an arc between the endpoints of~$e$. Could we have obtained such an arc directly by using compactness in a more subtle way, e.g.\ by choosing the edge sets to which we applied the infinity lemma more carefully?

\proclaim Problem \xxxDirectArc.
Given two vertices $x,y\in G$ and, for all $n$ large enough, an $x$--$y$ path $P_n\sub G_n$, construct by a direct compactness argument an $x$--$y$ arc in $|G|$ with edges in~$\bigcup_n E(P_n)$.

This problem is still unsolved. What we can do is construct a topological $x$--$y$ path in $|G|$ by direct compactness. This path will traverse vertices at most once, but it may be non-injective at ends (as in Figure~\figBadLimit). In order to obtain the desired arc, we therefore still have to invoke a theorem from general topology: that the image of any $x$--$y$ path in a Hausdorff space contains an $x$--$y$ arc.\looseness=-1

Before we apply this technique to obtain such a partial solution to Problem~\xxxDirectArc, let me illustrate it by an exercise. Let $T$ be the infinite binary tree. The task is to construct a topological path $\sigma\: [0,1]\to |T|$ that starts and ends at the root and traverses every edge exactly twice, once in each direction. For $k\in\N$, let $T^k$ denote the finite subtree of $T$ consisting of its first $k$ levels. We shall obtain $\sigma$ as a limit of analogous paths $\sigma_k$ around these~$T^k$, as follows.

Let $\sigma_0$ be the constant map from $[0,1]$ to the 1-vertex graph~$T^0$. Assume inductively that we have defined for $i=0,\dots,k$ paths $\sigma_i\: [0,1]\to T^i$ such that the inverse image under $\sigma_i$ of each leaf of $T^i$ is a non-trivial interval, that the inverse image of any other point of $T^i$ is a single point, and that $\sigma_{i+1}$ agrees with $\sigma_i$ on all points of $[0,1]$ that $\sigma_i$ maps to points of $T^i$ other than leaves. (Thinking of $\sigma_i$ as a journey with $[0,1]$ measuring `time', we thus ask that $\sigma_i$ should pause for a while at every leaf of~$T^i$, but not elsewhere.) In order to obtain $\sigma_{k+1}$ from~$\sigma_k$, we just expand every pause at a leaf $t$ of $T^k$ to a little path $tt'tt''t$ in~$T^{k+1}$ from $t$ to its two successors $t',t''$ and back, pausing again at each of $t'$ and~$t''$.%
   \COMMENT{}

To define our desired limit $\sigma$ of these~$\sigma_k$, we first look at all those $x\in [0,1]$ for which the $\sigma_k$ eventually agree, and let $\sigma(x) := \sigma_k(x)$ for these~$x$ (with $k$ large enough). For any other~$x$, every $\sigma_k(x)$ will be a leaf of~$T^k$, and as $k\to\infty$ these $\sigma_k(x)$ form a ray in~$T$; we let $\sigma(x)$ be the end of that ray. It is straightforward from the definition of the topology of~$|T|$, but still instructive, to check that this map $\sigma\: [0,1]\to |T|$ is continuous at these latter points~$x$ mapped to ends. Indeed, every basic open neighbourhood of an end~$\omega = \sigma(x)$ has the following form: for some vertex $t$ on the ray $R$ in $\omega$ starting at the root, it is the set of all points in $|T|$ strictly above~$t$. Let $t'$ be the upper neighbour of $t$ on~$R$. If $k$ is the level of~$t'$ in~$T$, then $t'$ is a leaf of~$T^k$, and all the points in the subinterval of $[0,1]$ that $\sigma_k$ maps to this leaf $t'$ will be mapped by all $\sigma_\ell$ with $\ell > k$ to points of $G$ that still lie above~$t$. Hence $\sigma$ will map all such points to the closure of this set (but not to~$t$), i.e.\ to the open neighbourhood of $\omega$ we started with.

This technique for constructing a topological path directly as a limit of maps can be combined with compactness, as follows. Consider Problem~\xxxDirectArc. Every $x$--$y$ walk in~$G_m$, including the given path~$P_m$, induces for every $n<m$ an $x$--$y$ walk in~$G_n$ consisting of those of its edges that already lie in~$G_n$, traversed in the same order and direction. In order to obtain an $x$--$y$ arc in~$|G|$, we first use the infinity lemma to select from these walks a sequence $(W_n)_{n\in\N}$ of $x$--$y$ walks $W_n$ each inducing the previous. We then parametrize these walks in turn for $n=1,2,\dots$, as in our binary tree example: by topological paths~$\sigma_n$ `pausing' at every dummy vertex of~$G_n$ (but only there),%
   \COMMENT{}
   and so as to agree with $\sigma_{n-1}$ at all points that $\sigma_{n-1}$ does not map to a dummy vertex.%
   \COMMENT{}
   The limit of these~$\sigma_n$ will be a topological $x$--$y$ path $\sigma$ in $|G|$, defined as earlier, since sequences of dummy vertices of the $G_n$ ($n\to\infty$) that correspond to nested components of $G-S_n$ define a unique end, and these components form a neighbourhood basis for that end.

Note that if we insist that the lengths of pauses tend to zero as $n$ tends to infinity, the inverse image of an end under the limit path $\sigma$ will be a totally disconnected subset of~$[0,1]$.%
   \COMMENT{}
   But I do not know how to set up the infinity lemma in such a way that $\sigma$ will be injective at ends: we have found an $x$--$y$ path in~$|G|$, but not the desired arc.

The construction of a topological Euler tour in~$|G|$ for the proof of Theorem \xxxCBasics\th(ii) works exactly like this, too. Georgakopoulos has bundled the essence of this approach into an out-of-the-box tool [\the\refAgelosHotchpotch, Lemma~1.4] that provides a well-defined limit path for a suitable subsequence of a given sequence of $x$--$y$ paths. It deals once and for all with the construction and continuity of the limit path, and can be used in many situations.

\section\secApplications. Applications

A particularly satisfying feature of the topological approach to infinite graphs via ends is its applicability. Considering arcs and circles in $|G|$ instead of paths and cycles in $G$ enables us to extend classical theorems from finite graph theory whose naive extension to infinite graphs fails.%

In this section we give an overview of three areas in which this phenomenon is particularly striking: an emerging `extremal' theory of infinite graphs; applications of the cycle space (eg.\ to characterizations of planarity), and infinite electrical networks. Let $G$ be a locally finite, connected infinite graph.

\subsection\secApplications.1 Extremal infinite graph theory

Following Bollob\'as~[\the\refEGT], let us take the wider notion of {\it extremal\/} graph theory to refer to how graph invariants (usually numerical) interact. A~typical question is how one invariant can be forced up or down in an arbitrary graph by making assumptions about what values another invariant takes on that graph. Particular emphasis is usually given to pairs of invariants of which one is `local' and the other `global'. The following variant of Tur\'an's original such question is a good case in point: for which function~$f(r)$, if any, can we force the presence of a $K_r$ minor in an arbitrary graph by assuming only that it has average degree at least~$f(r)$?

This type of question makes no immediate sense for an infinite graph, because there is no obvious notion of `average degree'.%
   \COMMENT{}
   So let us replace it by `minimum degree'. Now the question does make sense, but unlike in the finite case it has a trivial negative answer: there are locally finite {\it trees\/} of arbitrarily large minimum degree. Clearly, these have no $K_r$ minor for $r > 2$.

What, then, makes infinite graphs so different in this respect? Since there are finite graphs of both minimum degree and girth arbitrarily large~[\the\refBook], also finite graphs of large minimum degree can look {\em locally\/} like trees. However, a large finite tree also has many leaves, whose low degrees push its average degree back below~2. Requiring a high average or minimum degree in a finite graph therefore also has another effect, in addition to just fast growth: it forces us to `wrap up' those leaves. The classical theorem of Mader~[\the\refMaderHomEigenschaften] that a large enough average degree forces any desired minor can thus be restated as saying that, no matter how we do this `wrapping up' (whether by adding edges, say, or by identifying vertices), we will always create such a minor.

The key idea for infinite graphs now is that compactifying a locally finite graph by adding its ends wraps it up in a similar way, albeit at infinity. And this wrapping-up can help us restore properties of finite graphs that are lost in the naive transition from finite to infinite.

In fact, most of the time this restoration happens automatically: often, all we have to do is allow circles when we used to wish for a cycle, and arcs when we used to desire a path. These objects tend to exist as limits (or inside the limits) of the cycles or paths whose existence was asserted by the finite theorems that now fail. More precisely, this works whenever the cycles or paths whose existence is claimed in the finite theorem can be found there by a process of `focusing in', as discussed in Section~\secTechniques.%
   \COMMENT{}
   We shall see some more examples of this phenomenon later in this section.%
   \COMMENT{}

First, however, let us take advantage of ends as `wrapping' in a more direct way: by explicitly requiring them to have large degrees and seeing what substructures this can force.\?{define end degrees earlier}

Consider again the aim of forcing a $K_r$ minor, for $r=4$ say. In other words, we wish to force a $TK_4$ subgraph by assuming some large enough minimum degree. Assume for simplicity that our graph $G$ is 2-connected.%
   \COMMENT{}
   If $G$ is finite, a minimum degree of $\delta(G)\ge 3$ will force a $TK_4$ subgraph. Indeed, as $\kappa(G)\ge 2$ we can construct $G$ starting from a cycle and adding paths, one at a time, that share only their endvertices with the graph constructed so far. If $G\not\supe TK_4$, those two vertices lie on a common earlier path, and they do not `cross' other paths grafted later on to this same earlier path; see Figure~\figTKfour. This makes it easy to see that after every construction step there is still a vertex of degree~2.%
   \COMMENT{}

 \figure \figTKfour.
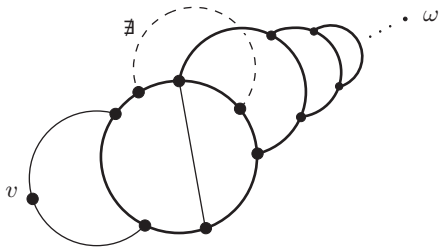
 A vertex or end of degree 2 in a graph $G\not\supe TK_4$
 (TKfour; 1000)
 
If $G$ is infinite, however, we can go on grafting new paths on to any vertices of degree~2, and after $\omega$ steps all such vertices have disappeared.%
   \COMMENT{}
   But any graph obtained in this way will have an end of vertex-degree~2. (In fact, every end will have vertex-degree~2.)%
   \COMMENT{}
   Indeed, using the infinity lemma we can find a sequence of paths each grafted on to the previous one. These paths converge to an end, which is easily seen to have vertex-degree~2. (Its edge-degree may be higher.)%
   \COMMENT{}
   We can therefore get an analogue of the finite theorem, that $\delta\ge 3$ forces a $T K_4$ subgraph, after all: every graph whose vertices {\it and ends\/} all have \hbox{(vertex-)} degree at least~3 contains a $TK_4$.

So what kind of substructures can be forced by assuming that both vertices and ends have large degree? The first theorem in this vein is due to Stein:

\proclaim Theorem~\xxxEndConnStein.~[\the\refMayaEndDeg] Let $G$ be a locally finite graph.
 \pitem{i} If $\delta(G)\ge 2k^2 + 6k$ and every end of $G$ has vertex-degree at
least ${2k^2+2k+1}$, then $G$ has a $(k+1)$-connected subgraph.
 \pitem{ii} If $\delta(G)\ge 2k$ and every end of $G$ has edge-degree~$\ge 2k$, then $G$ has a
$(k+1)$-edge-connected subgraph.
 \enditem

\noindent
   The bounds in Theorem~\xxxEndConnStein\ are close to best-possible; see~[\the\refMayaEndDeg].
   \COMMENT{}

Theorem~\xxxEndConnStein\ uses vertex and end degrees to force not a concrete desired subgraph but just some subgraph from a desired class. This is an interesting variant of the original extremal problem, which is to force a concrete subgraph (or minor etc.) by global assumptions such as on average or minimum degrees.\problem{} Forcing a concrete finite minor in an infinite graph by assuming that its vertex and end degrees are large, unfortunately, takes us little further than the $K_4$ example we saw earlier: already a $K_5$ minor cannot be forced in this way. This is because there are planar graphs with arbitrarily large degrees and only one end, of infinite degree: just take a regular tree of large degree, and add edges forming circuits~$D_i$, one for every $i\in\N$, through all the vertices at distance~$i$ from the root.

Forcing $K_r$ minors for $r\ge 5$ may be possible, however, with another notion of end degrees: one is motivated less by topological considerations (such as to give every end on a circle degree~2 there) than extremal ones. Call an induced subgraph $C$ of $G$ a {\it region\/} of $G$ if both $C$ and $G-C$%
   \COMMENT{}
   are connected and its {\it edge-boundary\/}~$B_{\rm e} (C)$, the set of edges of $G$ between $C$ and~$G-C$, is finite. Then also the {\it vertex boundary\/} $B_{\rm v} (C)$ of~$C$, the set of vertices of $C$ incident with an edge in~$B_{\rm e}(C)$, is finite. Let us say that a nested sequence $C_1\supe C_2\supe \dots$ of regions {\it defines\/} an end $\omega$ of~$G$ if the sets $\hat C_i$ form a neighbourhood basis of $\omega$ in~$|G|$, i.e., if $\omega$ belongs to every~$C_i$ and $\bigcap_i C_i =\es$. Now let the {\it relative maximum%
   \Footnote{One can also define a {\it relative average degree\/}, e.g.\ with reference to $|B_{\rm e}(C)|/|B_{\rm v}(C)|$.}
   degree\/} of $\omega$ be defined as the infimum, taken over all its defining sequences $C_1, C_2,\dots$ of regions, of the numbers
 $$\limsup_{i\to\infty} \mi_{v\in B_{\rm v}(C_i)}  d_{G[B_{\rm v}(C_i)\cup N(C_i)]}(v)\,.\eqno(*)$$
Thus, an end has relative maximum degree~$\ge k$ if and only if each of its defining sequences $C_1, C_2,\dots$ contains infinitely many $C_i$ such that every vertex in $B_{\rm v}(C_i)$ has at least $k$ neighbours in $B_{\rm v}(C_i)\cup N(C_i)$.

The first notion of such relative end degrees was introduced by Stein~[\the\refMayaBanff]. Based on her notion, she proves a theorem similar to the following:

\proclaim Theorem \xxxEndDegH.
If every end in a locally finite graph $G$ has relative maximum degree at least $k\in\N$ and $\delta(G)\ge k$, then $G$ has a finite subgraph of minimum degree at least~$k$.

\noindent
Since large enough average degree in a finite graph forces any given finite minor, Theorem~\xxxEndDegH\ implies that infinite graphs satisfying its premise have such minors too. There is also version of Theorem~\xxxEndDegH\ in terms of relative average degrees of ends.\looseness=-1

Theorem \xxxEndDegH\ is not difficult to prove. The idea is to start with a connected finite set $S$ of vertices and consider the components $C$ of~$G-S$. If every component $C$ is `good' in the sense that the local degrees of its boundary vertices as measured in~$(*)$ are at least~$k$, then $S$ together with all the vertex boundaries of components of~$G-S$ forms the desired finite subgraph. If $G-S$ has some `bad' components, we extend $S$ into these to some larger connected finite set~$S'$, and consider the components $C'$ of $G-S'$ in the same way. If this process continues infinitely long, the infinity lemma will hand us a defining sequence of regions for an end%
   \COMMENT{}
   of relative maximum degree~$<k$.

However, Theorem \xxxEndDegH\ has a serious snag: it is not clear how many graphs, if any, are such that all their ends have relative maximum degree~$\ge k$.%
   \Footnote{Indeed, if we drop the requirement in the definition of a region $C$ that $G-C$ be connected, one can show that every end in a locally finite graph has relative maximum degree~1.}
   It may be tempting, therefore, to change the definition of the relative maximum degree of an end by taking not the infimum but the supremum of the numbers in~$(*)$ over all its defining sequences of regions. Then an end would have relative maximum degree~$\ge k$ as soon as there {\em existed\/} a sequence $C_1\supe C_2\supe \dots$ of regions defining it in which for infinitely many $C_i$ all boundary vertices have large degrees in~$(*)$.

\proclaim Problem \xxxEndDegSubspaces. For which graphs $G$ does Theorem~\xxxEndDegH\ continue to hold if we replace the infimum in the definition of relative maximum end degrees with a supremum?%
   \COMMENT{}

\noindent For graphs with only countably many ends, Problem \xxxEndDegSubspaces\ has a positive answer. The general answer is negative.%
   \COMMENT{}

\bigbreak
What about forcing infinite minors? From our earlier example we know that by assuming large vertex and end degrees (non-relative) we cannot force non-planar minors. But we can force most planar minors. Indeed, by a result of Halin~[\the\refHalinInfGrid,\th\the\refBook] every graph $G$ with an end of infinite vertex-degree contains the half-grid $\N\times\Z$ as a minor,%
   \COMMENT{}
   and similarly every graph with an end of large enough finite vertex-degree contains an $n\times\N$ grid (and hence any given finite planar graph) as a minor.%
   \COMMENT{}
   On the other hand, there are planar graphs of arbitrarily large (finite) minimum degree and minimum vertex-degree for ends that do not contain the half-grid as a minor~[\the\refMayaBanff].%
   \COMMENT{}
   It would be interesting to find some natural strengthening of the degree assumption on ends that would force a planar graph to contain the full grid, or even to contain every locally finite planar graph, as a minor (cf.~[\the\refMinorUniversal]).\problem{see Comment}%
   \COMMENT{}
   More generally:

\proclaim Problem~\xxxEndDegInf.
   What infinite minors can be forced by assuming large vertex and end degrees (of any type)?

Let us return to our theme of how certain paths or cycles (with some desired properties) whose existence in a finite graph is proved by some focusing process can fail to exist in an infinite graph~$G$, because (in $G$ itself) that process need not converge. Our aim will then be to show that such paths or cycles in finite minors $G_n$ of~$G$ can tend to a limit that is an arc or circle in $|G|$ with the desired properties.

For example, consider in a finite 3-connected plane graph a maximal sequence of nested cycles (not necessarily disjoint). This sequence will end with a cycle that bounds a face.%
   \COMMENT{}
   When we delete this cycle, it will not disconnect the graph.%
   \COMMENT{}
   Cycles, in any graph, whose deletion does not reduce the connectivity of a graph by more than~3 are called {\it connectivity-preserving\/}. Cycles with the property that deleting their edges does not reduce the edge-connectivity of the graph by more than~2 are {\it edge-connectivity-preserving\/}. Such cycles exist in every finite graph, and they can be found by a process of `focusing in', just as in the planar case; this was proved by Thomassen~[\the\refCTConPres] for connectivity and by Mader~[\the\refMaderReduktion,\th\the\refMaderEdgeConPres] for edge-connectivity.%
   \Footnote{We remark that the connectivity-preserving cycle $C$ found by Thomassen~[\the\refCTConPres] is induced. Hence if $C\sub G$ with $\kappa(G) = k+3$, say, then every vertex of $C$ sends at least $k+1$ edges to $G-C$, which is $k$-connected. Thus, we also have the `mixed connectivity' result that deleting only the {\it edges\/} of $C$ reduces the ({\it vertex-\/}) connectivity of $G$ by at most~3.}%
   \COMMENT{}

In an infinite graph such cycles need not exist. Let us show this by constructing a counterexample, due to Aharoni and Thomassen~[\the\refAharoniCT]. This graph will be locally finite, and will depend on a given integer~$k$; so let us call it~$AT(k)$. The graph $AT(k)$ combines two properties that no finite graph can have at the same time:

\medskip
\item{$\bullet$} it is $k$-connected (where $k$ is as large as we like);
\item{$\bullet$} deleting any cycle, or the edges of any cycle, disconnects the graph.

\smallskip\noindent
   As a result of these two properties, the graphs $AT(k)$ are counterexamples to a number of statements which, for finite graphs, are well-known theorems. The existence of connectivity-preserving cycles is one of these, and we shall meet another below.

We shall construct $AT(k)$ inductively from copies of some fixed finite $k$-connected graph~$H$. Let us choose $H$ of girth at least~$k^2$; then $H$ has a set $X$ of $k$ vertices at distance at least~$k$ from each other. (For example, spread $X$ around a shortest cycle in~$H$.) The idea now is to build $AT(k)$ as a union of finite graphs $G_0\sub G_1\sub\dots$, where each $G_{n+1}$ is obtained from $G_n$ by grafting on to any non-separating cycles of $G_n$ some new copies of~$H$, to make them separating.

Formally, we begin with a copy $G_0$ of~$H$. Let us assume inductively that we have constructed~$G_n$ in such a way that it is $k$-connected, that any $G_{n-1}$-path%
   \Footnote{See~[\the\refBook] for the notion of an $H$-path in $G$ for a subgraph $H\sub G$.}
   in $G_n$ has length at least~$k$, and that the edges of any cycle contained in~$G_{n-1}$ separate~$G_n$. We now consider separately every cycle $C$ in~$G_n$ that does not lie in~$G_{n-1}$. By our second assumption about~$G_n$, the cycle $C$ has some $k$ edges that do not lie in~$G_{n-1}$; subdivide these once. (Thus, on an edge that lies on 3 such cycles of~$G_n$ we insert 3 subdividing vertices.) We now take $k$ fresh copies of~$H$ specific to our choice of~$C$, and identify their $k$-vertex sets $X$ with those $k$ subdividing vertices inserted on~$C$. It is easy to check that $G_{n+1}$ is again $k$-connected,%
   \COMMENT{}
   that every $G_n$-path has length at least~$k$ (because it links two vertices from $X$ in a copy of~$H$), and that the edges of any cycle contained in the newly subdivided~$G_n$ (including those inside~$G_{n-1}$) separate~$G_{n+1}$. Clearly, $AT(k) = G_0\cup G_1\cup\dots$ has the two desired properties.

\medbreak

The Aharoni-Thomassen graph $AT(k)$ for $k\ge 4$ has neither connectivity-preserving nor edge-connectivity-preserving cycles. However, it has edge-connectivity-preserving circles, as indeed does every locally finite graph:

\proclaim Theorem \xxxEdgeConnPresCircles. {\rm [\the\refTreeEnds]}
Let $G$ be a locally finite graph, and $k\in\N$. If $G$ is $(k+2)$-edge-connected,%
   \COMMENT{}
   then $|G|$ contains a circle $C$ such that the subspace of $|G|$ obtained by deleting the edges in~$E(C)$ is $k$-edge-connected.

\noindent
   (Recall that, by Lemma~\xxxConnected, the graph $G$ is $(k+2)$-edge-connected if and only if $|G|$ is, and likewise for vertex-connectivity and general connectivity.)%
   \COMMENT{}
   Theorem~\xxxEdgeConnPresCircles\ extends to standard subspaces~[\the\refTreeEnds].

\medbreak

The analogues of Theorem~\xxxEdgeConnPresCircles\ for deleting vertices and/or ends are open:

\proclaim Problem \xxxVxConnPresCircles.
\pitem{i} If $G$ is $(k+3)$-connected, does $|G|$ contain a circle $C$ such that the standard%
   \COMMENT{}
   subspace of $|G|$ obtained by deleting the vertices of $C$ and all their incident edges is $k$-vertex-connected?\problem{}
\pitem{ii} If $G$ is $(k+3)$-connected, does $|G|$ contain a circle $C$ such that the subspace of $|G|$ obtained by deleting $C$ and the edges incident with its vertices is $k$-connected?
   \enditem

\medbreak

There are also versions of these theorems and problems for deleting paths and arcs. Mader~[\the\refMaderReduktion] proved that any two vertices of a finite graph $G$ are linked by a path whose edges we can delete without reducing the edge-connectivity of $G$ by more than~2. A~construction very similar to that of $AT(k)$ provides a counterexample to this statement for locally finite graphs,%
   \COMMENT{}
   but the corresponding statement for arcs in~$|G|$ (joining two given vertices or ends) is true~[\the\refTreeEnds]. For vertex-connectivity, it is a well-known open problem of Lov\'asz whether or not there is even a function $f\colon\N\to\N$ such that any two vertices of any $f(k)$-connected finite graph can be linked by an induced path whose deletion leaves a $k$-connected graph.%
   \COMMENT{}
   The corresponding statement for arcs in~$|G|$, in the spirit of Problem~\xxxVxConnPresCircles, is also unknown.

\medbreak

Two classical theorems from finite extremal graph theory, due, respectively, to Halin and Mader, say that if a finite graph is edge-minimal with the property of being $k$-connected~[\the\refHalinMinVertex], or of being $k$-edge-connected~[\the\refMaderMinVertex], then it has a vertex of degree only~$k$. Since every finite $k$-connected or $k$-edge-connected graph contains an edge-minimal such graph, these are fundamental results about the structure of all finite $k$-connected or $k$-edge-connected graphs.

Unlike its edge-connectivity version,%
   \COMMENT{}
   the vertex-connectivity version of the above theorem remains true for infinite graphs~[\the\refHalinMinVertex] and ${k\ge 2}$:%
   \COMMENT{}
   every edge-minimal%
   \COMMENT{}
   $k$-connected graph has a vertex of degree~$k$. However, it is no longer that interesting: as the double ladder shows for $k=2$, an infinite $k$-connected graph need not contain an edge-minimal such graph. The justification for studying these minimal graphs, therefore, collapses.

However, if we extend our class of objects from graphs to all their standard subspaces, then edge-minimal objects exist by Lemma~\xxxkCon.%
   \Footnote{Having noticed that a given $k$-connected or $k$-edge-connected graphs need not contain an edge-minimal such subgraph, Halin~[\the\refHalinMinimization,\th\the\refHalinInfMinimization,\th\the\refHalinProblems] posed various problems to determine and study those that do. The shift to subspaces, coupled with Lemma~\xxxkCon, solves this problem in a much more satisfactory way.}
   For example, we can now delete all the rungs in the double ladder: what remains is disconnected as a subgraph, but 2-connected in both senses as a subspace that includes the ladder's ends. The basic objects to investigate for a study of the $k$-connected or $k$-edge-connected locally finite graphs, therefore, are their edge-minimal $k$-vertex-connected or $k$-edge-connected standard subspaces%
   \COMMENT{}%
   ---or their edge-minimal $k$-connected subspaces if they exist:

\proclaim Problem \xxxMinkCon.
Can we delete edges from any $k$-connected standard subspace of~$|G|$ to obtain an edge-minimal $k$-connected standard subspace?\problem{}

The graph of Figure~\figMinK, for example, does have an edge-minimal $k$-connected standard subspace: just delete all but one of the broken edges in Figure~\figMinK, or delete all the broken edges in Figure~\figMinKalt.

 \figure \figMinKalt.
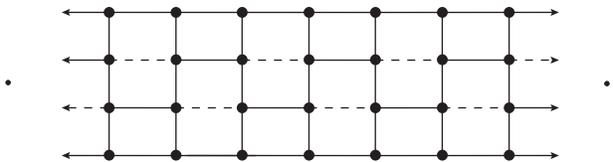
 Deleting all the broken edges yields an edge\rlap{-} minimal 3-connected standard subspace.
 (MinKalt; 1000)

\hskip-2pt
Edge-minimal $k$-edge-connected or $k$-vertex-connected standard subspaces need not have a vertex of degree only~$k$. Indeed, consider for $k=2$ the cartesian product of a 3-regular tree $T$ with an edge~$e$; deleting all the `rungs' (the edges projecting to~$e$) leaves a subspace consisting of the two copies of $T$ glued together at the ends. This subspace $X$ is edge-minimally 2-edge-connected and 2-vertex-connected, but every vertex is incident with 3 edges of~$X$. Note, however, that all the ends have edge- and vertex-degree~2 in~$X$.

This observation suggests the following infinite analogues to Halin's and Mader's finite theorems. Given~$k$, let `$k$-highly connected' mean any one of `$k$-edge-connected', `$k$-vertex-connected' or `$k$-connected', and let `degree' (for an end) mean any one of `edge-degree', `vertex-degree' or `relative degree'.

\proclaim Problem \xxxDegk.
Given $k\in\N$, does every edge-minimal $k$-highly connected standard subspace of $|G|$ contain a vertex or end of degree at most~$k$?\problem{}

Halin's finite theorem has been strengthened in various ways, and one can ask about infinite analogues also of those strengthenings. For example, Mader~[\the\refMaderMinVertices,\th\the\refMaderMinVerticesInf] proved that in an edge-minimal $k$-connected graph every \hbox{cycle} contains a vertex of degree~$k$.%
   \COMMENT{}
   Among other things this implies that every subgraph of an edge-minimal $k$-connected graph has a vertex of degree {\em at most}~$k$: either on a cycle, or else as a leaf. If Problem~\xxxDegk\ has a positive answer, it will be natural to ask (in the same informal terminology as above):

\proclaim Problem \xxxCircleDegk.
Given~$k$, does every circle in an edge-minimal $k$-highly connected standard subspace $X$ of $|G|$ contain a vertex or end whose degree in~$X$ is at most~$k$?%
   \COMMENT{}

\noindent
   For finite cycles, Problem \xxxCircleDegk\ has been answered positively by Stein~[\the\refMayaBanff].

\medbreak

Little is known about graphs that are minimally $k$-connected \wrt\ deleting vertices rather than edges (let alone vertices or ends in subspaces). Does every $k$-connected locally finite graph have a $k$-connected subgraph that is minimal in the sense of not having a $k$-connected proper subgraph?\problem{Counterex? See below.}

Perhaps this notion of minimality is too strong. Let us call a $k$-connected graph $G$ {\it minimal\/} if every $k$-connected $H\sub G$ also satisfies~$H\supe G$, that is, has a subgraph isomorphic to~$G$. An analogous definition can be adopted for subspaces.%
   \Footnote{Embeddings between subspaces should map vertices to vertices.}

Once we do have a minimal $k$-connected graph or standard subspace, what can we say about its structure? Must such a graph be finite?%
   \COMMENT{}
   Can it be $(k+1)$-connected?\problem{102 (See Comment)}%
   \COMMENT{}
   Must it have a vertex or end whose degree is small in terms of~$k$?\problem{$\to$Comment}%
   \COMMENT{}
   (By a theorem of Lick~[\the\refLick], every minimal $k$-connected finite graph has a vertex of degree at most $(3k-1)/2$.)%
   \COMMENT{}

\proclaim Problem \xxxVxMinimal. Does every $k$-connected graph or standard subspace contain a minimal such object, for some suitable notion of `minimal'? If so, what are its properties?\problem{}

See Stein~[\the\refMayaBanff] for more on this topic.

\medbreak

Another nice example of how arcs and circles in $|G|$ provide the natural setting for a classical finite theorem is tree-packing. The finite theorem here, due independently to Nash-Williams~[\the\refNWTreePacking] and Tutte~[\the\refTutteTreePacking], says that a finite graph contains $k$ edge-disjoint spanning trees unless its vertex set admits a partition, into $\ell$ sets say, such that $G$ has fewer than $k(\ell-1)$ {\it cross-edges\/}, edges between different partitions sets. (See~[\the\refBook].)

For infinite graphs, the Aharoni-Thomassen graph $AT(2k)$ provides a counterexample to this statement, even to its corollary that $2k$-edge-connected graphs have $k$ edge-disjoint spanning trees~[\the\refBook]. Indeed, since $AT(2k)$ is $2k$-edge-connected but the edge set of every cycle separates it, there can be no more than two such trees: the edges of a fundamental circuit of one tree obtained by adding an edge of another tree would separate the graph, so no third tree could be spanning.

However, the finite tree-packing theorem has a topological analogue:

\proclaim Theorem \xxxTreePacking. {\rm [\the\refBook]}
The following statements are equivalent for all $k\in\N$ and locally
finite multigraphs~$\mo G$:
   \pitem{i} $G$ has $k$ edge-disjoint topological spanning trees.
   \pitem{ii} For every finite partition of~$V(G)$, into $\ell$ sets say,
$G$~has at least ${k\,(\ell-1)}$ cross-edges.
   \enditem

\noindent
   In particular, if $G$ is $2k$-edge-connected it has $k$ edge-disjoint \TST s.

\medbreak

Theorem \xxxTreePacking\ has an interesting history: while Nash-Williams had conjectured (incorrectly) that the finite tree-packing theorem ought to extend to countable graphs verbatim, Tutte anticipated Theorem~\xxxTreePacking, even though he could not express it in the now natural topological language. See~[\the\refCyclesExpository] for details of the story.

The problem of tree-packing is closely related to the {\it arboricity\/} of a graph: the least number of forests that will cover its edges. For a finite graph~$G$, another classical theorem of Nash-Williams~[\the\refNWArboricity] says that the edges of $G$ can be covered by $k$ forests if no set of $\ell$ vertices spans more than $k(\ell-1)$ edges, the number of edges on an $\ell$-set that $k$ forests can at most provide. This theorem extends verbatim to infinite graphs,
by compactness.%
   \COMMENT{}

In our topological setting, however, it is natural to ask for more: that we can cover the edges of $G$ by $k$ {\it topological forests\/}, standard subspaces of $|G|$ that contain no circle. Interestingly, the above local sparseness condition no longer implies this; see [\the\refMayaTreePacking] for a counterexample of Bruhn. However, if we require that our graphs are also `sparse at infinity', by bounding their end degrees from above, we get the following result of Stein:

\proclaim Theorem \xxxMayaArboricity.~[\the\refMayaTreePacking]
Let $k\in\N$, and let $G$ be a locally finite graph. If no set of (say) $\ell$ vertices of $G$ induces more than $k(\ell-1)$ edges and every end of $G$ has edge-degree $<2k$, then $|G|$ contains $k$ topological forests covering all its edges.

\noindent
The bound of $2k$ in Theorem~\xxxMayaArboricity\ is sharp~[\the\refMayaTreePacking].

\medbreak

Before we leave the subject of connectivity, let us briefly summarize what is known about Menger's theorem in~$|G|$. When $a,b$ are vertices in a locally finite graph~$G$, it is well known and easy to see that if we cannot separate $a$ from $b$ by a set $S$ of fewer than $k$ vertices then there are $k$ independent $a$--$b$ paths in~$G$. This result does not gain from admitting $a$--$b$ arcs in~$|G|$ instead of just paths in~$G$, since those arcs will still meet~$S$ (Lemma~\xxxJumpingArc). Hence neither the minimum size of an $a$--$b$ separator nor the maximum number of $a$--$b$ paths (arcs) changes. The problem becomes more interesting for standard subspaces of~$|G|$, and has been solved: since standard subspaces are locally connected (Lemma~\xxxArcs), results of Whyburn~[\the\refWhyburnMenger] imply that both the point-to-point and the set-to-set version of Menger's theorem hold in them, as long as $k$ is finite. See~[\the\refCTVella] for what is known about Menger's theorem in more general 1-dimensional spaces.

For infinite~$k$, the correct version of Menger's theorem to consider is the set-to-set version in the form suggested by \Erdos\ (see~[\the\refBook]). This states that given two sets $A,B$ of vertices there is a set of disjoint $A$--$B$ paths and an $A$--$B$ separator consisting of a choice of one vertex from each of these paths. This was proved for countable graphs by Aharoni~[\the\refAharoniCtbleEM], and for arbitrary graphs by Aharoni and Berger~[\the\refAharoniBergerEM]. These results have been extended to versions where $A$ and $B$ are allowed to contain ends as well as vertices~[\the\refBruhnSteinDiestelEM], and these have to be connected by paths, rays, or double rays (not arbitrary arcs in~$|G|$). With these assumptions, Menger's theorem holds if $A\cap\overline B = \es = \overline A\cap B$, and this condition is also necessary. If we allow arbitrary arcs in~$|G|$, an example of K\"uhn (see~[\the\refCtbleEM]) shows that one even has to require that $\overline A\cap \overline B = \es$. In that case, however, we once more have the situation that arcs cannot avoid separators that meet all connecting paths, rays or double rays, so again the topological version offers nothing new~[\the\refCtbleEM].

No versions of Menger's theorem are known for standard subspaces of $|G|$ with $k$ infinite.\problem{}

\penalty-2000\medbreak

A popular area of finite graph theory which, traditionally, has no infinite counterpart is the theory of Hamilton cycles.%
   \Footnote{This is not to say that there have been no attempts. Nash-Williams and others sought to replace Hamilton cycles by spanning rays or double rays. This approach works to some extent for graphs with only one end (in which spanning double rays form circles). But it runs into difficulties as soon as the graph has more than two ends, since no ray or double ray can pass through a finite separator infinitely often. Realizing these difficulties, Halin suggested replacing Hamilton cycles with `end-faithful' spanning trees: spanning trees that contain from every end exactly one ray starting at the root (which can be chosen arbitrarily). This notion has led to some interesting problems that are still open---see e.g.~[\the\refNST]---but not to any theory related to that of finite Hamilton cycles.}
   When we replace `cycle' with `circle', however, hamiltonicity problems immediately make sense. So let us call a circle in $|G|$ a {\it Hamilton circle\/} of~$G$ if it contains every vertex of~$G$. Since circles are compact and hence closed in~$|G|$, Hamilton circles also contain every end.

What does a Hamilton circle look like? The answer to this question is somewhat daunting: as soon as the graph has uncountably many ends (which is the rule rather than the exception), any Hamilton circle must be wild~[\the\refTreeEnds], as in Figure~\figWildCircle. Still, the notion of a Hamilton circle seems to be just the right one to generalize hamiltonicity problems to infinite graphs. Let us look at some of these.

\penalty-2000

From the extremal graph theory point of view, a particularly interesting problem is how local conditions can force the (global) existence of a Hamilton cycle. Most popular among these are minimum degree conditions. When the degrees needed are large not just in absolute terms but in terms of the order of~$G$ (as in most classical results such as Dirac's theorem~[\the\refBook]), such theorems are hard to generalize to infinite~$G$.

But there are also local degree conditions that force a Hamilton cycle in a finite graph. For example, Asratian and Khachatrian~[\the\refAsratianKhachatrianLocalization] found a number of local Hamiltonicity conditions, all implying Dirac's theorem, of which the simplest version has the following infinite analogue:

\proclaim Conjecture \xxxLocalHamilton.
A~connected locally finite graph $G$ of order at least~3 has a Hamilton circle if $$d(u)+d(w) \geq |N(u)\cup N(v)\cup N(w)|$$
   for every induced path~$uvw$.

\noindent
   See also~[\the\refBook] for a proof of the finite result.

\medbreak
   Another local condition, due to Oberly and Sumner~[\the\refOberlySumner], says that a connected finite graph has a Hamilton cycle if the neighbours of each vertex span a {\it well-connected\/} subgraph: one that is connected and has independence number at most~2. (Thus, such graphs are `claw-free'.) For infinite graphs, this led Stein~[\the\refMayaBanff] to pose the following problem:

\proclaim Conjecture \xxxLocalHamiltonOS.
A~connected locally finite graph of order at least~3 has a Hamilton circle if all its vertex neighbourhoods span well-connected subgraphs.

\medbreak

A~classical sufficient local Hamiltonicity condition that does generalize to Hamilton circles is Fleisch\-ner's theorem: the {\it square\/} $G^2$ of a 2-connected finite graph~$G$ has a Hamilton cycle. (The $n$th power $G^n$ of $G$ is the graph on $V(G)$ with edges joining any pairs of vertices that have distance at most~$n$ in~$G$.) While this finite theorem is not easy, its infinite counterpart, conjectured in~[\the\refCyclesExpository] and proved by Georgakopoulos, is perhaps the deepest result about $|G|$ to date:\looseness=-1

\proclaim Theorem \xxxInfiniteFleischner.~[\the\refAgelosInfiniteFleischner] Let $G$ be a locally finite connected graph.
\pitem{i} $G^3$ has a Hamilton circle.
\pitem{ii} If $G$ is 2-connected, then $G^2$ has a Hamilton circle.

\noindent
   Thomassen~[\the\refCTFleischner] had previously proved (ii) for 1-ended graphs (in which a Hamilton circle is a spanning double ray).

\medbreak

Georgakopoulos~[\the\refAgelosOWReportFleischner] conjectured that Theorem~\xxxInfiniteFleischner\ should extend to countable graphs that are not locally finite (see Section~\secOutlook\ for subtleties about~$|G|$). This is interesting, because at first glance it seems impossible.  For since $\Omega(G)$ will be a closed subset of any Hamilton circle of~$G$, it must be compact for the conjecture to be true. But $\Omega(G)$ is compact (if and) only if no finite separator $S\sub V(G)$ splits $G$ into infinitely many components containing rays---a~property of locally finite graphs that usually fails in a countable graph. But, fortuitously, it always holds in $G^2$ (and in~$G^3$): since $G$ is connected, any component of $G^2-S$ sends a $G$-edge to~$S$, but no two components can send a $G$-edge to the same vertex of~$S$, since this would create an edge of $G^2$ between those components. Hence $G^2-S$ has at most $|S|$ components.

There are numerous other local density conditions that force a Hamilton cycle in a finite graph. For example, line graphs are `locally dense'. Thomassen~[\the\refCTLineGraphs] conjectured that every 4-connected line graph is hamiltonian, and 7-connected finite line graphs indeed are~[\the\refZhanLineGraphs]. Also, the line graph of a 4-edge-connected finite graph is hamiltonian~[\the\refCatlinLineGraphs].

\proclaim Problem \xxxLineGraphs.
Does sufficient connectivity force the line graph of a locally finite connected graph to have a Hamilton circle?\problem{}

\noindent
   Georgakopoulos~[\the\refAgelosOWReportFleischner] conjectured that this should be true with the same connectivity assumptions as are currently known for finite graphs.

\medbreak

Another classical result about finite Hamilton cycles is Tutte's hamiltonicity theorem for 4-connected planar graphs. It was the following conjecture of Bruhn (see~[\the\refCyclesExpository]) that first advanced the notion of a Hamilton circle:

\proclaim Conjecture \xxxInfiniteTutte.
Every 4-connected locally finite planar graph has a Hamilton circle.

\noindent
  Like the extension of Fleischner's theorem, this appears to be a hard problem. Partial results have been obtained by Bruhn~\&~Yu~[\the\refBruhnYuHamilton] and by Cui, Wang \& Yu~[\the\refYuHamilton].

\medbreak

Another famous hamiltonicity problem for finite graphs is the toughness conjecture (see~[\the\refBook]), and its analogue for Hamilton circles is equally intriguing. The problem might become easier if we ask only that the circle must pass through all the ends, but not necessarily through all the vertices:%
   \COMMENT{}

\proclaim Problem \xxxToughness.
  \pitem{i} Is there an integer $t$ such that every $t$-tough locally finite graph contains a Hamilton circle?
  \pitem{ii} Is there an integer $t$ such that if deleting $tk$ vertices from a locally finite graph $G$ never leaves more than $k$ infinite components then $|G|$ contains a circle through all its ends?
   \enditem

\noindent
   Note that 1-tough\-ness, let alone the assumption that deleting $k$ vertices never leaves more than $k$ infinite components, is not enough to ensure that $|G|$ contains a circle through all its ends.%
   \Footnote{Pick a vertex in a complete graph~$K_4$ and turn each of its three incident edges into a ladder, the original edge becoming its first rung.%
   \COMMENT{}
   This graph contains no circle through its three ends, but deleting at most $k$ vertices never leaves more than $k$ components.}
   Asking for `hamiltonicity for ends' as in~(ii) may be interesting also as a weakening of other hamiltonicity conjectures.\problem{}

\medbreak

There are also very interesting results and conjectures about Hamilton cycles in finite graphs asserting the sufficiency of conditions for their existence that are themselves global. Think of hamiltonicity problems for Cayley graphs,%
   \COMMENT{}
   for sparse expanders,%
   \COMMENT{}
   or for products of graphs.%
   \COMMENT{}
   Some of these may be extendable to Hamilton circles, perhaps under additional assumptions.\problem{}

Let us close this section by remarking that most of the theorems and problems we discussed have meaningful analogues in arbitrary standard subspaces of~$|G|$ rather than just $|G|$ itself. Sometimes, these extensions are easy and can be obtained by imitating the proof for~$
|G|$.%
   \COMMENT{}
   But at other times they can be challenging.%
   \COMMENT{}
   The reader is invited to explore this further.\problem{}

\subsection\secApplications.2 Cycle space applications

The earliest, and so far the most successful, applications of our topological approach to locally finite graphs have been results which, for finite graphs, relate the cycle space of $G$ to it structural properties. We have seen one such example: Theorem \xxxCBasics, which says that a standard subspace of $|G|$ can be covered by a topological Euler tour if and only if its edge set lies in~$\cyclespace(G)$. There had been a number of earlier attempts to generalize Euler's finite theorem to infinite graphs, based on double rays as infinite analogues of finite cycles, but these attempts were hampered from the outset by the handicap that a double ray cannot visit more than $|S|+1$ components of $G-S$ for any finite set $S$ of vertices, and hence could not really succeed.

Another early application is MacLane's planarity criterion. Call a set $\D\subseteq \C(G)$ {\it sparse\/} if no edge of $G$ lies in more than two elements of~$\D$. MacLane's theorem says that a finite graph is planar if and only if its cycle space has a sparse generating subset. If the graph is 3-connected, these generators will necessarily be its peripheral circuits (all but at most one), the face boundaries in any drawing. (See~[\the\refBook], and [\the\refMacLaneArbitrarySurfaces] for a generalization to arbitrary surfaces.) For our infinite~$G$, this fails unless we allow infinite circuits: the 3-connected graph $G$ in Figure~\figBadEdge, for example, has no finite face boundary containing the edge~$e$, so its finite peripheral circuits do not even generate $\C\fin(G)$.%
   \COMMENT{}

 \figure \figBadEdge.
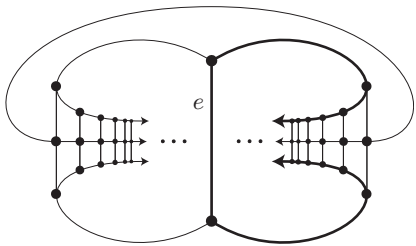
 The edge $e$ lies on no finite face boundary, so these do not generate the cycle space
 (BadEdge; 1000)

Solving a long-standing problem of
Wagner~[\the\refWagnerBook], Bruhn and Stein extended MacLane's theorem to infinite graphs, using~$|G|$:

\proclaim Theorem \xxxMacLane.~[\the\refBruhnSteinMacLane]
   $G$ is planar if and only if $\C(G)$ has a sparse generating subset.

Just as MacLane turned the algebraic properties of the face boundaries of a plane graph into a planarity criterion, Archdeacon, Bonnington and Little~[\the\refABLLeftRightTours] found an algebraic planarity criterion in terms of the `left-right-tours' of finite plane graphs. Bruhn et al~[\the\refBruhnBicycles] extended this to locally finite graphs, based on possibly infinite `left-right-tours' in~$|G|$.

\medbreak

The planarity criterion of Kelmans and Tutte says that a finite 3-connected graph is planar if and only if its set of peripheral circuits is sparse. This follows from MacLane's theorem and the theorem of Tutte that, in any 3-connected finite graph (planar or not), the peripheral circuits generate the cycle space.

The Kelmans-Tutte theorem, too, fails for infinite graphs unless we allow circuits to be infinite: sparseness of the finite peripheral circuits alone does not imply planarity~[\the\refBruhnSteinMacLane]. But Bruhn~[\the\refBruhnPeripheral] extended Tutte's generating theorem to the topological cycle space (Theorem \xxxGenerating\th(ii)). Hence for~$\C(G)$ rather than~$\C\fin(G)$, the infinite MacLane theorem implies the Kelmans-Tutte criterion:

\proclaim Theorem \xxxKelmans.
   If $G$ is 3-connected, then $G$ is planar if and only if every edge lies in
at most two peripheral circuits.

Another classical result in this context is Whitney's duality theorem for finite graphs. It is often thought of
as a planarity criterion, but can equally be viewed as a topological
characterization of the graphs that have an (algebraic) dual: that these are precisely the planar graphs. As explained in Section~\secTheory.5, duality for infinite graphs, and in particular any analogue of Whitney's theorem, will work only in the class of finitely separable graphs~$G$ , and for the topological cycle and cut spaces of the compactifications $\tilde G$ obtained by adding only the undominated ends as new points (and making any other rays converge to their dominating vertex). In that setting, however, Whitney's theorem does generalize smoothly:

\proclaim Theorem~\xxxWhitney.~[\the\refDuality]
   A finitely separable graph has a dual if and only if it is planar.

By colouring-flow duality (see~[\the\refBook]), the four-colour theorem can be
rephrased as saying that the edge set of any finite planar bridgeless graph $G$ is the union of two algebraic cycles: this is equivalent to the existence of a ($\Z_2\times \Z_2$)-flow on~$G$, which in turn is equivalent the the 4-colourability of any dual of~$G$.%
   \COMMENT{}
   Since infinite planar graphs are 4-colourable by the \Erdos~--~de~Bruijn theorem~[\the\refBook], Theorem~\xxxWhitney\ has the following double cover corollary:

\proclaim Corollary~\xxxFCT.~[\the\refDuality]
   Assume that $G$ is finitely separable, and let $\cyclespace(G)$ be based on the space~$\tilde G$. If $G$ is planar and bridgeless, then $E(G)$ is the union of two algebraic cycles.

\subsection\secApplications.3 Flows in infinite graphs and networks

Most of this section deals with electrical flows in infinite networks: an exciting field with rich connections to other branches of mathematics, and one where the study of $|G|$ has a natural well-recognized place. Towards the end of the section we also discuss algebraic (group-valued) and (non-electrical) network flows. Basing these on $|G|$ rather than just $G$ is likely to lead to an extension of most of the known finite theory to locally finite graphs, just as in Sections~\secApplications.1--2. 

An {\it electrical network\/} is a locally finite connected graph $(V,E)$ whose (undirected) edges~$e$ have real {\it resistances\/} $r(e)>0$ assigned to them, and which has two specified vertices $s$ and~$t$, called the {\it source\/} and the {\it sink\/}. A {\it flow\/} in this network is a real function $f\:\vE\to\R$ on the set of orientations $\ve$ and $\ev$ of these edges $e$ such that $f(\ve) = -f(\ev)$ for every~$e$ and
 $$\sum_v f(\ve) = 0 \hbox{ \sl for every vertex }v\notin\{s,t\}, \eqno{\rm (KH1)}$$
where the sum $\sum_v f(\ve)$ ranges over all edges $\ve$ at $v$ oriented away from~$v$. The {\it value\/} of this flow $f$ is the number $\sum_s f(\ve)$.

Let us assume for the moment that $G$ is finite. We can then use (KH1) to show that the value of~$f$ is equal to $\sum_{\ve\in\vF} f(\ve)$ for every cut $\vF$ separating $s$ from~$t$, oriented from the side containing $s$ to that containing~$t$. And for every $i\in\R$ there is now a unique flow $f$ on $G$ of value~$i$ that also satisfies
 $$\sum_{\ve\in\vC} f(\ve) r(e) = 0 \hbox{ \sl around every oriented circuit $\vC$ in}~G. \eqno{\rm (KH2)}$$
Indeed, the existence of some flow of value~$i$ is clear: we just send a flow of value~$i$ along an $s$--$t$ path and let $f(\ve)=0$ elsewhere. To find a flow of value~$i$ that also satisfies~(KH2), one can use the {\it total energy\/} of~$f$: the number
 $$w(f) := \sum_{e\in E(G)} f^2(e) r(e)\,,$$
where $f^2(e) := f^2(\ve) = f^2(\ev)$. It is easy to show that among all the flows of value~$i$ there is one of minimum total energy, and a short calculation shows that this flow satisfies~(KH2). Uniqueness now comes for free, since there cannot be two flows $f,f'$ of the same value that both satisfy~(KH2). Indeed, their difference $g = f-f'$ would be a flow of value~0 such that $g(\ve_0)>0$ for some oriented edge~$\ve_0$. Using that $\sum_v g(\ve) = 0$ at every vertex~$v$ (including $s$ and~$t$ now, since $g$ has value zero), we `move greedily along' oriented edges of~$G$, starting at~$\ve_0$, to find an oriented circuit $\vC$ with $g(\ve) > 0$ for all $\ve\in\vC$. Hence $g = f-f'$ does not satisfy~(KH2), but it should, since $f$ and~$f'$ do.

We remark that condition (KH2) is equivalent to the existence of a {\it potential function\/} inducing~$f$, a function $p\: V\to\R$ related to $f$ via {\it Ohm's law\/} that $f(\ve)r(e) = p(v)-p(u)$ for every oriented edge $\ve = uv$ of~$G$: this is easily seen by considering a spanning tree.

For the case that $r(e)=1$ for all~$e$, Ohm's law says in topological terms (cf.\ Section~\secOutlook.3) that $f = \delta^0 p$, where $\delta^0$ is the coboundary operator on 0-chains. Then being `induced by a potential function via Ohm's law' is the same as lying in the image of~$\delta^0$. In our terminology, the image of~$\delta^0$ is the oriented cut space: every oriented cut $f = \lambda\vE(X,Y)$ has the form~$\delta^0 p$ (let $p$ assign $\lambda$ to the vertices in~$Y$ and 0 to those in~$X$), and conversely we can write any function of the form $f = \delta^0 p$ as a sum of oriented atomic bonds~$\lambda\vE(v)$.%
   \COMMENT{}
   The fact that it is precisely these functions~$f\:\vE\to\R$%
   \COMMENT{}
   that satisfy (KH2) then is the oriented version of Theorem~\xxxOrthogonality\th(ii). Finally, such a function $f = \delta^0 p$ satisfies~(KH1) if and only if $p$ is harmonic at every vertex $v\notin \{s,t\}$: as
 $$\sum_v f(\ve) = \sum_v (\delta^0 p) (\ve) = \bigg(\sum_{u\in N(v)} p(u)\bigg)  -\ d(v)\, p(v)\,,$$
   the net flow out of a vertex $v$ is zero if and only if $p(v)$ equals the average of the $p$-values at its neighbours.%
   \COMMENT{}

So what about infinite networks? Given any $i\in\R$, there will still be a flow of value~$i$ (just along an $s$--$t$ path), there will be a flow of value~$i$ with minimum total energy (defined as before)~[\the\refAgelosUniqueFlows], and this flow will satisfy~(KH2) for finite circuits.%
   \COMMENT{}
   However, there may now be more than one such flow of value~$i$: Figure~\figElusive\ shows a particularly striking example, which satisfies (KH2) because the graph contains no circuit at all. (Note that its second flow `from $s$ to~$t$' would even exist if we deleted the edge $st$, disconnecting $s$ from~$t$.)

 \figure \figElusive.
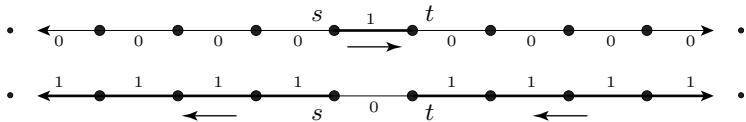
 Two flows of value 1 satisfying~{(KH2)}%
   \COMMENT{}
 (Elusive; 1000)
 
Such `flows' are not always desirable, and one would like to be able to amend the definition so as to exclude them. But it has not been clear until recently how this can be done without also excluding desirable flows. From our topological perspective, however, this is much clearer: the requirement of $\sum_v f(\ve) = 0$ should be applied not only to vertices $v$, but also, somehow, to ends. (In the second flow of Figure~\figElusive, a net flow of~1 disappears into the left end, while a net flow of~1 emerges from the right end, violating the intended requirement at both ends.)

The simplest way to make such an additional requirement for ends is without mentioning ends directly: instead of (KH1) we require a condition that is equivalent to (KH1) for finite networks, but stronger for infinite ones: that $f$ should sum to zero across {\em any finite\/} oriented cut not separating $s$ from~$t$. Formally:
 \textno $\sum_{\ve\in\vF} f(\ve) = 0$ for every finite oriented cut $\vF$ of $G$ that does not separate $s$ from~$t$. &\rm (KH1')

\smallskip\noindent
Let us call flows satisfying (KH1$'$) {\it non-elusive\/}. It is not hard to show that the value $\sum_s f(\ve)$ of a non-elusive flow $f$ equals the total value $\sum_{\ve\in\vF} f(\ve)$ of $f$ across any finite oriented cut~$\vF$ that does separate $s$ from~$t$, oriented from its side containing $s$ to that containing~$t$.%
   \COMMENT{}

For reference, let us also restate (KH2) with explicit reference to finite circuits only:
 $$\sum_{\ve\in\vC} f(\ve) r(e) = 0 \hbox{ \sl around every finite oriented circuit $\vC$ in}~G. \eqno{\rm (KH2')}$$
 
For networks whose {\it total resistance\/} $\sum_{e\in E(G)} r(e)$ is finite, Georgakopoulos~[\the\refAgelosUniqueFlows] proved that non-elusive flows that satisfy (KH2$'$) and have finite total energy are indeed unique for any given value. This is interesting, since (KH2$'$) makes demands only on finite circuits, while the difference between two flows of the same value might have an infinite circuit as its carrier.%
   \Footnote{By the oriented version of Theorem~\xxxOrthogonality\th(i), it will lie in the oriented cycle space~$\vCC(G)$.}
   In particular, we cannot prove uniqueness directly from~(KH2$'$), as for finite graphs, but have to use a limit construction as discussed in Section~\secTechniques.

\proclaim Theorem \xxxFlows.~[\the\refAgelosUniqueFlows]
Let $G$ be a locally finite network of finite total resistance. For every $i\in\R$ there is a unique non-elusive flow%
   \COMMENT{}
   in $G$ of value~$i$ and finite total energy that satisfies~{\rm (KH2$'$).} This flow also satisfies~(KH2) for infinite circuits~$\vC$.%
   \Footnote{In particular, the infinite sum $\sum_{\ve\in\vC} f(\ve)\, r(e)$ is well defined and finite.}

Theorem \xxxFlows\ is best possible in that uniqueness fails if we do not require finite total resistance.%
   \COMMENT{}
   And the assumption of finite total resistance is not, in fact, unnatural. Physically, it means that the entire network can be `cut out of a single finite piece of wire'. Mathematically, it means that the network admits $|G|$ as a topological model. More precisely, if we interpret the numbers $r(e)$ as edge lengths (as in Figure~\figFlow), then $r$ makes $G$ into a metric space. Let $\overline G$ be its completion. Georgakopoulos~[\the\refAgelosLTop] showed that if $\sum r(e) < \infty$ then $\overline G$~is homeomorphic to~$|G|$, by a homeomorphism that is the identity on~$G$ (Theorem~\xxxLtopModG). We shall look at metric completions of $G$ in more detail in Section~\secOutlook.

Consider any locally finite connected graph~$G$, a function $r\: E(G)\to \R^+$ such that $\sum_{e\in E(G)} r(e) < \infty$, and a function $p\: V(G)\to\R$ such that $f:= \delta^0 p / r$%
   \COMMENT{}
   has finite total energy and satisfies
 \textno $\sum_{\ve\in\vF} f(\ve) = 0$ for every finite oriented cut $\vF$ of~$G$. &\rm (KH1'')

\noindent
   Then for any choice of $s,t\in V(G)$, our function $f$ is a non-elusive flow in the network $(G,r,s,t)$ of value~$0$ that satisfies~(KH2$'$).%
   \COMMENT{}
   Since the constant function $\vE\to\{0\}$ is also such a flow, Theorem~\xxxFlows\ implies that our $f$ is in fact constant with value~0, and hence $p$ too must have been constant.

However if we relax (KH1$''$) to~(KH1), there may be more such functions $f$ and~$p$. And since we did not specify a source or sink among the vertices, these functions are more natural than the example of Figure~\figElusive: $f$~may be viewed as a flow with (possibly many) sources and sinks at infinity, as in Figure~\figFlow.

 \figure \figFlow.
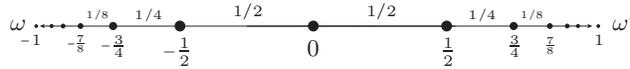
 {A flow of value~1 from $\omega$ to~$\omega'$, with resistances (top) and potentials (bottom)}
 (Flow; 1000)

To rule out other pathological examples, however, we should require now that $p$ extends continuously to the ends:\problem{proved by Agelos?}

\proclaim Problem \xxxBoundarySourceSink.
Let $G$, $r$ and $p$ be given as above, with $f = \delta^0 p/r$ satisfying (KH1) but not necessarily~(KH1$''$). Does $p$ extend continuously to~$\Omega(G)$?%
   \COMMENT{}

The converse of Problem~\xxxBoundarySourceSink\ is the (discrete) {\it Dirichlet problem at infinity\/}:\problem{ask Agelos or Woess who proved what.}

\proclaim Problem \xxxDirichlet.
Which continuous functions $p\:\Omega(G)\to\R$ extend continuously to~$V(G)$ so that $\delta^0 p /r$ satisfies~(KH1)?

\noindent
   The Dirichlet problem has been widely studied; see, e.g.,~[\the\refBenjaminiSchramm,\th \the\refWoessDirichlet,\th \the\refWoessBook].

\penalty-1000\medbreak

Let us close this section with a glance at how non-electrical network and algebraic flow theory extend to locally finite graphs. As before, we wish to exclude `flows' that issue from or dissipate to infinity, and will adjust our definitions accordingly. Let $G$ be a locally finite connected graph. Call a function $f$ on~$\vE(G)$ {\it symmetrical\/} if $f(\ve) = -f(\ev)$ for all oriented edges~$\ve\in\vE(G)$.

Specifying two vertices $s,t$ as source and sink turns $G$ into a (non-electrical) network. A~{\it flow\/} through this network is a symmetrical function on~$\vE(G)$; it is {\it non-elusive\/} if it satisfies~(KH1$'$).

\proclaim Problem \xxxNflows.
Extend the non-algorithmic aspects of network flow theory to non-elusive flows in locally finite graphs.

\noindent
   The max-flow min-cut theorem has been extended in this way; see the last few pages of~[\the\refMFMC].

\medbreak

Turning now to algebraic flows, let $H$ be any abelian group. Call a symmetrical function $f\:\vE(G)\to H$ an $H$-{\it flow\/} on~$G$ if $f(\ve)\ne 0$ for all $\ve\in\vE$, and call it {\it non-elusive\/} if it satisfies~(KH1$''$). Such a function induces $H$-flows also on the finite minors $G_n$ of~$G$ defined in Section~\secTechniques.1. It seems that, using compactness as explained in Section~\secTechniques.2. (and \TST s wherever spanning trees are needed), one can extend most---but not all%
   \Footnote{Infinite bipartite cubic graph need not have a 3-flow, even when their ends have vertex- and edge-degree~3~[\the\refTheoFlows].}%
   \COMMENT{}%
   ---the standard results about $H$-flows in finite graphs to locally finite graphs:

\proclaim Problem \xxxHflows. Extend the algebraic flow theory of finite graphs to non-elusive algebraic flows in locally finite graphs.\problem{Easy; see fil\rlap{e} \rlap{InfiniteFlows.pdf}}

\section\secOutlook. Outlook

The topological approach to studying locally finite graphs discussed in this paper can be taken further, in various directions. In this section we address some of these:
  \item{$\bullet$} graphs that are not locally finite;
  \item{$\bullet$} locally finite graphs with finer `boundaries' than ends, such as boundaries of (Cayley graphs of) hyperbolic groups in the sense of Gromov;
  \item{$\bullet$} homology of non-compact spaces other than graphs;
  \item{$\bullet$} infinite matroids.
\medbreak

\subsection \secOutlook.1 Graphs with infinite degrees

Consider the two graphs shown in Figure~\figFans, the {\it fan\/} and the {\it double fan\/}. In each of theses graphs there is a ray $R$ that is dominated by a vertex (respectively, by two vertices).%
   \Footnote{See the end of Section~\secTheory\ for the formal definition of `dominate'.}

\vskip0pt\penalty-200
\figure \figFans.
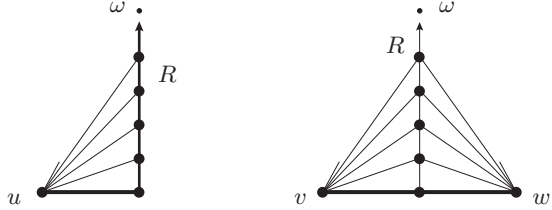
 Should the fat edge sets be circuits?
 (Fans; 1000)

In both these graphs, the sum (mod-2) of all the triangle circuits equals the set of heavy edges. These edges do not look like circuits: they are the edge sets of a ray, or of a path of length~2. However, as thin sums of finite circuits, these edge sets ought to be elements of the cycles space $\cyclespace(G)$: remember that we need to allow infinite thin sums in $\cyclespace$ in order to make important applications of $\C$ work, such as MacLane's theorem. But if these sets of edges lie in~$\cyclespace$, they ought to be disjoint unions of circuits%
   \vadjust{\penalty-200}
   (if we want Theorem~\xxxGenerating\th(i) to extend to these graphs). Hence we feel compelled to admit those edge sets as circuits---or to abandon much of the cycle space theory we just established for locally finite graphs.

So can we make sense of the notion that those edge sets should be circuits? Indeed it seems we can; and the answer has to do with what topology we choose for $|G|$ when $G$ is not locally finite.

Let us be more formal. Let $G$ be any connected graph.%
   \COMMENT{}
   The ends of $G$ are defined as earlier, as equivalence classes of rays. As the topology on $G$ itself we take slightly fewer basic open sets than the identification topology for a 1-complex does:%
   \Footnote{That topology leads to a topology for $|G|$ called \fnTop. Historically, this was the first topology for $|G|$ to be considered, but it has few advantages over the topologies \fnMTop\ and \fnVTop\ discussed here. In particular, $|G|$ under \fnTop\ is neither compact nor metrizable as soon as a vertex has infinite degree. See~[\the\refSpanningTrees,\th \the\refCyclesOne] for more.}
   around every vertex $v$ we take as basic open sets only the open stars $E_\epsilon(v)$ of length~$\epsilon$ for arbitrary but fixed $\epsilon > 0$, the same~$\epsilon$ for every edge at~$v$. (Then the sets $E_{1/n}(v)$ form a countable neighbourhood basis of~$v$. In the 1-complex topology, vertices of infinite degree have no countable neighbourhood basis.) As basic open sets around ends $\omega$ we take for every finite set $S$ of vertices and every $\epsilon$ with $0 < \epsilon \le 1$ the open $\epsilon$-collar $\hat C_\epsilon(S,\omega)$ of~$\overline{C(S,\omega)}$.%
   \Footnote{More formally, let $\hat C_\epsilon(S,\omega)$ be obtained from the set $\hat C(S,\omega)$ defined in Section~\secConcepts\ by replacing each open $C$--$S$ edge $\interior e = (c,s)$ with its initial open segment of length~$\epsilon$, assuming that $e$ itself has length~1.}
   Note that, as for vertices, these $\epsilon$-collars have uniform width.

Let us call this topology \MTop.%
   \Footnote{The `M' comes from fact that \fnMTop\ makes $|G|$ metrizable, at least for countable connected~$G$ (see below).}
   Under this topology, $|G|$ is clearly Hausdorff, in fact normal~[\the\refPhilippNormal]. However, unless $G$ is locally finite, $|G|$~will no longer be compact. (Indeed, consider an open cover that consists the middle half of $e$ for every edge $e$ at some fixed vertex~$v$, and one further open set to cover the rest of~$|G|$ but none of the mid-points of edges at~$v$. If $v$ has infinite degree, this cover has no finite subcover.)

In order to give $|G|$ a chance to be compact, we have to take a coarser topology still: we take the same open stars $E_\epsilon(v)$ around vertices, but around ends $\omega$ we take only the open sets $\hat C_\epsilon (S,\omega)$ with $\epsilon=1$, our original sets~$\hat C(S,\omega)$. This topology is called~\VTop. 

\proclaim Theorem \xxxVTop.~[\the\refSpanningTrees] 
   For every graph~$G=(V,E)$, not necessarily connected or locally finite, the following statements are equivalent under \VTop{\rm:}
   \pitem{i} $|G|$ is compact.
   \pitem{ii} For every finite $S\sub V$, the graph $G-S$ has only finitely many
componen\rlap{ts.}
   \pitem{iii} Every closed set of vertices is finite.
   \enditem

For $G$ countable and connected, Theorem \xxxVTop\ implies that $|G|$ is compact under \VTop\ if and only if $G$ has a locally finite spanning tree, which can be chosen normal~[\the\refSpanningTrees]. Note that since \VTop\ agrees with \Top\ and~\MTop\ on its closed subspace~$\Omega(G)$, condition~(ii) of Theorem~\xxxVTop\ characterizes, for any of these topologies, the graphs whose end space $\Omega$ is compact.

We complete our excursion to the various topologies of $|G|$ by mentioning one more, known as \ETop. Strictly speaking, this is not a topology on~$|G|$,%
   \Footnote{\dots unless $G$ is locally finite, in which case its edge-ends coincide with its usual ends.}
   but on a similar space obtained by adding to $G$ its {\it edge-ends\/} rather than its ends: the equivalence classes of rays under which two rays are equivalent if no finite set of {\it edges\/} separates them. The basic open sets of this space are the components of $G-X$ left by deleting a finite set $X$ of inner points of edges (and an edge-end belongs to the component in which all its rays have a tail).%
   \Footnote{Such a `component' is a graph-theoretical component left by deleting from $G$ the edges containing points from~$X$, together with the segments of these edges (after deleting~$X$) that have an endpoint in such a component.}
   If $G$ is connected, as we assume, then this space is compact~[\the\refSchulzEdgeEnds].

On the set of vertices and edge-ends, being separated by such a finite set~$X$ (i.e., by a finite set of edges) is an equivalence relation. The open sets in \ETop\ cannot distinguish such equivalent points. But identifying equivalent points (including different vertices, such as the vertices $v$ and~$w$ in Figure~\figFans) yields a compact Hausdorff space, which is metrizable if $G$ was countable.%
   \COMMENT{}
   See [\the\refHahnEdgeEnds, \the\refSchulzEdgeEnds, \the\refAgelosLTop] for more on edge-ends and \ETop.

\subsection \secOutlook.2 The identification topology

Let us turn back to our two graphs from Figure~\figFans. Under \MTop, the sets of fat edges are still far from being circuits. Under \VTop\ and~\ETop, however, they almost are: these topologies cannot distinguish an end from vertices dominating it,%
   \Footnote{More generally, it is easy to see that a space $|G|$~is Hausdorff under \fnVTop\ if and only if no end of $G$ is dominated.}%
   \COMMENT{}
   and if we identified the end with their dominating vertices, the fat edge sets would indeed become circuits. In the single fan on the left, this does indeed look right: why introduce, when forming $|G|$ from~$G$, a new limit point for the ray $R$ if a natural limit point, the vertex~$u$, already exists? In the double fan on the right, however, one may have more qualms: identifying both $v$ and $w$ with the end~$\omega$ would result in the identification of two vertices, changing the graph.

So shall we identify ends with their dominating vertices, or shall we not? Before we describe a radical and intriguing answer to this question due to Georgakopoulos~[\the\refAgelosLTop], let us deal with the case when the dilemma does not arise: the case that no end is dominated by more than one vertex. This is the case, for example, for the finitely separable graphs defined in Section~\secTheory.5, a~natural class that occurs quite independently of this problem. For such graphs~$G$, let $\tilde G$ denote the identification space obtained from $|G|$ by identifying each vertex with all the ends it dominates. This is a Hausdorff space---even under \VTop\ when $|G|$ fails to be Hausdorff~[\the\refTST]. Other features of the topology of~$\tilde G$, such as compactness, depend on the topology used for~$|G|$. If \VTop\ was used on~$|G|$, the identification topology on $\tilde G$ is known as \ITop.

Finitely separable 2-connected graphs $G$ are countable,%
   \COMMENT{}
   and no finite set of vertices separates~$G$ into infinitely many components.%
   \COMMENT{}
   Then $|G|$ is compact under \VTop\ (Theorem~\xxxVTop), and hence so is~$\tilde G$ under~\ITop.%
   \COMMENT{}
   Answering a question of Spr\"ussel~[\the\refPhilippNormal], Richter and Vella~[\the\refRichterVella] extended%
   \COMMENT{}
   this by showing that $\tilde G$ is in fact `Peano':

\proclaim Theorem \xxxRichterVellaPeano.~[\the\refRichterVella] Let $G$ be a 2-connected finitely separable graph.%
   \Footnote{This assumption can be weakened: it is only necessary to assume that no two vertices can be linked by infinitely many independent paths~[\the\refRichterVella]. Under this weaker assumption, however, $\tilde G$~can contain arcs and circles consisting entirely of vertices or ends. Such spaces were studied in~[\the\refRichterVella].%
   \COMMENT{}}
   Then $\tilde G$ under \ITop\ is a compact, connected, locally connected, metric space.

For more on~$\tilde G$, see~[\the\refTST,\th\the\refRichterVella,\th\the\refCTVella].

\subsection \secOutlook.3 Compactification versus metric completion

As mentioned already in Theorem~\xxxCompactMetricG, the compactification $|G|$ of a connected and locally finite graph $G$ is metrizable. When $G$ is not locally finite, $|G|$~with our usual topology, \VTop, fails to be even Hausdorff (as soon as there is a dominated end), and no vertex of infinite degree has a countable neighbourhood basis in the usual 1-complex topology.%
   \COMMENT{}
   However, if we consider $|G|$ with \MTop, we have the following result:

\proclaim Theorem~\xxxMetricNST.~[\the\refSpanningTrees]
For a connected graph~$G$, and $|G|$ endowed with~\MTop, the space $|G|$ is metrizable if and only if $G$ has a \NST.

\noindent
   Unless otherwise mentioned, the space $|G|$ will in this section always carry the topology \MTop.

\medbreak

Theorem~\xxxMetricNST\ is interesting in both directions: as a characterization of the graphs $G$ for which $|G|$ is metrizable, but also as an unexpected topological characterization of the graphs admitting a \NST. These are known to include all countable graphs, and they have been characterized by two types of forbidden substructure~[\the\refNST]. But these substructures themselves are not fully understood, and classifying them more accurately remains a challenging open problem.

The metric which a \NST\ $T$ of $G$ induces on~$|G|$ is easy to describe. We first assign length $\ell(e) := 2^{-n}$ to every tree-edge $e$ from level $n-1$ to level~$n$, the root being at level~0. This defines a metric on~$V(G)\cup\Omega(G)$, via
 $$d_\ell(x,y) := \sum_{e\in xTy} \ell(e)\,,$$
 where $xTy$ is the path, ray or double ray in $\overline T$ from $x$ to~$y$. We then extend this metric to one on $|G|$ by mapping every edge $xy$ homeomorphically to a real interval of length~$d_\ell (x,y)$. As one can check~[\the\refSpanningTrees], this metric induces the topology \MTop\ on~$|G|$.

Instead of defining distances between ends and other points explicitly, we could also have defined them only for points of $G$ itself (as above), and then taken the completion of this metric space~$(G,d_\ell)$. Since two Cauchy sequences of points in $(G,d_\ell)$ are equivalent%
   \COMMENT{}
   if and only if in~$|G|$ they converge to the same end, this would have yielded the same result: the complete metric space on $G\cup\Omega(G)$ we defined explicitly above, and whose metric induces \MTop\ on~$|G|$.

Georgakopoulos~[\the\refAgelosLTop] showed that this is not just a feature of our particular metric. Given any function $\ell\: E(G)\to \R^+$ assigning positive `edge lengths', defining distances between vertices $u,v$ via
$$d_\ell (u,v)\ =\ \inf\! \Big\{\sum_{e\in P}\ell(e)\ \big|\
        \hbox{$P$ an $u$--$v$ path in $G$}\Big\}$$%
   \COMMENT{}
   and identifying any $u,v$ with $d_\ell(u,v) = 0$ defines a metric on the resulting quotient space of~$V(G)$, which extends to a metric $d_\ell$ on the corresponding quotient of the entire graph $G$ once we fix homeomorphisms between the edges and real intervals of their respective lengths. Let us denote this metric space as~$(G,d_\ell)$, and the topology it induces on its completion as $\ell$-\Top. The metric subspace induced by the completion points we added is the $\ell$-\Top\ {\it boundary\/} of~$G$.

\proclaim Theorem~\xxxLtopModG.~[\the\refAgelosLTop]
If $G$ is any countable connected graph and $\ell\: E(G)\to \R^+$ satisfies $\sum_{e\in G} \ell(e) < \infty$, then completing the metric space $(G,d_\ell)$ yields the edge-ends of $G$ as completion points and \ETop\ as the induced topology of the completion. If $G$ is locally finite, this coincides with~$|G|$.

In view of our earlier \NST\ example one should expect that, for $G$ locally finite, the finiteness condition in Theorem~\xxxLtopModG\ can be relaxed considerably without losing that the metric completion of $G$ coincides with~$|G|$:

\proclaim Problem \xxxFiniteSums.
Given a countable connected graph~$G$, characterize the functions $\ell\: E(G)\to \R^+$ for which the completion of $(G,d_\ell)$ coincides with~$|G|$.\problem{}%
   \COMMENT{}

For functions $\ell$ not satisfying $\sum_{e\in G} \ell(e) < \infty$, the metric-completion approach opens up a wide range of new possibilities over the purely topological compactification approach, even for locally finite graphs. For example, the {\it hyperbolic compactification\/}%
   \COMMENT{}
   of a locally finite hyperbolic graph~$G$, as introduced by Gromov~[\the\refGromovHyperbolicGroups], is sometimes defined in a somewhat roundabout way: in an intermediate step one endows $G$ with a different (bounded) {\it Gromov metric} (which on $G$ itself is somewhat odd%
   \Footnote{For example, it induces the discrete topology on~$G$, including edges!}),%
   \COMMENT{}
   whose sole purpose is that its completion has a compact boundary that can then be combined with $G$ in its {\em original\/} topology to yield an interesting compact space.%
   \Footnote{For example, this {\it hyperbolic boundary\/} refines the `end boundary' $\Omega(G)$ of $G$ induced by~$|G|$, and unlike $\Omega(G)$ it can have non-trivial connected components.}
   But Gromov~[\the\refGromovHyperbolicGroups] also proved that the hyperbolic compactification of $G$ can, alternatively, be obtained as the completion of $(G,d_\ell)$ for a natural function $\ell$ of edge lengths (see [\the\refCoornaertDelzantPapadopoulos] for a detailed proof, or~[\the\refAgelosLTop] for an indication),%
   \COMMENT{}
   a result that makes its use---e.g.\ for the study of hyperbolic groups---more accessible to a graph-theoretical approach. Figure~\figHyperbolic\ shows a hyperbolic graph $G$ with an assignment $\ell$ of edge lengths for which the $\ell$-\Top\ boundary of $(G,d_\ell)$ is its hyperbolic boundary.%
   \COMMENT{}

 \figure \figHyperbolic.
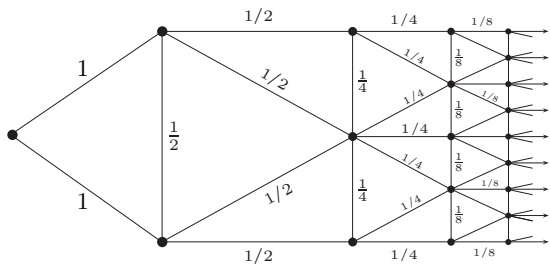
 A hyperbolic graph whose boundary is an~\rlap{arc}\nl {(the vertical bar on the right)}
 (Hyperbolic; 1000)

The flexibility of the metric-completion approach to defining a boundary of an infinite graph $G$ allows us to keep our options open regarding the opening question of Section~\secOutlook.1, the problem of how to define circuits---and the resulting homology---in a graph with infinite degrees. Depending on the graphs under investigation, we can choose our edge length function $\ell$ in whatever way seems best to describe those graphs: we can {\em choose\/} whether or not we want two rays from the same end to converge to a common point at infinity or not, or not even to converge at all.%
   \COMMENT{}

For example, for the fan $G$ shown on the left in Figure~\figFans\ we could let~$\ell$ assign length $2^{-n}$ to the (one or) two edges down or sideways from the $n$th vertex of the vertical ray. Then this ray converges in $(G,d_\ell)$ to the vertex~$u$, so $(G,d_\ell)$ is complete and induces \ITop\ on~$G$. The double fan of Figure~\figFans, with a similar function $\ell$ of edge lengths, will have the two vertices of infinite degree identified in~$(G,d_\ell)$ (which is again already complete), the vertical ray will converge to this identified new vertex, and the two bottom edges will form a circle.\looseness=-1

As another example, consider the powers $T^2$ and $T^3$ of the $\aleph_0$-regular tree~$T$. While $T^2$ still has all the `original' ends of~$T$ (and even some new ones corresponding to the levels of~$T$), the graph $T^3$ does not. Indeed, since every vertex $v$ has infinitely many neighbours at a lower level than its own, or at level~1 (the upper neighbours of the root), no finite set $S$ of vertices can separate $v$ from the infinitely many vertices at level~1, which form a complete subgraph in~$T^3$. Hence $T^3-S$ has only one component, and $T^3$ has only one end! In the context of hamiltonicity problems for graph powers, however, we may well want that $G^3$ retains the end structure of~$G^2$, but that its vertices and boundary points still form a compact set (so as to allow for Hamilton circles; see Section~\secApplications.1).%
   \COMMENT{}
   Both these can be achieved by choosing $\ell$ appropriately and considering the metric completion of $(G^3,d_\ell)$ instead of the Freudenthal compactification~$|G^3|$ of~$G^3$.%
   \COMMENT{}

But this flexibility comes at a cost: in order to benefit from our existing cycle space theory and its applications in the case of locally finite graphs, we need a homology theory that works in the generality of all metric spaces that can arise as such completions, while defaulting to our topological cycle space when the graph considered is locally finite. This homology theory will not just be a minor adaptation of the definition of~$\C(G)$, since the spaces that can arise as completions need not look like graphs.%
   \Footnote{Gromov~[\the\refGromovHyperbolicGroups] observed that every compact metric space arises as the hyperbolic boundary of a hyperbolic graph. Georgakopoulos~[\the\refAgelosLTop] showed that a metric space arises as the $\ell$-\fnTop\ boundary of a locally finite graph if and only if it is complete and has a countable dense subset.}%
   \COMMENT{}
   For example, the boundary of the hyperbolic graph of Figure~\figHyperbolic\ is an arc that contains no edge. In particular, there may be circles in this completion in which the edges they contain are not dense; we shall therefore not be able to represent homology classes by sets of edges, as we did in the case of~$\C(G)$.

\proclaim Problem \xxxMetricHomology.
Extend the topological cycle space theory of locally finite graphs to a homology theory for metric spaces that can arise as completions of spaces of the form~$(G,d_\ell)$.\problem{}

As a step towards this goal, it will help to recast the topological cycle space theory of Sections \secConcepts\ and~\secTheory, even for locally finite graphs, as a homology theory in the usual terms of algebraic topology~[\th\the\refHomologySpaces]. We shall look at this problem next. Some first steps towards Problem~\xxxMetricHomology\ itself are already taken in~[\the\refAgelosLTop].

\subsection \secOutlook.4 Homology of locally compact spaces with ends

In this section we look at $|G|$ from the viewpoint of algebraic topology. This section will be expanded to a more comprehensive survey of these aspects in~[\the\refHomSurvey]. Any undefined graph $G$ will again be an arbitrary connected and locally finite graph.

When the topological cycle space was first introduced~[\the\refCyclesExpository,\th\the\refCyclesOne], the motivation was to extend the standard notion of the cycle space of a finite graph in such a way that it could play a similar role for infinite graphs. In particular, the various theorems relating $\C(G)$ to other structural properties of~$G$, such as planarity, should extend to the new notion. The topological cycle space $\C(G)$, for locally finite~$G$, has achieved this goal superbly---see Section~\secApplications---but at a price: due to the ad-hoc manner of its definition, this new $\C(G)$, unlike~$\C\fin(G)$, is no longer obviously a special case of the first homology of either $G$ or~$|G|$ in any sense that would be standard in algebraic topology.

On the other hand, it seems desirable to recast our notion of $\C(G)$ in terms of standard homology. This is for two reasons: our understanding of $\C(G)$ might benefit from the vast body of knowledge available for standard homology theories; and conversely, the ideas that went into the notion of~$\C(G)$ might throw a new light on the homology of more general spaces once we can rephrase $\C(G)$ in such a way that it makes sense when $G$ is not a graph.

For a finite or infinite graph~$G$, its first simplicial homology group is the same as, in our notation,~$\C\fin(G)$. The elements of this group or vector space are finite sets of edges.%
   \Footnote{We still take all our coefficients from~$\F_2$, for easier compatibility with our earlier treatment of~$\C(G)$. However all the results we describe hold with integer coefficients too, which is indeed the setting in the papers we refer to.}
   Hence as soon as $G$ contains an infinite circuit, $\C(G)$~differs from~$\C\fin(G)$. But $\C(G)$ also differs from the simplicial homology we would obtain from just allowing infinite (thin) sums of edges as 1-chains: in this homology, the edges of a double ray would form an algebraic cycle, since every vertex lies on two such edges and hence the boundary of this edge set is zero, but this is not a cycle in our topological cycle space. Thus, $\C(G)$~cannot be described by an extension of simplicial homology obtained by just allowing infinite chains.%
   \Footnote{By contrast, $\cutspace(G)$, the set of all cuts, remains the image of the coboundary operator~$\delta^0$ if we allow arbitrary infinite 0-cochains.}

Another way to extend the simplicial homology of finite complexes to infinite ones is to take limits. In our situation, we could use the techniques from Section~\secTechniques\ to do this. Indeed, recall that for the finite minors $G_0\minor G_1\minor\dots$ of~$G$ defined in Section~\secTechniques.1 the algebraic cycles of $G_m$ induce algebraic cycles in~$G_n$, for all $n < m$. The first simplicial homology groups of the $G_n$ thus form an inverse system, whose limit is indeed~$\C(G)$:%
   \Footnote{This is tantamount to taking the {\it \vvvvv Cech homology\/} of $G$ or of~$|G|$;%
   \COMMENT{}
   see~[\the\refHomologyGraphs].}

\proclaim Theorem \xxxCasInverseLimit.
The maps $D\mapsto D\cap E(G_n)$ define a group isomorphism from $\C(G)$ to $\invlim\! (\C\fin(G_n))_{n\in\N}\,.$

While Theorem \xxxCasInverseLimit\ points out an interesting aspect of~$\C(G)$, it does not capture all its relevant aspects. For even if we assume that representing $\C(G)$ as $\invlim \C(G_n)$ captures all its relevant features as a group, the (abstract) homology group of a graph does not capture all the relevant aspects of its homology.%
   \Footnote{This is reflected by the fact that, as abstract groups, $\C(G)$~hardly depends on~$G$: by Theorem \xxxFundamental, it is the direct product of as many copies of the coefficient ring as a \TST\ of $G$ has chords---i.e., of $\aleph_0$ copies for most~$G$.}
   For example, recall that if $G$ is 3-connected, then its peripheral circuits generate~$\C(G)$ (Theorem \xxxGenerating\th(ii)), and $G$ is planar if and only if these circuits form a sparse set (Theorem \xxxKelmans). In order to prove this for infinite~$G$, we thus need to know when a given edge set $D\in\C(G)$ is a peripheral circuit. If all we know of $D$ are its intersections with the edge sets of the finite minors~$G_n$ (as we would if we viewed $D$ as an element of~$\invlim\! (\C\fin(G_n))_{n\in\N}$), this would be hard to decide; indeed it would not be obvious from this information whether $D$ is a circuit at all (see Section~\secTechniques.2).%
   \COMMENT{}

A~more subtle approach, which has been pursued in~[\the\refHomologyGraphs,\th\the\refHomologySpaces], is to see to what extent $\C(G)$ can be captured by the singular homology of~$|G|$. After all, $\C(G)$~was defined via (the edge sets of) circles in~$|G|$, which are just injective singular loops. Can we extend this correspondence between injective loops and circuits to one between $H_1(|G|)$ (singular) and~$\C(G)$?

There are two things to notice about~$H_1(|G|)$. The first is that we can subdivide a 1-simplex, or concatenate two 1-simplices into one, by adding a boundary. Indeed, if $\sigma\: [0,1]\to |G|$ is a path in $|G|$ from $x$ to~$y$, say, and $z$ is a point on that path, there are paths $\sigma'$ from $x$ to~$z$ and $\sigma''$ from $z$ to~$y$ such that $\sigma'+\sigma''-\sigma$ is the boundary of a singular 2-simplex `squeezed' on to the image of~$\sigma$. The second fact to notice is that inverse paths cancel in pairs: if $\sigma^+$ is an $x$--$y$ path in~$|G|$, and $\sigma^-$ an $y$--$x$ path with the same image as~$\sigma^+$, then $[\sigma^+ + \sigma^-] = 0\in H_1$.%
   \Footnote{To see that this sum is a boundary, subtract the constant 1-simplex $\sigma$ with value~$x$: there is an obvious singular 2-simplex of which $\sigma^+ + \sigma^- - \sigma$ is the boundary. Subtracting~$\sigma$ is allowed, since $\sigma = \sigma + \sigma - \sigma$, too, is a boundary: of the constant 2-simplex with value~$x$.}
   These two facts together imply that every homology class in $H_1$ is represented by a single loop: given any 1-cycle, we first add pairs of inverse paths between the endpoints of its simplices to make its image connected in the right way,%
   \COMMENT{}
   and then use Euler's theorem to concatenate the 1-simplices of the resulting chain into a single loop~$\sigma$.

To establish the desired correspondence between $H_1(|G|)$ and~$\C(G)$, we would like to assign to a homology class in~$H_1(|G|)$, represented by a single loop~$\sigma$, an edge set $f([\sigma])\in\cyclespace(G)$. Intuitively, we do this by counting for each edge $e$ of $G$ how often $\sigma$ traverses it entirely (which, since the domain of $\sigma$ is compact, is a finite number of times), and let $f([\sigma])$ be the set of those edges $e$ for which this number is odd. Using the usual tools of homology theory, one can make this precise in such a way that $f$ is clearly a well-defined homomorphism $H_1(|G|)\to\edgespace(G)$,%
   \Footnote{For each edge~$e$, let $f_e\: |G|\to S^1$ be a map wrapping $e$ once round~$S^1$ and mapping all of $|G|\sm\interior e$ to one point of~$S^1$. Let $\pi$ denote the group isomorphism $H_1(S^1;\Z_2)\to\Z_2$. Given $h\in H_1(|G|)$, let $f(h) := \{\,e\mid (\pi\circ (f_e)_*)(h) = 1\in\Z_2\,\}$. See~[\the\refHomologyGraphs] for details.}
   and whose image is easily seen to be~$\C(G)$. What is not clear at once is whether $f$ is 1--1 and onto.

Surprisingly, $f$~is indeed surjective---and this is not even hard to show. Indeed, let an edge set $D\in\C(G)$ be given. Our task is to find a loop~$\sigma$ that traverses every edge in~$D$ an odd number of times, and every other edge of $G$ an even number of times. As a first approximation, we let $\sigma_0$ be a path that traverses every edge of some fixed \NST\ of $G$ exactly twice, once in each direction; see Section~\secTechniques.3 for how to construct such a loop. Moreover, we construct $\sigma_0$ in such a way that it pauses at every vertex~$v$---more precisely, so that $\sigma_0^{-1}(v)$ is a union of finitely many closed intervals at least one of which is non-trivial. Next, we write $D$ as a thin sum $D = \sum_i C_i$ of circuits; such a representation of $D$ exists by definition of~$\C(G)$. For each of these $C_i$ we pick a vertex $v_i\in\overline{C_i}$, noting that no vertex of $G$ gets picked more than finitely often, because it has only finitely many incident edges and the $C_i$ form a thin family. Finally, we turn $\sigma_0$ into the desired loop~$\sigma$ by expanding the pause at each vertex $v$ to a loop going once round every $\overline{C_i}$ with $v=v_i$. Using the methods from Section~\secTechniques.3 it is not hard to show that $\sigma$ is continuous~[\the\refHomologyGraphs], and clearly it traverses every edge of $G$ the desired number of times.

Equally surprisingly, perhaps, $f$~is usually not injective (see below). In summary, therefore, the topological cycle space $\C(G)$ of $G$ is related to the first singular homology group of~$G$ as follows:

\proclaim Theorem \xxxSingular.~[\the\refHomologyGraphs]
The map $f\colon H_1(|G|)\to\edgespace(G)$ is a group homomorphism onto~$\C(G)$, which has a non-trivial kernel if and only if $G$ contains infinitely many (finite) circuits.

An example of a non-null-homologous loop in~$|G|$ whose homology class maps to the empty set $\es\in\C(G)$ is easy to describe. Let $G$ be the one-way infinite ladder~$L$ (with its end on the right), and define a loop~$\rho$ in~$L$, as follows. We start at time~$0$ at the top-left vertex, $v_0$~say, and begin by going round the first square of $L$ in a clockwise direction. This takes us back to~$v_0$. We then move along the horizontal edge incident with~$v_0$, to its right neighbour~$v_1$. From here, we go round the second square in a clockwise direction, back to~$v_1$ and on to its right neighbour~$v_2$. We repeat this move until we reach the end~$\omega$ of~$L$ on the right, say at time ${1\over2}\in [0,1]$. So far, we have traversed the first vertical edge and every bottom horizontal edge once (in the direction towards~$v_0$), every other vertical edge twice (once in each direction), and every top horizontal edge twice in the direction towards the end. From there, we now use the remaining half of our time to go round the infinite circle formed by the first vertical edge and all the horizontal edges one and a half times, in such a way that we end at time~1 back at~$v_0$ and have traversed every edge of~$L$ equally often in each direction. Clearly, $f$~maps (the homology class of) this loop~$\rho$ to $0\in\C(G)$.

 \figure \figKringelsimplex.
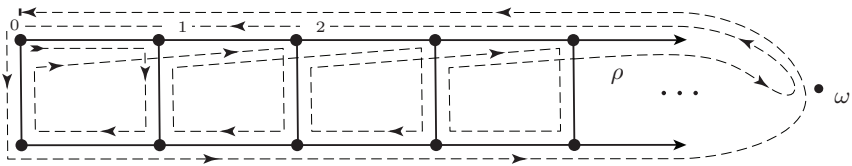
 The loop $\rho$ is not null-homologous
 (Kringelsimplex; 1000)

The loop $\rho$ is indeed not null-homologous~[\the\refHomologyGraphs], but it seems non-trivial to show this. To see why this is hard, let us compare $\rho$ to a loop winding round a finite ladder in a similar fashion, traversing every edge once in each direction. Such a loop~$\sigma$ is still not null-homotopic, but it is null-homologous. To see this, we subdivide it into single edges: we find a finite collection of 1-simplices~$\sigma_i$, two for every edge, such that $[\sigma] = \big[\sum_i\sigma_i\big]$ and every $\sigma_i$ just traverses its edge. Next, we pair up these~$\sigma_i$ into cancelling pairs: if $\sigma_i$ and~$\sigma_j$ traverse the same edge~$e$ (in opposite directions), then $[\sigma_i + \sigma_j] = 0$. Hence $[\sigma] = \big[\sum_i\sigma_i\big] = 0$, as claimed. But we cannot imitate this proof for $\rho$ and our infinite ladder~$L$, because homology classes in $H_1(|G|)$ are still finite chains: we cannot add infinitely many boundaries to subdivide $\rho$ infinitely often.

As it happened, the proof of the seemingly simple fact that $\rho$ is not null-homologous took a detour via the solution of a much more fundamental problem: the problem of understanding the fundamental group of~$|L|$, or more generally, of~$|G|$ for a locally finite graph~$G$. In order to distinguish $\rho$ from boundaries, we looked for a numerical invariant $\Lambda$ of 1-chains%
   \COMMENT{}
   that was non-zero on~$\rho$ but both linear and additive (so that $\Lambda(\sigma_1\sigma_2) = \Lambda(\sigma_1 + \sigma_2) = \Lambda(\sigma_1) + \Lambda(\sigma_2)$ for concatenations of 1-simplices~$\sigma_1,\sigma_2$) and invariant under homotopies (so that $\Lambda(\sigma_1\sigma_2) = \Lambda(\sigma)$ when $\sigma\sim \sigma_1\sigma_2$). Then, given a 2-simplex $\tau$ with boundary $\partial\tau = \sigma_1 + \sigma_2 - \sigma$, we would have $\Lambda(\partial\tau) = \Lambda(\sigma_1\sigma_2) - \Lambda(\sigma) = 0$, so~$\Lambda$ would vanish on all boundaries%
   \COMMENT{}
   but not on~$\rho$. We did not quite find such an invariant~$\Lambda$, but a collection of similar invariants which, together, can distinguish loops like~$\rho$ from boundaries.

In order to find such functions on 1-chains that are invariant under homotopies, it was necessary to find a combinatorial description for the homotopy types of loops---that is, for the fundamental group of~$|G|$. Such a combinatorial characterization is given in~[\the\refHomotopy], where $\pi_1(|G|)$ is characterized as a group of (infinite) words of chords of a \TST---in the spirit of the usual description of~$\pi_1(G)$ for a finite graph as the free group generated by such chords---and as a subgroup of the inverse limit of these finitely generated free groups. (The group $\pi_1(|G|)$ itself is not free, unless $G$ is essentially finite.)

Not surprisingly, our topological cycle space $\C(G)$ can, as a group, be viewed as the {\it infinite abelianization\/} of $\pi_1(|G|)$: the factor group of $\pi_1(|G|)$ obtained by declaring two (reduced) words as equivalent if each letter occurs in both the same number of times (regardless of position).

Let us return to our original goal: to see to what extent $\C(G)$ can be captured by the singular homology of~$|G|$. In view of Theorem~\xxxSingular, the goal might be phrased as follows:

\proclaim Problem \xxxSingularAlt. Devise a singular-type homology theory for locally compact spaces with ends that coincides with $\C(G)$ when applied to $|G|$ in dimension~1.

Some first steps in this direction were already taken in~[\the\refHomologyGraphs]. The approach there was to allow only singular simplices whose vertices lie in~$G$ (i.e., are not ends); to allow infinite chains that are locally finite at points of~$G$---i.e., every point in $G$ has a neighbourhood meeting only finitely many simplices in the given chain---but not necessarily at ends (reflecting the fact that, for example, the union of all the vertical double rays in the grid defines an algebraic cycle of which infinitely many 1-simplices contain the end); and to restrict the set of cycles to those chains in the kernel of the boundary operator that could be written as a (possibly infinite) sum of finite cycles.

This approach solves Problem~\xxxSingularAlt\ in an ad-hoc sort of way: it permits the definition of homology groups for arbitrary locally compact spaces with ends, it defaults to the topological cycle space for graphs in dimension~1---but it is not a homology theory in the sense of the usual axioms~[\the\refHatcher,\th\the\refEilenbergSteenrod]. To achieve the latter, one has to find a way of implementing the required restrictions as conditions on chains rather than on cycles. This was done in~[\the\refHomologySpaces]. However, this is no more than a beginning, and more translation work remains to be done---for example, of the duality%
   \COMMENT{}
   theory indicated for $\C(G)$ in Section~\secTheory.

\subsection \secOutlook.5 Infinite matroids

Traditionally, infinite matroids are defined like finite ones, with the additional axiom that an infinite set is independent as soon as all its finite subsets are independent. This reflects the notion of linear independence in vector spaces, and also the absence of the usual (finite) cycles in a graph: the bases of the cycle matroid of an infinite graph are then the edge sets of its (ordinary) spanning trees. We shall call such matroids {\em finitary\/}. Note that the circuits in a finitary matroid, the minimal dependent sets, are necessarily finite.%
   \COMMENT{}

An important and regrettable feature of such finitary matroids is that the additional axiom restricting the infinite independent sets spoils duality, one of the key features of matroid theory. For example, every bond of a graph would be a circuit in any dual of its cycle matroid: a set of edges that is minimal with the property of not lying in the complement of a spanning tree, i.e.\ of containing an edge from every spanning tree.%
   \COMMENT{}
   Since finitary matroids have no infinite circuits, the cycle matroid of a graph with an infinite bond thus cannot have a finitary dual.

Our theory, however, suggests an obvious solution to this problem: shouldn't infinite matroids be defined in such a way that infinite circuits in a graph can become matroid circuits, and \TST s become bases? Indeed, infinite circuits are not contained in the edge set of any \TST\ (although they {\em are\/} contained in the edge set of an ordinary spanning tree), while if we delete any edge from an infinite circuit, then its remaining edges can be extended to a \TST\ by Lemma~\xxxTSTexistence.

There are two main challenges in devising axioms for such a non-finitary theory of infinite matroids: to avoid the mention of cardinalities, and to take care of limits.%
   \Footnote{If one wants to have basis and circuit axioms, one has to ensure that maximal independent sets and minimal dependent sets exist: with infinite sets, this is no longer clear.}
   In~[\the\refInfiniteMatroidAxioms] such a theory has been proposed. It can be stated in terms of any of five equivalent sets of axioms: independence, basis, circuit, closure or rank axioms. They are shown to be equivalent to the `B-matroids' explored in the late 1960s by Higgs~[\the\refHiggsMatroidsDuality,\th\the\refHiggsBMatroids,\th\the\refHiggsMatroids], who had defined them in terms of a different and rather more complicated set of closure requirements. As just one of a plethora of alternatives for a possible concept of infinite matroids considered at the time, these `B-matroids' had gone largely unnoticed, although a workable combination of independence and exchange axioms was later found by Oxley~[\the\refOxleyInfiniteMatroidsPaper,\th \the\refOxleyInfiniteMatroidsSurvey].

With any of the said five sets of axioms, duality works as expected from finite matroids: it is exemplified by dual planar graphs (see below), there is a well-defined notion of minors with contraction and deletion as dual operations, and so on.

So what non-finitary matroids are there in graphs? As had been our motivation, the circuits (finite or infinite) of a locally finite graph~$G$ form the circuits of a matroid in this theory, as do the finite circuits of~$G$. Let us denote these matroids by $M_{\rm C}(G)$ and~$M_{\rm FC}(G)$, respectively. Similarly, the bonds (finite or infinite) of $G$ form the circuits of a matroid~$M_{\rm B}(G)$, just as its finite bonds form the circuits of a matroid~$M_{\rm FB}(G)$. Clearly, $M_{\rm FC}(G)$ and $M_{\rm B}(G)$ form a pair of dual matroids,%
   \COMMENT{}
   and by Lemmas \xxxConnected\ and~\xxxTSTeq\ so do $M_{\rm C}(G)$ and~$M_{\rm FB}(G)$. 

The same is true for the slightly more general identification spaces $\tilde G$ under \ITop, defined in Section~\secOutlook.2. Thus if $G$ and $G^*$ are a pair of dual finitely separable graphs, and we take as circuits the edge sets of circles in the spaces $\tilde G$ and~$\tilde{G^*}$, then $M_{\rm FC}(G)$ and $M_{\rm C}(G^*)$ form a dual pair of matroids, as do $M_{\rm FB}(G)$ and~$M_{\rm B}(G^*)$.

Call a (finite or infinite) matroid {\it graphic\/} if it is the cycle matroid $M_{\rm C}(G)$ of a graph~$G$, and {\it finitely graphic\/} if it is the finite-cycle matroid $M_{\rm FC}(G)$ of a graph~$G$. The infinite version of Whitney's theorem (Theorem~\xxxWhitney) can now be restated in matroid terms (see~[\the\refInfiniteMatroidAxioms]):

\proclaim Theorem \xxxMatroidWhitney.
A finitely separable graph $G$ is planar if and only if its cycle matroid has a finitely graphic dual.

\noindent
The conditions in Theorem \xxxMatroidWhitney\ should be equivalent also to saying that the finite-cycle matroid of $G$ has a graphic dual, but this has no been proved.

\medbreak

Much of the attractiveness of finite matroids stems from the fact that they provide a unified framework for some essential common aspects of otherwise disparate branches of mathematics. Whether or not the same can be said for infinite matroids, axiomatized in this way, will depend on concrete examples that have yet to be found---if possible, from as different areas of mathematics as possible.

There is no doubt that non-finitary matroids are plentiful. Indeed, a~finitary matroid has a finitary dual only if it is the direct sum of finite matroids~[\the\refLasVergnasMatroids]. Since all our matroids have duals, the duals of all the other finitary matroids (e.g., of all connected matroids~[\the\refBeanFinitary]) thus form a large class of non-finitary matroids. However, perhaps there are natural `primary' matroids that are non-finitary and have therefore gone unnoticed---for example, in the context of Banach or Hilbert spaces?

Here is an interesting concrete problem. An example in~[\the\refInfiniteMatroidAxioms] shows that there are matroids with both infinite circuits and cocircuits, indeed matroids in which all these are infinite. However, we do not know the answer to the following:

\proclaim Problem \xxxMatroidExamples. Is the intersection of a circuit and a cocircuit always finite?\problem{}

Theorem~\xxxDualTopTrees\ seems to suggest that, given a locally finite graph~$G$ (or a finitely separable one), there might be a matroid on the set $E(G)\cup\Omega(G)$ in which the sets $F\cup\Psi$ with $F\sub E(G)$ and $\Psi\sub\Omega(G)$ for which $\bigcup F\cup\Psi$ contains no circle form the independent sets. However, as soon as we try to apply basis or circuit elimination axioms to ends, we see that this fails. For example, if two circles meet in exactly one end, we cannot delete the end and find another circle in the rest (the union of all the edges of the two circles). Similarly, in the double ladder we could choose as a basis $B_1$ the union of one double ray, all the rungs, and both ends, and as another basis $B_2$ all the edges that are not rungs and one of the two ends. If we delete the other end, $\omega$~say, from~$B_1$, we cannot find an element of $B_2\sm B_1$ (which would be an edge) that we could add to $B_1 - \omega$ to form another basis.

More generally, the duality which Theorem~\xxxDualTopTrees\ expresses for graphs with ends cannot be expressed in terms of matroid duality. Indeed, suppose that, with the notation of Theorem~\xxxDualTopTrees, there is a matroid $M$ on $E\cup\Omega$ whose bases are the sets $F\cup\Psi$ (with $F\sub E$ and $\Psi\sub\Omega$) that form the subspaces $X$ of~$|G|$ for which $\tilde X$ is a spanning tree of~$\tilde G$.%
   \COMMENT{}
   Then $M/\Omega = M_{\rm C}(G)$ and $M\bs\Omega = M_{\rm FC}(G)$. But one can show that, given any 2-connected finitely separable graph~$G$, there is no matroid $M$ on any set $E(G)\cup X$ with $X\cap E(G) = \es$ such that $M/X = M_{\rm C}(G)$ and $M\bs X = M_{\rm FC}(G)$. See~[\the\refTreeEnds] for details.

Richter et al.\ used a non-finitary matroid similar to $M_{\rm C}(G)$ for their proof of an extension of Theorem~\xxxWhitney\ (Whitney) to `graph-like spaces'; see [\the\refRichterGraphLikeSpacesSurvey].



{\eightpoint
\beginsection References

\ref\refAharoniCtbleEM
   R.\th Aharoni, Menger's theorem for countable graphs, \JCTB43 (1987), 303--313.

\ref\refAharoniBergerEM
   R.\th Aharoni \& E.\th Berger, Menger's theorem for infinite graphs, \Inv176 (2009), 1--62.

\ref\refMFMC
   R.\th Aharoni, E.\th Berger,  A.\th Georgakopoulos, A.\th Perlstein~\& P.\th Spr\"ussel, The max-flow min-cut theorem for countable networks, \JCTB101 (2011), 1--17.

\ref\refAharoniCT
   R.\th Aharoni \& C.\th Thomassen, Infinite highly connected digraphs with no two arc-disjoint spanning trees, \JGT13 (1989)

\ref\refABLLeftRightTours
   D.\th Archdeacon, P.\th Bonnington, and C.\th Little, An algebraic characterization
of planar graphs, \JGT19 (1995), 237--250.

\ref\refAsratianKhachatrianLocalization
   A.S.\th Asratian \& N.K.\th Khachatrian, Some localization theorems on hamiltonian circuits, \JCTB49 (1990), 287--294.

\ref\refBeanFinitary
   D.W.T.\th Bean, A connected finitary co-finitary matroid is finite, In {\em Proceedings of the Seventh Southeastern Conference on Combinatorics, Graph Theory and Computing\/}, {\em Congressus Numerantium~\bf 17} (1976), 115--119.

\ref\refBenjaminiSchramm
   I.\th Benjamini \& O.\th Schramm, Harmonic functions on planar and almost
planar graphs and manifolds, via circle packings, \Inv126 (1996), 565--587.

\ref\refBergerBruhnEndDegrees
   E.\th Berger \& H.\th Bruhn, Eulerian edge sets in locally finite graphs, \Comb31 (2011), 21--38.

\ref\refBiggsPotential
   N.\th Biggs, Algebraic potential theory on graphs, \BLMS29 (1997), 641--682.

\ref\refEGT
   B.\th Bollob\'as, {\sl Extremal Graph Theory}, Academic Press, London 1978.

\ref\refBruhnPeripheral
   H.\th Bruhn, The cycle space of a 3-connected locally finite graph is generated by its finite and infinite peripheral circuits, \JCTB92 (2004), 235--256.

\ref\refBruhnPersonal
   H.\th Bruhn, personal communication 2009.

\ref\refDuality
   H.{\tie}Bruhn and R.{\tie}Diestel, Duality in infinite graphs, \CPC15 (2006), 75--90.

\ref\refMacLaneArbitrarySurfaces
   H.{\tie}Bruhn and R.{\tie}Diestel, MacLane's theorem for arbitrary surfaces, \JCTB99 (2009), 275--286.

\ref\refInfiniteMatroidAxioms
   H.{\tie}Bruhn, R.{\tie}Diestel, M.\th Kriesell \& P.\th Wollan, Axioms for infinite matroids, preprint arXiv:1003.3919 (2010).

\ref\refTreeEnds
   H.\th Bruhn, R.\th Diestel, and J.\ Pott, Dual trees must share their ends, preprint 2011.

\ref\refCyCoCy
   H.{\tie}Bruhn, R.{\tie}Diestel~\& M.{\tie}Stein, Cycle-cocycle partitions and faithful cycle covers for locally finite graphs, \JGT50 (2005), 150--161.

\ref\refBruhnSteinDiestelEM
   H.{\tie}Bruhn, R.{\tie}Diestel~\& M.{\tie}Stein, Menger's theorem for infinite graphs with ends, \JGT50 (2005), 199--211.

\ref\refAgelosHenningLA
   H.\th Bruhn \& A.\th Georgakopoulos, Bases and closures under infinite sums, preprint 2008.

\ref\refBruhnBicycles
   H.\th Bruhn, S.\th Kosuch~\& M.\th Win Myint, Bicycles and left-right tours in locally finite graphs, \EJC30 (2009), 356--371.

\ref\refBruhnSteinMacLane
   H.{\tie}Bruhn and M.{\tie}Stein, MacLane's planarity criterion for locally finite graphs, \JCTB96 (2006), 225--239.

\ref\refBruhnSteinEndDeg
   H.{\tie}Bruhn and M.{\tie}Stein, On end degrees and infinite circuits in locally finite graphs, \Comb27 (2007), 269--291.

\ref\refBruhnSteinEndDuality
   H.{\tie}Bruhn and M.{\tie}Stein, Duality of ends, \CPC12 (2009), 47--60.

\ref\refBruhnYuHamilton
   H.{\tie}Bruhn and X.\th Yu, Hamilton circles in planar locally finite graphs, \SIAM22 (2008), 1381--1392.

\ref\refCasteelsRichterBicycles
   K.\th Casteels and B.\th Richter, The bond and cycle spaces of an infinite
graph, \JGT59 (2008), 126--176.

\ref\refCatlinLineGraphs
   P.A.\th Catlin, Supereulerian graphs: a survey, \JGT16 (1992), 177--196.

\ref\refCoornaertDelzantPapadopoulos
   M.\th Coornaert, T.\th Delzant \& A.\th Papadopoulos, Geometrie et theorie des groupes. Les groupes hyperboliques de Gromov. {\em Springer Lecture Notes in Mathematics~\bf 1441}, Springer 1990.

\ref\refYuHamilton
   Q.\th Cui, J.\th Wang and X.\th Yu, Hamilton circles in infinite planar graphs, \JCTB99 (2009), 110--138.

\ref\refCtbleEM
   R.\th Diestel, The countable \Erdos-Menger conjecture with ends, \JCTB87 (2003), 145--161.

\ref\refCyclesExpository
   R.\th Diestel, The cycle space of an infinite graph, \CPC14 (2005), 59--79.

\ref\refSpanningTrees
   R.\th Diestel, End spaces and spanning trees, \JCTB96 (2006), 846--854.


\ref\refTopSurveyI
   R.\th Diestel, Locally finite graphs with ends: a topological approach. I.~Basic theory. \DM311 (special volume on infinite graph theory, 2011), 1423--1447.

\ref\refTopSurveyII
   R.\th Diestel, Locally finite graphs with ends: a topological approach. II.~Applications. \DM311 (Carsten Thomassen 60 special volume, 2010), 2750--2765.

\font\ttt=cmtt8
\ref\refBook
     R.\th Diestel, {\it Graph theory\/}, 4th edition, Springer-Verlag
     2010. Electronic edition available at {\ttt http://diestel-graph-theory.com/}

\ref\refMinorUniversal
   R.\th Diestel \& D.\th K\"uhn, A~universal planar graph under the minor relation, \JGT32 (1999), 191--206.

\ref\refCyclesOne
   R.\th Diestel \& D.\th K\"uhn, On infinite cycles~I, \Comb24 (2004), 69--89.

\ref\refCyclesTwo
   R.\th Diestel \& D.\th K\"uhn, On infinite cycles~II, \Comb24 (2004), 91--116.

\ref\refTST
   R.\th Diestel \& D.\th K\"uhn, Topological paths, cycles and
spanning trees in infinite graphs, \EJC25 (2004), 835--862.

\ref\refNST
   R.\th Diestel \& I.B.\th Leader, Normal spanning trees, Aronszajn trees and excluded minors, \JLMS63 (2001), 16--32.

\ref\refHomologyGraphs
   R.\th Diestel \& P.\th Spr\"ussel, The homology of a locally finite graph with ends, \Comb30 (2010), 681--714 (mit P.\th Spr\"ussel).

\ref\refHomotopy
   R.\th Diestel \& P.\th Spr\"ussel, The fundamental group of a locally finite graph with ends, \Advances226 (2011) 2643--2675.

\ref\refHomologySpaces
   R.\th Diestel \& P.\th Spr\"ussel, On the homology of locally compact spaces with ends, {\sl Topology and its Applications\/}~{\bf 158} (2011), 1626--1639 (mit Ph.\th Spr\"ussel).

\ref\refHomSurvey
   R.\th Diestel \& P.\th Spr\"ussel, Locally finite graphs with ends: a topological approach. III.~Fundamental group and homology, \DM312 (special volume on algebraic graph theory, 2011), 21--29 (mit Ph.\th Spr\"ussel).

\ref\refEilenbergSteenrod
   S.\th Eilenberg and N.\th Steenrod, {\it Foundations of Algebraic Topology\/}, Princeton University Press 1952.

\ref\refAgelosPathConnected
A.\th Georgakopoulos, Connected but not path-connected subspaces of infinite graphs, \Comb27 (2007), 683--698.

\ref\refAgelosOWReportFleischner
A.\th Georgakopoulos, Fleischner's theorem for infinite graphs, {\sl Oberwolfach reports \bf 4} (2007).

\ref\refAgelosHotchpotch
   A.\th Georgakopoulos, Topological circles and Euler tours in locally finite graphs, \EJ16 (2009), \#R40.

\ref\refAgelosInfiniteFleischner
   A.\th Georgakopoulos, Infinite Hamilton cycles in squares of locally finite graphs, \Advances220 (2009), 670--705.

\ref\refAgelosLTop
   A.\th Georgakopoulos, Graph topologies induced by edge lengths, preprint arXiv:0903.1744 (2009).

\ref\refAgelosUniqueFlows
   A.\th Georgakopoulos, Uniqueness of electrical currents in a network of finite total resistance, preprint arXiv:0906.4080 (2009).

\ref\refAgelosPersonal
   A.\th Georgakopoulos, personal communication 2009.

\ref\refAgelosPhilippGeodetic
   A.\th Georgakopoulos \& Ph.\th Spr\"ussel, Geodetic topological cycles in locally finite graphs, \EJ16 (2009), \#R144.

\ref\refGromovHyperbolicGroups
   M.\th Gromov, Hyperbolic Groups, in: {\it Essays in group theory\/} (S.M.\th Gersten, ed), MSRI series vol.~8, pp.~75--263, Springer, New York, 1987.

\ref\refHahnEdgeEnds
   G.\th Hahn, F.\th Laviolette and J.\th Siran, Edge-Ends in countable graphs, \JCTB70 (1997), 225--244.

\ref\refHalinInfGrid
   R.\th Halin, \"Uber die Maximalzahl frem\-der unendlicher Wege, \MN30 (1965), 63--85.

\ref\refHalinMinVertex
   R.\th Halin, A theorem on $n$-connected graphs, \JCTB7 (1969), 150--154.

\ref\refHalinMinimization
   R.\th Halin, Unendliche minimale $n$-fach \zh e Graphen, \Abh36 (1971), 75--88.

\ref\refHalinInfMinimization
   R.\th Halin, Minimization problems for infinite $n$-connected graphs, \CPC2 (1993), 417--436.

\ref\refHalinProblems
   R.\th Halin, Miscellaneous problems on infinite graphs, \JGT35 (2000), 128--151.

\ref\refHatcher
   A.\th Hatcher, {\it Algebraic Topology\/}, Cambridge University Press 2002.

\ref\refHiggsMatroidsDuality
   D.A.\th Higgs, Matroids and duality, {\em Colloq.\ Math.~\bf 20} (1969), 215--220.

\ref\refHiggsBMatroids
   D.A.\th Higgs, Equicardinality of bases in $B$-matroids, {\em Can.\ Math.\ Bull.~\bf 12} (1969), 861--862.

\ref\refHiggsMatroids
   D.A.\th Higgs, Infinite graphs and matroids, {\it Recent Prog.\ Comb., Proc.\ 3rd Waterloo Conf.} (1969), 245--253.

\ref\refLick
   D.\th Lick, Critically and minimally $n$-connected graphs, in (G.\th Chartrand \& S.F.\th Kapoor, eds.): The many facets of graph theory, {\sl Lecture Notes in Mathematics\/~\bf 110} (Springer-Verlag 1969), 199--205.

\ref\refMaderHomEigenschaften
   W.\th Mader, Homomorphieeigenschaften und mittlere Kantendichte von Graphen, \MA174 (1967), 265--268.

\ref\refMaderMinVertex
   W.\th Mader, Minimale $n$-fach zusammenh\"angende Graphen, \MA191 (1971), 21--28.

\ref\refMaderMinVertices
   W.\th Mader, Ecken vom Grad $n$ in minimalen $n$-fach zusammenh\"angenden Graphen, \Archiv23 (1972), 219--224.

\ref\refMaderMinVerticesInf
   W.\th Mader, \"Uber minimal $n$-fach zusammenh\"angende, unendliche Graphen und ein Extremalproblem, \Archiv23 (1972), 553--560.

\ref\refMaderReduktion
   W.\th Mader, Eine Reduktionsmethode f\"ur den Kantenzusammenhang in Graphen, \MN93 (1979), 187--204.

\ref\refMaderEdgeConPres
   W.\th Mader, Paths in graphs reducing the edge-connectivity only by two, \GC1 (1985), 81--89.

\ref\refTheoFlows
   Th.\th M\"uller, personal communication 2007.

\ref\refNadlerContinuumTheory
   B.\th Nadler, {\it Continuum theory\/}, Dekker 1992.

\ref\refNWTreePacking
   C.St.J.A.\th Nash-Williams, Edge-disjoint spanning trees of
finite graphs, \JLMS36 (1961), 445--450.

\ref\refNWArboricity
   C.St.J.A.\th Nash-Williams, Decompositions of finite graphs into forests, \JLMS39 (1964), 12.

\ref\refOberlySumner
   D.J.\th Oberly and D.P.\th Sumner. Every connected, locally connected nontrivial
graph with no induced claw is hamiltonian, \JGT3 (1979), 351--356.

\ref\refOxleyInfiniteMatroidsPaper
   J.G.\th Oxley, Infinite matroids, \PLMS37 (1978), 259--272.

\ref\refOxleyInfiniteMatroidsSurvey
   J.G.\th Oxley, Infinite matroids, in: (N.\th White, ed) {\it Matroid Applications\/}, Cambridge Univ.\ Press 1992.

\ref\refRichterGraphLikeSpacesSurvey
   B.\th Richter, Graph-like spaces: an introduction, preprint 2009.

\ref\refRichterVella
   B.\th Richter and A.\th Vella, Cycle spaces in topological spaces, preprint 2006.

\ref\refSchulzEdgeEnds
   M.\th Schulz, Der Zyklenraum nicht lokal-endlicher Graphen, {\it Diplomarbeit\/}, Univ.\ Hamburg 2005.

\ref\refPhilippNormal
   P.\th Spr\"ussel, End spaces of graphs are normal, \JCTB98 (2008), 798--804.

\ref\refMayaTreePacking
   M.\th Stein, Arboricity and tree-packing in locally finite graphs, \JCTB96 (2006), p.302--312.

\ref\refMayaEndDeg
   M.\th Stein, Forcing highly connected subgraphs in locally finite graphs, \JGT54 (2007), 331--349.

\ref\refMayaBanff
   M.\th Stein, Extremal infinite graph theory, preprint 2009.

\ref\refCTFleischner
  C.{\tie}Thomassen, Hamiltonian paths in squares of infinite locally finite blocks,
  {\sl Ann.\ Discrete Math.~\bf 3} (1978), 269--277.

\ref\refCTPlanDualInf
   C.{\tie}Thomassen, Planarity and duality of finite and infinite graphs, \JCTB29 (1980), 244--271.

\ref\refCTConPres
   C.{\tie}Thomassen, Nonseparating cycles in $k$-connected graphs, \JGT5 (1981), 351--354.

\ref\refCTLineGraphs
   C.\th Thomassen, Reflections on graph theory, \JGT10 (1986), 309--324.

\ref\refCTVella
   C.{\tie}Thomassen and A.{\tie}Vella, Graph-like continua and Menger's theorem, \Comb28 (2009).\?{Add details once known. Springer released the e-version in 2008.}

\ref\refTutteTreePacking
   W.T.\th Tutte, On the problem of decomposing a graph into $n$ connected factors, \JLMS36 (1961), 221--230.

\ref\refVellaThesis
   A.\th Vella, A fundamentally topological perspective on graph theory, PhD thesis, Waterloo~2004.

\ref\refLasVergnasMatroids
   M.~Las Vergnas, Sur la dualit\'e en th\'eorie des matro\"{\i}des, In {\em Th\'eorie des Matro\"{\i}des\/}, \SLNM211 (1971), 67--85.

\ref\refWagnerBook
   K.~Wagner, {\it Graphentheorie\/}, BI-Hochschultaschenb\"ucher, Biblio\-gra\-phi\-sches Institut, Mannheim 1970.

\ref\refWhyburnMenger
   G.T.\th Whyburn, On $n$-arc connectedness, \TAMS63 (1948), 452--456.

\ref\refWoessDirichlet
   W.\th Woess, Dirichlet problem at infinity for harmonic functions on graphs, in (J.~Kral et.\ al, eds): International conference on potential theory (1994), Proceedings, de Gruyter (1996), 189--217. 

\ref\refWoessBook
   W.\th Woess, {\it Random walks on infinite graphs and groups\/}, Cambridge University Press 2000.

\ref\refZhanLineGraphs
   S.\th Zhan, On hamiltonian line graphs and connectivity, \DM89 (1991), 89--95.

}

\bigskip\ninepoint\obeylines\parindent=0pt
Mathematisches Seminar\hfill Version 7.7.2012
Universit\"at Hamburg
Bundesstra\ss e 55
D - 20146 Hamburg
Germany

\bye